\documentclass[11pt]{article}
\usepackage{amsmath}
\usepackage{amssymb}

\newtheorem{thm}{Theorem}[section]

\newtheorem{lem}{Lemma}[section]

\makeatletter
   
          \@addtoreset{equation}{section}
\makeatother

\title{Generalized eigenfunctions of relativistic
Schr\"odinger operators I}
\author{Tomio Umeda\thanks{Research supported by 
     Grant-in-Aid for Scientific Research (C) 
    No. 09640212, 
   Japan Society for the Promotion of Science.}\\
       Department of Mathematics  \\
       Himeji Institute of Technology, Japan\\
  {\it E-mail address}: umeda@sci.himeji-tech.ac.jp }

\date{}

\begin{document}

\maketitle

\begin{abstract}
Generalized eigenfunctions of the
3-dimensional relativistic Schr\"o\-dinger
operator $\sqrt{\Delta} + V(x)$ with
$|V(x)|\le C \langle x \rangle^{{-\sigma}}$,
$\sigma > 1$, are considered. 
We construct the generalized eigenfunctions
by exploiting results on the 
limiting absorption principle.
We compute explicitly the 
integral kernel of
$(\sqrt{-\Delta}-z)^{-1}$,
$z \in {\mathbb C}\setminus [0, \, +\infty)$,
which has nothing in common with
 the integral kernel of $({-\Delta}-z)^{-1}$,
but the leading term of the integral
kernels of the boundary
values $(\sqrt{-\Delta}-\lambda \mp i0)^{-1}$,
$\lambda >0$, turn out to be 
the same, up to a constant, as the integral
kernels of the boundary values 
$({-\Delta}-\lambda \mp i0)^{-1}$.
This fact enables us to show that
the asymptotic behavior, 
as $|x| \to +\infty$, of 
the generalized eigenfunction
of $\sqrt{\Delta} + V(x)$ is equal to
the sum of a plane wave and a spherical wave
when $\sigma >3$.
\end{abstract}

\vspace{20pt}

\tableofcontents

\vspace{10pt}

\section{Introduction}
This is 
the first part of a paper,
consisting of two parts, on the operator
\begin{eqnarray}
\sqrt{-\Delta} + V(x), \quad x \in 
{\mathbb R}^3,
\label{eqn:00}
\end{eqnarray}
with a short range potential $V(x)$,
the operator which we shall call the 
relativistic
Schr\"odin\-ger operator.
The first part, the present paper, 
is concerned with asymptotic behaviors,
as $|x| \to + \infty$, of
the generalized eigenfunctions
of $\sqrt{-\Delta} + V(x)$, 
whereas the second part \cite{Umeda4} 
deals with the completeness of 
the generalized eigenfunctions, i.e.,
the eigenfunction expansion for 
the absolutely continuous spectrum. 
Part of the present and 
 coming papers was announced in \cite{Umeda3}.

We remark here that
a prototype of generalized eigenfunction
expansions is provided by the
Fourier inversion formula
\begin{eqnarray*}
u(x) = (2\pi)^{-n/2}
\int_{{\mathbb R}^n} e^{ix\cdot k}
 \, \hat u(k) \, dk,
\end{eqnarray*}
where $e^{ix \cdot k}$ should be 
regarded as 
a generalized eigenfunction of the
Laplace operator $-\Delta_x$ in the sense
that $e^{ix\cdot k}$ satisfies
$-\Delta_x e^{ix\cdot k} = |k|^2 e^{ix\cdot k}$,
but does not belong to $L^2({\mathbb R}^n_x)$.
It has to be noted that the absolutely
continuous spectrum of $-\Delta$ is 
given by the interval $[0, +\infty)$. 

Although
 relativistic Schr\"odinger operators 
have received a substantial amount of
attention in recent years, 
 there have been only 
a few works on the decay of eigenfunctions 
associated to the discrete spectra of 
these operators; see Nardini\cite{Nardini},
 \cite{Nardini2}
 Carmona-Masters-Simon\cite{CarmonaMS}
and Helffer-Parisse\cite{HelfferPar}).
And it is a surprise that 
up to now
there seems to have been
no results on asymptotic behaviors of the
generalized eigenfunctions of
these operators and on
the completeness of the generalized 
eigenfunctions.

For the purpose of making a comparison,
let us briefly recall some results
of Ikebe\cite{Ikebe} on the 
asymptotic behaviors of the
generalized eigenfunctions of 
the Schr\"odinger operator
\begin{eqnarray*}
-\Delta + V(x), \quad x \in
{\mathbb R}^3
\end{eqnarray*}
in connection with the eigenfunction
expansion for the 
absolutely continuous spectrum.
In \cite{Ikebe}, the generalized eigenfunction of
$-\Delta + V(x)$ was constructed 
as a solution to the
Lippmann-Schwinger equation
\begin{eqnarray}
\varphi(x, \, k) = e^{ix\cdot k}
 - \frac{1}{\, 4 \pi \,}
\int_{{\mathbb R}^3}
  \frac{\, e^{i|k| |x-y| \,}}{|x - y|}\,
  V(y) \, \varphi(y, \, k) \, dy,
\label{eqn:lpmnschw}
\end{eqnarray}
the solution being unique 
if $\varphi(x, \,k) - e^{ix\cdot k}$ belongs
to $C_{\infty}({\mathbb R}^3_x)$, 
the space of all continuous functions
vanishing at infinity.
Then the generalized Fourier transform,
of which kernel is
 the generalized eigenfunctions obtained,
was introduced and the
generalized Fourier inversion formula,
i.e., the eigenfunction expansion for
the absolutely continuous spectrum of
the operator $-\Delta + V(x)$ was established.

Ikebe's discussions on asymptotic behaviors
of the generalized eigenfunctions were based 
upon the Lippmann-Schwinger equation
(\ref{eqn:lpmnschw}). 
Roughly speaking, we see that his assumption
on the potential function is that
$V(x)$ is locally H\"older continuous and
$V(x) = O(|x|^{-\sigma})$, \ $\sigma >2$,
at infinity
(see Ikebe\cite[\S 1]{Ikebe} for 
the precise description of his assumption).

It is apparent that the term
\begin{eqnarray*}
\frac{1}{\, 4 \pi \,}
  \cdot \frac{\, e^{i|k| |x-y| \,}}{|x - y|}
\end{eqnarray*}
in (\ref{eqn:lpmnschw}) comes from
the integral kernel of the resolvent
of $-\Delta$:
\begin{eqnarray*}
(-\Delta - z)^{-1} u(x) =
\frac{1}{\, 4 \pi \,}
\int_{{\mathbb R}^3}
 \frac{\, e^{i \sqrt{z} \,|x-y| \,}}{|x - y|}\,
  u(y) \, dy,
\quad \hbox{Im}\, \sqrt{z\,} >0
\end{eqnarray*}
for 
$z \in {\mathbb C} \setminus [0, \, +\infty)$.
In other words,
the limiting absorption principle for
$-\Delta$ shows that
the boundary value of the resolvent
$(-\Delta -z)^{-1}$, as
$z= \lambda + i \mu \;\, 
(\lambda, \, \mu >0)$ tends to $\lambda + i 0$,
is expressed as the integral
operator
\begin{eqnarray*}
(-\Delta - \lambda - i0)^{-1} u(x) =
\frac{1}{\, 4 \pi \,}
\int_{{\mathbb R}^3}
 \frac{\, e^{i \sqrt{\lambda} \,|x-y| \,}}
  {|x - y|}\,
  u(y) \, dy.
\end{eqnarray*}

It was also shown in \cite{Ikebe}, 
by appealing to the
Lippmann-Schwinger equation (\ref{eqn:lpmnschw}),
that if $\sigma >3$ then the 
generalized eigenfunction has the asymptotics
\begin{eqnarray}
\varphi(x, \,k)=
e^{ix\cdot k} + 
f(|k|, \, \omega_x,  \, \omega_k ) \,
  \frac{\, e^{i|k| |x| \,} }{|x|} 
+ o(\frac{1}{\,|x|\,})
\label{eqn:planesherical}
\end{eqnarray}
as $|x| \to +\infty$, where
$\omega_x=x/|x|$, and $\omega_k=k/|k|$. 
From the view point of
physics, (\ref{eqn:planesherical}) is 
interpreted to mean that $\varphi(x, \,k)$ 
is asymptotically equal to a superposition
of the incoming plane wave $e^{ix\cdot k}$ and 
the outgoing spherical wave $e^{i|k| |x|}/|x|$
 (cf. Yafaev\cite[\S 1.3]{Yafaev}).

What we have recalled above indicates that
computing the integral kernel of
$(\sqrt{-\Delta} - z)^{-1}$ is 
naturally a starting point  to investigate 
asymptotic behaviors of 
the generalized eigenfunctions of
$\sqrt{-\Delta} + V(x)$.
Our computations show that the integral kernel
of the resolvent of $\sqrt{-\Delta}$ is
given by
\begin{eqnarray*}
(\sqrt{-\Delta} - z)^{-1} u(x) = 
  \int_{{\mathbb R}^3} g_z(x-y) \, u(y) \, dy
\end{eqnarray*}
for 
 $z \in {\mathbb C}\setminus[0, \, +\infty)$,
where
\begin{eqnarray*}
\lefteqn{g_z(x) =  \frac{\, 1 \,}{2 \pi^2|x|^2} \, +  } \\
\noalign{\vskip 4pt}
 && \quad
 \frac{\, z \,}{2 \, \pi^2 \, |x|}
\big[ \, \sin (z|x|) \, {\rm ci}(-z|x|) 
       - \cos(z|x|) \,{\rm si}(-z|x|)\,\big]                          
\end{eqnarray*}
(see Section 2).
For the definitions of the 
cosine  and sine 
integral functions ${\rm ci}(z)$ and
 ${\rm si}(z)$, see Subsection A.1 in Appendix.

The integral kernel $g_z (x-y)$ has nothing 
in common with the integral kernel of
$(-\Delta -z)^{-1}$, but if we take the limit
of $g_z (x-y)$ as $z$ approaches the positive
half of the real axis 
($z= \lambda + i\mu \to \lambda + i0$), 
then the term
\begin{eqnarray*}
\frac{\lambda}{\, 2 \pi \,}
  \cdot \frac{\, e^{i\lambda |x-y| \,}}{|x - y|}
\end{eqnarray*}
emerges as the leading term of 
$g_{\lambda + i0}(x-y)$, which
is actually the integral kernel 
of the boundary value 
$(\sqrt{-\Delta} - \lambda - i0)^{-1}$:
\begin{eqnarray*}
(\sqrt{-\Delta} - \lambda - i0)^{-1} u(x) = 
  \int_{{\mathbb R}^3} 
g_{\lambda + i0}(x-y) \, u(y) \, dy,
\quad \lambda >0,
\end{eqnarray*}
where
\begin{eqnarray}
g_{\lambda + i0}(x)&=&
\frac{\lambda}{\, 2 \pi \,}
  \cdot \frac{\, e^{i\lambda |x| \,}}{|x|} 
  + \frac{1}{\, 2\pi^2 |x|^2 \,}  
  + m_{\lambda}(x), 
\label{eqn:0g012}\\
\noalign{\vskip 4pt}
m_{\lambda}(x) &=& O(|x|^{-2})
\quad \hbox{as} \;\, |x| \to + \infty.
\nonumber
\end{eqnarray}
This fact enables us to investigate
asymptotic behaviors of the 
generalized eigenfunctions of 
$\sqrt{-\Delta} +V(x)$ by
utilizing the integral equation which we shall
call the modified Lippmann-Schwinger equation.

Unfortunately, the term $1/(2\pi^2 |x|^2)$ in
(\ref{eqn:0g012}) is quite troublesome.
The reason for this is that our generalized
eigenfunctions must be  bounded functions of
$x$ since they are expected to be  distorted
plane waves in physics terminology. 
However, the integral operator
\begin{eqnarray*}
\frac{1}{\, 2\pi^2 \,}
\int_{{\mathbb R}^3}  
 \frac{1}{\, |x-y|^2 \,} \,  u(y) \, dy,
\end{eqnarray*}
which is known as the Riesz potential,
cannot be a bounded operator
from $L^p({\mathbb R}^3)$ to
$L^{\infty}({\mathbb R}^3)$ for any $p\ge1$
(see Stein\cite[p.119]{Stein}).
To overcome this difficulty, 
we shall introduce a few inequalities for 
the Riesz potentials in Section 5.

We should like to remark here that one might
ignore the formula\linebreak
$-\Delta_x \, e^{ix\cdot k}= |k|^2 e^{ix\cdot k}$ 
plays a significant role in discussing
the generalized eigenfunction expansion for 
the Schr\"odinger operator $-\Delta + V(x)$,
because the formula is so trivial. 
On the contrary, it is far from trivial to 
show that
\begin{eqnarray}
\sqrt{-\Delta_x} \, e^{ix\cdot k} 
    = |k|e^{ix\cdot k} \quad
\hbox{ in the distribution sense}.
\label{eqn:rsplnwv}
\end{eqnarray}
Indeed, 
the left hand side of (\ref{eqn:rsplnwv}) 
is formally defined by
\begin{eqnarray*}
\int e^{ix\cdot\xi} \, |\xi| 
 \, \delta(\xi -k) \, d\xi,
\quad (\, \delta(\cdot)  
\hbox{ is the delta function}\,), 
\end{eqnarray*}
while the symbol $|\xi|$ of $\sqrt{-\Delta}$ 
is singular at the origin $\xi = 0$. 
Therefore, making sense of the expression
$\sqrt{-\Delta_x} \, e^{ix\cdot k}$ is one of the
main tasks in the present paper, and 
it will be accomplished in Section 8
with the aid of a theorem 
in Section 6.

\vspace{10pt}
{\bf Assumption } Throughout the paper
 we shall assume that
$V(x)$ is a real-valued 
 measurable function on 
${\mathbb R}^3$ satisfying
\begin{eqnarray}
|V(x)|  \le C \langle x \rangle^{-\sigma},  
\qquad \sigma >1,
\label{eqn:Vdecay}
\end{eqnarray}
though $\sigma$ will be required to satisfy
the assumption $\sigma >2$
when we investigate asymptotic 
behaviors of the generalized eigenfunctions
in precise manners.
We emphasize that we do not require any smoothness
assumption on the potential $V$. 
Although we could allow some local
singularities of $V$ in the sense that
$V(x)$ behaves like $|x - x_0|^{-\beta}$ with
$0 < \beta <1$ near some isolated points $x_0$'s,
 we shall not do so for the sake of
simplicity.

\vspace{10pt}

The plan of the paper is as follows. 
In Section 2, we compute the integral kernel
of the resolvent $(\sqrt{-\Delta} - z)^{-1}$
for 
$z \in {\mathbb C}\setminus[0, \, +\infty)$.
In Section 3, we derive expressions of the
boundary values 
$(\sqrt{-\Delta} - \lambda \mp i0 )^{-1}$
on the half positive axis in terms of
the boundary values 
$(-\Delta - \lambda \mp i0 )^{-1}$.
The expressions will be used in Section 6.
In Section 4, we compute the integral kernels
of $(\sqrt{-\Delta} - \lambda \mp i0 )^{-1}$.
In order to show that our generalized 
eigenfunctions are bounded functions, 
we shall prove some inequalities,
 in Section 5, for the Riesz potential
and the integral operator appearing as 
a part of
$(\sqrt{-\Delta} - \lambda \mp i0 )^{-1}$.
In Section 6, we establish the radiation
conditions for $\sqrt{-\Delta}$, which
implies that the second term of the
generalized eigenfunction of 
$\sqrt{-\Delta} + V(x)$ is a 
spherical wave in a certain sense.
In Section 7, we establish
the radiation
conditions for $\sqrt{-\Delta}+V(x)$,
which is of some interest on its own.
We construct the generalized eigenfunctions
of $\sqrt{-\Delta} + V(x)$, and
characterize them as unique solutions 
to the modified Lippmann-Schwinger equations
in Section 8. 
In Section 9, we show that the generalized
eigenfunctions are bounded functions of
$x$, and continuous functions of the
both variables $x$ and $k$.
Our discussions here are based on the
modified Lippmann-Schwinger equations.
In Section 10, we give estimates on the
difference between the generalized eigenfunction
and the plane wave when $\sigma >2$.
Also, we give estimates on the
difference between the generalized eigenfunction
and the sum of a plane wave 
and a spherical wave when $\sigma >3$.
In Appendix, we illustrate some properties
of the cosine and sine integral functions,
and prove inequalities for a convolution
which are used several times in the present
paper.

It is worthwhile to mention that all the
results and the discussions in Sections 3, 6
and 7 remain valid for the $n$-dimensional
case with $n\ge2$ with
trivial changes. However, we shall confine
our attention, throughout the present paper,
to the $3$-dimensional case for the sake
of clarity of description.

\vspace{10pt}
{\bf Notation}
We introduce the notation which will be used in the
present paper. Although the discussions in the
present paper will be made for the 3-dimensional
case, the notation introduced here are given
in the $n$-dimensional setting. 

For $x \in {\mathbb R}^n$,
$|x|$ denotes the Euclidean norm of $x$ and
\begin{eqnarray*}
\langle x \rangle = \sqrt{1 + |x|^2}  . 
\end{eqnarray*}
The Fourier transform of a function $u$ is 
denoted by ${\cal F}u$ or $\hat u$, and
defined by
$$
[{\cal F}u](\xi)= \hat u (\xi) = (2\pi)^{-n/2} 
   \int_{{\mathbb R}^n} e^{-ix\cdot \xi} \, u(x)  \, dx.
$$
For $s$ and $\ell$ in $\mathbb R$,  we define the
weighted $L^2$-space and the weighted Sobolev 
space by
\begin{eqnarray*}
L^{2, \,s}({\mathbb R}^n) = 
\{ f \; | \;  \langle x \rangle^s f \in 
L^2({\mathbb R}^n)   \;\}
\end{eqnarray*}
and
\begin{eqnarray*}
H^{\ell, \, s}({\mathbb R}^n) = 
\{ f \; | \;   \langle x \rangle^s
  \langle D \rangle^{\ell} f\in 
L^2({\mathbb R}^n)  \;\} 
\end{eqnarray*}
respectively, where $D$ stands for 
$-i\partial /\partial x$ and
$\langle D \rangle =
  \sqrt{1 + |D|^2}= \sqrt{1 -\Delta}$.
When $s=0$, we  write 
$L^2({\mathbb R}^n) =L^{2, \,0}({\mathbb R}^n)$ and
$H^{\ell}({\mathbb R}^n) =H^{\ell , \, 0}({\mathbb R}^n)$.
The inner products and the norms in
$L^{2, \,s}({\mathbb R}^n)$ and 
$H^{\ell, \, s}({\mathbb R}^n)$ are given by
\begin{equation*}
\begin{cases}  
(f, \, g)_{L^{2, \, s}} =
 \displaystyle \int_{{\mathbb R}^n} 
 \langle x \rangle^{2s} \, f(x) \,
    \overline{g(x)}  \, dx      &            \\
{}    & \\
\Vert f  \Vert_{L^{2, \, s}} = 
    \{ (f, \, f)_{L^{2, \, s}}   \}^{1/2}   &
\end{cases}
\end{equation*}
and
\begin{equation*}
\begin{cases}  
(f, \, g)_{H^{\ell, \, s}} =
 \displaystyle \int_{{\mathbb R}^n} 
  \langle x \rangle^{2s} 
  \,  \langle D \rangle^{\ell}f(x) \,
 \overline{ \langle D \rangle^{\ell}g(x)}  \, dx  &\\
{}    &\\
\Vert f  \Vert_{H^{\ell, \, s}} = 
    \{ (f, \, f)_{H^{\ell, \, s}}   \}^{1/2}   &
\end{cases}
\end{equation*}
respectively.

By $C_0^{\infty}({\mathbb R}^n)$ we mean the space of 
$C^{\infty}$-functions
of compact support. 
By ${\cal S}({\mathbb R}^n)$ we mean the Schwartz space
of rapidly decreasing functions, and
by ${\cal S}^{\prime}({\mathbb R}^n)$ the space of tempered
distributions. 
For a pair of  
$f \in {\cal S}^{\prime}({\mathbb R}^n)$ and
$\psi \in {\cal S}({\mathbb R}^n)$, we denote
the duality bracket by $\langle f, \, \psi\rangle$.
For a pair of 
$f \in L^{2, \,-s}({\mathbb R}^n)$ and
$g \in L^{2, \,s}({\mathbb R}^n)$, we define
the anti-duality bracket by
\begin{eqnarray*}
(f, \, g)_{-s, s} := \int_{{\mathbb R}^n} 
 f(x) \, \overline{g(x)} \, dx .
\end{eqnarray*}

For a pair of Hilbert spaces ${\cal H}$ and
 ${\cal K}$, ${\bf B}( {\cal H}, \,{\cal K})$
denotes the Banach space of all bounded linear 
operators from ${\cal H}$ 
to ${\cal K}$.
We set 
${\bf B}({\cal H})={\bf B}( {\cal H}, \,{\cal H})$. 

For a selfadjoint operator $T$ in a Hilbert space,
$\sigma(T)$ and $\rho(T)$ denote the spectrum of
$T$ and the resolvent set of $T$ respectively.
The point spectrum, i.e., the set of all eigenvalues 
of $T$, will be denoted by $\sigma_{\rm p}(T)$.
The essential spectrum, the continuous spectrum and
the absolutely continuous spectrum of $T$ will be
denoted by $\sigma_{\rm ess}(T)$,
$\sigma_{\rm c}(T)$ and
$\sigma_{\rm ac}(T)$ respectively.

\newpage
\section{Integral kernels
 of the resolvents of
\boldmath
$\mathnormal{H_0}$
\unboldmath
}

This section is devoted to the computation of the resolvent
kernel of $H_0=\sqrt{-\Delta}$ on ${\mathbb R}^3$. 
We shall start with the definition of the operator
$H_0$, and the description of its basic
properties from the view point of spectral theory.  

Let $H_0$ be the selfadjoint operator in 
$L^2({\mathbb R}^3)$ given by
$$
H_0 := \sqrt{-\Delta}  \quad 
\hbox{ with domain } H^1({\mathbb R}^3) .
$$
Since $H_0$ is unitarily equivalent, through the
Fourier transform $\cal F$, to the multiplication
operator by $|\xi| \times$ in $L^2({\mathbb R}^3_{\xi})$,
it follows from Kato\cite[p. 520, Example 1.9]{Kato}
that $H_0$ is absolutely continuous, and that
$$
\sigma(H_0) =\sigma_{\rm ac}(H_0) = [0, \, \infty).
$$
Furthermore, we see that $H_0$ restricted on 
$C_0^{\infty}({\mathbb R}^3)$
 is  essentially selfadjoint.
Indeed, with a $C^{\infty}-$function $\chi$
satisfying
\begin{equation*}
\chi(\xi) = 
\begin{cases}  
1  & \text{ if $\;\; |\xi| \le 1$, } \\
{}    & \\
0  & \text{ if $\;\;   |\xi| \ge 2$,} \\
\end{cases}
\end{equation*}
 we can decompose $\sqrt{-\Delta}$
into a regular part and a singular part:
$$
\sqrt{-\Delta} = (1 - \chi (D) )\sqrt{-\Delta}
  + \chi (D)\sqrt{-\Delta},
$$
which enables us to regard $\sqrt{-\Delta}$
as a sum of a
essentially selfadjoint operator
on $C_0^{\infty}({\mathbb R}^3)$ 
(see Nagase and Umeda\cite[Theorem 3.4]{NagaseUmeda})
and a bounded selfadjoint operator. The resolvent of 
$H_0$ will be 
denoted by
$$
R_0(z) = (H_0 - z)^{-1}   
     \qquad ( z \in 
        \rho(H_0) = {\mathbb C}\setminus [0, \, \infty) \; ) .
$$

By virtue of the fact that 
$R_0(z)={\cal F}^{-1} (|\xi| - z )^{-1} {\cal F}$.
 it would be possible to obtain the resolvent
kernel, i.e., the integral kernel of 
$R_0(z)$ 
by direct computation 
of $[{\cal F}^{-1}(|\xi| - z )^{-1}](x)$.
We shall, however, avoid this computation. 
Instead, we take advantage of 
 the fact 
that the strongly continuous semigroup generated 
by $-H_0$ is expressed as a convolution with the 
Poisson kernel 
(Stein\cite[p.61]{Stein}, 
Strichartz\cite[p.50]{Strichartz}):
$$
e^{-tH_0}u(x) = P_t * u (x)
  = \int_{{\mathbb R}^3} P_t(x -y ) u(y) \, dy,  
\quad t>0, \; u \in L^2({\mathbb R}^3),
$$
where 
\begin{equation}
P_t(x) = \frac{\, t \,}{\pi^2 \, ( t^2 + |x|^2 )^2 }.
\label{eqn:poisson2} 
\end{equation}
We then take the Laplace transform of $e^{-tH_0}$ 
to get the resolvent:
$$
R_0(z) = \int_0^{\infty} e^{tz}e^{-tH_0} \, dt   
      \qquad \hbox{if \ Re}\;z <0.
$$
Thus we need the following prerequisite.

\vspace{20pt}

\begin{lem} \label{lem:lt}
 If \ ${\rm Re} \, z <0$, then
\begin{eqnarray*}
\lefteqn{\int_0^{+\infty} e^{tz} 
\, \frac{\, t \,}{\pi^2 \, ( t^2 + a^2 )^2 }  \, dt  
 = \frac1{\, 2 \, \pi^2 \, a^2 \,} + }    \\
&&  \quad\quad\quad\quad\quad\quad
\frac{z}{\, 2 \, \pi^2  \, a \,} 
[ \, \sin (za) \, {\rm ci}(-za) 
            - \cos(za) \, {\rm si}(-za) \, ] ,    
\end{eqnarray*}
where $a$ is a positive constant.
\end{lem}
  
{\it Proof.} Since 
$$
\frac{t}{\, ( t^2 + a^2 )^2 \,} =
 \frac{d}{dt}\Big\{ 
    -\frac{1}{\,2( t^2 + a^2 )\,}  \Big\},
$$
we get, by integration by parts, 
\begin{eqnarray}
\int_0^{+\infty} e^{tz} 
\, \frac{\, t \,}{ \, ( t^2 + a^2 )^2  \,}  \, dt  
= \frac{1}{\, 2a^2 \,}   
  +  \frac{z}{\, 2 \,} 
\int_0^{\infty} e^{tz} 
\, \frac{\, 1 \,}{ \,  t^2 + a^2  \, }  \, dt .
\label{eqn:2a}       
\end{eqnarray}
Applying the formula (A.4) in Appendix to the 
integral on the right-hand side of (\ref{eqn:2a})
and noting the remark after
the formula (A.4), we obtain the lemma.
\hfill$\square$

\vspace{20pt}
In accordance 
with Lemma \ref{lem:lt},
we need to introduce  
two functions, which constitute 
the integral kernel of $R_0(z)$ as we shall see 
in Theorem \ref{thm:rk} below.

\vspace{20pt}
{\it Definition 2.1\ }
For 
$z \in {\mathbb C}\setminus [0,\,+\infty)$, 
we define
\begin{eqnarray}
\ell_z(x) &:=&  \frac{\, z \,}{2 \, \pi^2 \, |x|}
\big[ \, \sin (z|x|) \, {\rm ci}(-z|x|) 
            - \cos(z|x|) \, {\rm si}(-z|x|) \, \big] ,        
\label{eqn:2c}  \\
\noalign{\vskip 4pt}
g_z(x) &:=&  
\frac{\, 1 \,}{2 \pi^2|x|^2}+ \ell_z(x).                          
\label{eqn:2d}
\end{eqnarray}
By $G_0$ we mean the operator defined by
\begin{eqnarray}
G_0u(x):=\frac{\, 1 \,}{2 \, \pi^2 \, }
\int_{{\mathbb R}^3}\frac{\, 1 \,}{|x-y|^2} u(y) \, dy.
\label{eqn:2e}     
\end{eqnarray}
By  $G_z$ 
we mean the operator defined by
\begin{eqnarray}
 G_z u(x) := G_z*u (x) =
\int_{{\mathbb R}^3} g_z(x - y) u(y) \, dy.  
\label{eqn:2f}
\end{eqnarray}

\vspace{15pt}
\noindent
Note that $G_0$ is the Riesz potential.
See Stein\cite[p.117]{Stein}, in which $I_1$ is
the same as the operator $G_0$ in the present paper.  
Note also that (\ref{eqn:2c}), (\ref{eqn:2d})
and Lemma \ref{lem:lt} yield 
\begin{eqnarray}
\int_0^{+\infty} e^{tz} 
\, \frac{\, t \,}{\pi^2 \, ( t^2 + |x|^2 )^2 }  \, dt  
 =g_z(x)  
\label{eqn:2g}  
\end{eqnarray}
if $\hbox{Re}\, z <0$.

\vspace{15pt}

\begin{thm}\label{thm:rk}
If $z \in {\mathbb C}\setminus [0,\,+\infty)$,
then  
\begin{eqnarray*}
R_0(z) u  = G_z u                               
\end{eqnarray*}
for all 
$u \in C_0^{\infty}({\mathbb R}^3)$.
\end{thm}

{\it Proof.} It is sufficient to show
that
\begin{eqnarray}
(R_0(z) u, \, v)_{L^2}  = (G_z u, \, v)_{L^2}
\label{eqn:2h}                                           
\end{eqnarray}
for all
$z \in {\mathbb C} \setminus [0, \,+\infty)$
and all $u$, 
$v \in C_0^{\infty}({\mathbb R}^3)$. 

As mentioned before Lemma \ref{lem:lt}, we have
\begin{eqnarray}
( R_0(z)u, \, v)_{L^2} = 
  \int_0^{+\infty} e^{tz} \, 
( e^{-tH_0}u, \, v)_{L^2} \, dt  
   \qquad \qquad \qquad \qquad \qquad  \nonumber\\
\noalign{\vskip 3pt}
\qquad \qquad\qquad 
=  \int_0^{+\infty} e^{tz} \, 
\Big\{ \int_{{\mathbb R}^3} \!
\Big( \int_{{\mathbb R}^3} P_t(x -y ) u(y) \, dy \Big) \,
\overline{v(x)} \, dx  \Big\}  \, dt 
\label{eqn:2z10}   
\end{eqnarray}
if ${\rm Re}\, z < 0$. In order to 
make a change of order of integration in (\ref{eqn:2z10}),
we shall show
that the function
 $e^{tz}P_t(x -y ) u(y) \overline{v(x)}$
is absolutely integrable with respect the variables 
$x$, $y$ and $t$ if $\hbox{Re}\,z <0$
and $u$, 
$v \in C_0^{\infty}({\mathbb R}^3)$. 
To this end, we see 
(by integration by parts as in (\ref{eqn:2a})) that
\begin{eqnarray*}
\int_0^{+\infty} e^{t (\hbox{\scriptsize Re}\, z)} 
\, \frac{\, t \,}{ \, ( t^2 + a^2 )^2  \,}  \, dt  
&=& \frac{1}{\, 2a^2 \,}   
  +  \frac{\hbox{Re}\,z}{\, 2 \,} 
\int_0^{+\infty} e^{t(\hbox{\scriptsize Re}\, z)} 
\, \frac{\, 1 \,}{ \,  t^2 + a^2  \, }  \, dt \\
\noalign{\vskip 4pt}
&\le& \frac{1}{\, 2a^2 \,} + 
   \frac{|\hbox{Re}\,z|}{\, 2a^2 \,}
\int_0^{+\infty} 
    e^{t(\hbox{\scriptsize Re}\, z)} \, dt  \\
\noalign{\vskip 4pt}
&=& \frac{1}{\, a^2 \,}.
\end{eqnarray*}
This estimate, together with
(\ref{eqn:poisson2}), implies that
\begin{eqnarray*}
\lefteqn{ \int \!\!\!\! \int \!\!\!\! 
  \int_{{\mathbb R}^6 \times (0, \, \infty)}
\Big| \, e^{tz}P_t(x -y ) u(y) 
  \overline{v(x)}  \,\Big|\, dx \, dy \, dt }  \\
\noalign{\vskip 4pt}
& &
\le 
\frac{1}{\pi^2}
  \int \!\!\!\!  \int_{{\mathbb R}^6} 
  \frac{ \, | u(y)  v(x) | \,}{ |x - y|^2} \, dx \, dy    \\
\noalign{\vskip 4pt}
& &
=
\frac{1}{\pi^2}
 \int_{{\mathbb R}^3}   | v(x) | \, dx
   \Big(
\int_{|x -  y| \le 1} +  \int_{|x -  y| \ge 1}
  \Big)  \,
  \frac{ \, | u(y)| \,}{ |x - y|^2}  \, dy    \\
\noalign{\vskip 4pt}
& &
\le 
\frac{1}{\pi^2}
 \Vert v \Vert_{L^1}
   \Big(
\Vert u \Vert_{L^{\infty}}
\int_{|y| \le 1} \frac{1}{ \, |y|^2 \, }  \, dy
+ \Vert u \Vert_{L^1}
  \Big)  \,  < +\infty.
  \end{eqnarray*}
Therefore we can make a change of order of integration 
in (\ref{eqn:2z10}), 
and we get
\begin{eqnarray}
( R_0(z)u, \, v)_{L^2}   
 \qquad \qquad \qquad \qquad \qquad \qquad \qquad
\qquad \qquad\qquad   \nonumber\\
\quad 
= 
 \,\int_{{\mathbb R}^3}  
        \Big\{ \int_{{\mathbb R}^3}  
           \Big( \int_0^{+\infty} e^{tz}  
\frac{t}{\, \pi^2 ( t^2 + |x - y|^2 )^2 \,} \, dt  \Big) \,
u(y) \, dy \Big\}  
\overline{v(x)} \, dx 
\label{eqn:2y10}                
\end{eqnarray}
when $\hbox{Re }z <0$. If we apply Lemma \ref{lem:lt} to
the integral with respect to the $t$ 
variable in (\ref{eqn:2y10} ) and 
appeal to (\ref{eqn:2g}), 
we obtain
\begin{equation}
( R_0(z) u, \, v)_{L^2}= ( G_z u, \, v)_{L^2} 
   \quad \hbox{on }  
        \{ \, z  \in {\mathbb C}  
       \; | \; \hbox{Re }z < 0 \; \} .    
\label{eqn:2x11}
\end{equation}
Differentiating 
\begin{eqnarray*}
\iint_{{\mathbb R}^6} g_z(x-y)
  u(y) \, \overline{v(x)} \, dxdy
\end{eqnarray*}
with 
respect to $z$ under the 
sign of integration
(recall (\ref{eqn:2c}),
(\ref{eqn:2d}) and that
$u$, $v\in C_0^{\infty}({\mathbb R}^3)$),
we can deduce 
that $( G_z u, \, v)_{L^2}$ is 
a holomorphic function of $z$ in 
${\mathbb C} \setminus [0, \, +\infty)$. 
In view of the fact
that $( R_0(z) u, \, v)_{L^2}$ is also
a holomorphic function of $z$ in 
${\mathbb C} \setminus [0, \, +\infty)$,  
(\ref{eqn:2x11}) implies
that (\ref{eqn:2h}) holds 
on ${\mathbb C} \setminus [0, \, +\infty)$ for
all $u$, $v\in C_0^{\infty}({\mathbb R}^3)$.
\hfill$\square$

\vspace{15pt}
{\it Remark.} 
Since $R_0(z)$ is 
a bounded operator in $L^2({\mathbb R}^3)$ for 
any $z \in {\mathbb C}\setminus [0,\,+\infty)$,
Theorem \ref{thm:rk} implies that so is
$G_z$. On the other hand, 
it is a well-known 
fact (Stein\cite[Chapter V, \S 1.2]{Stein})
that the Riesz potential $G_0$ cannot 
be a bounded operator in $L^2({\mathbb R}^3)$.
This makes it difficult to show
directly from (\ref{eqn:2c})--(\ref{eqn:2f}) 
that $G_z$ is 
a bounded operator  in $L^2({\mathbb R}^3)$.

\newpage
\section{Properties of the resolvents of
\boldmath$\mathnormal{H_0}$
\unboldmath}

This section is devoted to investigating properties
of the resolvents of 
$H_0 = \sqrt{-\Delta}$.
We put  emphasis on 
expressions of  the extended resolvents 
$R^{\pm}_0(z)$ in the forms 
which will be useful for establishing
the radiation conditions for
$\sqrt{-\Delta}$ as well as
$\sqrt{-\Delta} + V(x)$.

We shall begin with the limiting
absorption principle for $\sqrt{-\Delta}$,
which assures the existence of the extended
resolvents $R^{\pm}_0(z)$, that is, the
existence of the boundary values of $R_0(z)$
on the positive axis.
The limiting
absorption principle for 
$\sqrt{-\Delta + m^2}$ was first proved 
by Umeda\cite{Umeda1} in the case
where $m>0$. The results in \cite{Umeda1}
were greatly generalized by
Ben-Artzi and Nemirovski\cite{BenArtzi-Nem}, where they were
able to treat
$\sqrt{-\Delta}$. Actually, Theorem \ref{thm:BA-N} below
is a corollary to results in Ben-Artzi and 
Nemirovski\cite[Section 2]{BenArtzi-Nem}, 
which is 
based on a general theory developed by Ben-Artzi and 
Devinatz\cite{BenArtzi-Dev}.

\begin{thm}[Ben-Artzi and Nemirovski] \label{thm:BA-N}
Let $s > 1/2$. Then 
\begin{description}
\item[\rm (i)] For any $\lambda >0$, there exist the limits
$$
R_0^{\pm}(\lambda) = \lim_{\mu \downarrow 0} 
   R_0(\lambda \pm i \mu)   \quad
  {\it in } \; {\bf B}(L^{2, \,s}, \; H^{1, \,-s}).
$$
\item[\rm (ii)] The operator-valued 
functions $R_0^{\pm}(z)$ defined by
\begin{equation*}
R_0^{\pm}(z) = 
\begin{cases}  
R_0(z) &\text{\it if $\;  z \in {\mathbb C}^{\pm}$}  \\
{}   & \\
R_0^{\pm}(\lambda)  &\text{\it if $\;  z=\lambda >0$ }
\end{cases}
\end{equation*}
are 
${\bf B}(L^{2, \,s}, \; H^{1, \,-s})$-valued 
continuous functions,
where ${\mathbb C}^+$ and ${\mathbb C}^-$ are 
the upper and the  lower
half-planes respectively{\rm :}
\begin{eqnarray*}
{\mathbb C}^{\pm} = 
 \{ \, z \in {\mathbb C} \; | \;  \pm\hbox{\rm Im } z >0 \; \}.
\end{eqnarray*}
\end{description}
\end{thm}

\vspace{15pt}

Theorem \ref{thm:lap0} below gives  representation formulae for 
the extended  resolvents $R_0^{\pm}(z)$ of
$\sqrt{-\Delta}$ in terms of 
the extended  resolvents 
${\varGamma}_0^{\pm}(z)$ of $-\Delta$
(see Agmon\cite[Section 4]{Agmon} for the limiting absorption
principle for
$-\Delta$). 
The advantage of Theorem \ref{thm:lap0} is that its representation
formulae are  convenient tools to
derive the radiation conditions for $\sqrt{-\Delta}$, which we shall
need in later sections.
It should be noted that Theorem \ref{thm:lap0} provides an 
alternative
proof of Theorem \ref{thm:BA-N}.

\vspace{15pt}

\begin{thm}\label{thm:lap0}
Let $s> 1/2$. 
Suppose that $b > a > 0$, and define 
$$
D_{ab} := \{ \, z= \lambda + i \mu \in 
  {\mathbb C} \;\, | \;\, a \le \lambda \le b, 
                \; |\mu| \le \frac{\, a \,}{2}
   \, \}.
$$
Then there exist operator-valued functions $A(z)$ and
$B(z)$ such that
\begin{description}
\item[\rm (i)]  $A(z)$ is a ${\bf B}(L^{2,s})$-valued 
continuous function on ${\mathbb C}$,
\item[\rm (ii)]  $B(z)$ is a
${\bf B}(L^{2,s}, \; H^{1,-s})$-valued 
continuous function on $D_{ab}$,
\item[\rm (iii)] $R_0^{\pm}(z) = 
{\varGamma}_0^{\pm}(z^2) \, A(z) + B(z)$
   for all $z \in D_{ab}^{\pm}$, where
$$
D_{ab}^{\pm}:= \{ \; z \in  D_{ab}  \; | \; 
                 \pm \hbox{\rm Im } z \ge 0  \; \}.
$$
\end{description}
\end{thm}

\vspace{10pt}
Following the idea in Umeda\cite[Section 2]{Umeda1}, 
we shall 
give a proof of Theorem \ref{thm:lap0} by means of
 a series of
lemmas. We first note that for 
$z \in {\mathbb C}^{\pm}$
\begin{eqnarray}
R_0(z) &=& {\cal F}^{-1} 
  \Big[ \frac{|\xi| + z}{\, |\xi|^2 - z^2 \,}  
\Big] {\cal F}   \nonumber\\
  &=& {\cal F}^{-1}  
           \Big[ \frac{1}{\, |\xi|^2 - z^2 \,}  \Big]
    {\cal F} \cdot 
     {\cal F}^{-1} \Big[ z + \gamma (\xi) |\xi|  
\Big] {\cal F}   \label{eqn:3a} \\  
            &&\quad  + \; {\cal F}^{-1} 
           \Big[ 
\frac{\, \big(1 - \gamma (\xi) \big)
   |\xi| \, }{|\xi|^2 - z^2}  
               \Big] 
                  {\cal F} ,   \nonumber
\end{eqnarray}
%
where $\gamma$ is a $C^{\infty}_0$-function, which
will be specified soon. 
It is easy to see that
\begin{eqnarray}
\frac{3}{\, 4 \,} a^2 \le \hbox{Re } z^2    \le b^2
                        \quad  \hbox{for }\forall z \in D_{ab}, 
\label{eqn:3b}
\end{eqnarray} 
and that
\begin{eqnarray}
 \pm \, \hbox{Im }z^2 > 0  \quad  \hbox{for }
      \forall z \in D_{ab} 
\cap {\mathbb C}^{\pm}
\label{eqn:3c}
\end{eqnarray}
In view of (\ref{eqn:3b}) and (\ref{eqn:3c}), 
we choose $\gamma \in C^{\infty}_0({\mathbb R}^3)$ so that
\begin{equation*}
\gamma(\xi) = 
\begin{cases}  
1   &\text{if $\;  
    \displaystyle{ \frac{\, 1 \,}{2} }
         a^2  \le |\xi|^2   
                        \le\frac{\, 3 \,}{2}  b^2$ }\\
{}    &\\
0  &\text{ if $\;  |\xi|^2 \le 
        \displaystyle{ \frac{1}{\, 4 \,} } a^2$
      or   $2 b^2  \le |\xi|^2$. }
\end{cases}
\end{equation*}
One can easily find that
\begin{equation}
\Big| \, |\xi|^2 - z^2 \, \Big|   \ge 
    \frac{1}{\, 4 \,} a^2  \quad \hbox{for } 
    \forall z \in D_{ab}, \;\; 
       \forall \xi \in \hbox{supp}[ 1- \gamma],
\label{eqn:3cc}
\end{equation}
and that
\begin{equation}
\Big| \, |\xi|^2 - z^2 \, \Big|   \ge 
    \frac{1}{\, 3 \,} |\xi|^2  \;\;\quad \hbox{if } 
     \, z \in D_{ab}, \;\; 
        |\xi|^2 \ge \frac{3}{\, 2 \,} b^2.
\label{eqn:3ccc}
\end{equation}
In accordance with (\ref{eqn:3a}), 
we now define $A(z)$ and $B(z)$ by
\begin{equation}
A(z) :={\cal F}^{-1} 
\Big[ z + \gamma (\xi) |\xi|  \Big] {\cal F}
 =zI + {\cal F}^{-1} 
\Big[ \gamma (\xi) |\xi|  \Big] {\cal F}
\label{eqn:3d}
\end{equation}
and
\begin{equation}
B(z) :={\cal F}^{-1} 
           \Big[ 
 \frac{\, (1 - \gamma (\xi) )|\xi| \, }{|\xi|^2 - z^2}  
        \Big] 
              {\cal F}
\label{eqn:3e}
\end{equation}
respectively. With
\begin{equation}
{\varGamma}_0(z) = ( - \Delta - z )^{-1},   
   \quad z \in  {\mathbb C} \setminus [0, \, +\infty),
\label{eqn:3f}
\end{equation}
we  have 
\begin{equation}
R_0(z) = {\varGamma}_0(z^2) \, A(z) + B(z) 
   \quad \hbox{for } \forall 
    z \in D_{ab} \hbox{ with Im } z \not= 0
\label{eqn:3g}
\end{equation}
by (\ref{eqn:3a}).  
In order to  treat  $A(z)$ and $B(z)$ in
weighted $L^2$-spaces and weighted Sobolev spaces, 
we need terminology and a boundedness result on
pseudo\-differential operators in 
these spaces.

\vspace{15pt}

\noindent
{\sc Definition.}   A $C^{\infty}$-function $p(x, \xi)$
on ${\mathbb R}^n \times{\mathbb R}^n$
 is said to be in the class  $S_{0,0}^{\, \mu \,}$
  \ $(  \, \mu \in {\mathbb R} \, ) $  if for any pair 
     $\alpha$  and $\beta$  of multi-indices 
there exists a constant
    $C_{\alpha\beta} \ge 0$ such that
$$
\Big| \Big(\frac{\partial}{\partial\xi}\Big)^{\alpha}                    
  \Big(\frac{\partial}{\partial x}\Big)^{\beta} 
p(x, \xi) \Big|
\le C_{\alpha\beta} \,
{\langle  \xi  \rangle}^{\mu}.      
$$
The class  $S_{0,0}^{\, \mu \,}$ is a Fr\'echet 
space equipped with 
     the seminorms
$$
|p|_{\ell}^{(\mu)} = \max_{|\alpha|,|\beta|\le\ell}
\sup_{x,\xi}\Big\{                                                   
   \Big| \Big(\frac{\partial}{\partial\xi}\Big)^{\alpha}                    
  \Big(\frac{\partial}{\partial x}\Big)^{\beta}         
p(x, \xi) \Big|
 {\langle  \xi  \rangle}^{-\mu}  \, \Big\}
    \quad\quad (\ell = 0, 1, 2, \cdots). 
$$       
For $p(x, \, \xi) \in S_{0,0}^{\, \mu \,}$, 
a pseudodifferential operator $p(x, \, D)$ is
 defined by
$$
p(x, \, D) u(x) = (2 \pi)^{-n/2} 
  \int e^{i x \cdot \xi} p(x, \, \xi) {\hat u} (\xi) \, d \xi.
$$

It is well-known 
(Kumano-go \cite[Theorem 1.3, p.57]{Kumano-go}) that
$p(x, \, D)$ maps ${\cal S} ( {\mathbb R}^n)$
 continuously into itself,
and by duality, maps 
${\cal S}^{\prime} ( {\mathbb R}^n)$ into itself.

\vspace{15pt}

\begin{lem}\label{lem:lap0}
 Let $p(x, \, \xi)$  belong to   
$S^{-m}_{0,\, 0}$
  for some integer $m \ge 0$,  and let 
 $s \in \mathbb R$.   
Then there exist a nonnegative constant $C = C_{ms}$ 
   and  a positive integer
  $\ell =  {\ell}_{ms}$  such that
$$
\Vert p(x, \, D) u \Vert_{H^{m,s}}  \,  \le   \,  
C  \, | p |_{\ell}^{(-m)}  \, 
   \Vert u \Vert_{L^{2,s}}    
$$
for all $u \in  {\cal S}({\mathbb R}^n)$.
\end{lem}

{\it Proof.}  We  first prove the lemma in the
case where $m=0$. 
If $s\ge 0$, the lemma
is a special case of \cite[Lemma 2.2]{Umeda1},
where 
${\langle  x  \rangle}^s p(x, \, D){\langle  x  \rangle}^{-s}$
was shown to be a bounded operator 
in $L^2({\mathbb R}^n)$, of which norm is estimated 
by a constant times $| p |_{\ell}^{(0)}$ with
some integer $\ell$. 

If $s<0$, we consider
${\langle  x  \rangle}^{-s} p^*(x, \, D)
{\langle  x  \rangle}^s$,
where
$p^*(x, \, D)$ is a formal adjoint operator
of $p(x, \, D)$ in the sense that
$$
(p(x, \, D)u, \, v)_{L^2} 
= (u, \, p^*(x, \, D)v)_{L^2},
\qquad u, \, v \in {\cal S}({\mathbb R}^n).
$$
It is well-known
that the symbol $p^*(x, \, \xi)$ of the operator
 $p^*(x, \, D)$ belongs to  $S^0_{0, \, 0}$ 
(see \cite[Theorem 2.6, p.74]{Kumano-go}), and 
 that each seminorm of $p^*(x, \, \xi)$ is 
estimated by a seminorm of $p(x, \, \xi)$
(see \cite[Theorem 2.5, p.73]{Kumano-go}).  Hence,
for all $u$ and $v$ in ${\cal S}({\mathbb R}^n)$, we have 
\begin{eqnarray*}
\Big|
({\langle  x  \rangle}^{s} 
p(x, \, D){\langle  x  \rangle}^{-s}u, \, v)_{L^2}
\Big|
&=&
\Big|
(u, {\langle  x  \rangle}^{-s} 
p^*(x, \, D){\langle  x  \rangle}^s v)_{L^2}
\Big|                                                  \\
& \le &
\Vert u \Vert_{L^2} \, 
 \Vert {\langle  x  \rangle}^{-s} 
p^*(x, \, D){\langle  x  \rangle}^{-(-s)} v\Vert_{L^2}     \\
& \le &
\Vert u \Vert_{L^2}
 \, C |p^*|^{(0)}_{\ell}  \, \Vert v \Vert_{L^2} \,
\qquad (\because -s >0)                                     \\
& \le &
\Vert u \Vert_{L^2} \, 
C^{\prime} |p|^{(0)}_{{\ell}^{\prime}}  \, 
\Vert v \Vert_{L^2},
\end{eqnarray*}   
where in the second inequality the result 
in the preceding paragraph was used. We have thus shown
that for $s <0$, the operator
${\langle  x  \rangle}^{s} p(x, \, D){\langle  x  \rangle}^{-s}$
is  bounded  
in $L^2({\mathbb R}^n)$, and its norm is estimated 
by a constant times $| p |_{{\ell}^{\prime}}^{(0)}$ with
some integer $\ell^{\prime}$.

All that remains is to 
prove the lemma in the case where $m$ is
a positive integer. This can be done in the same 
manner as in the proof of \cite[Lemma 2.2]{Umeda1}.
We omit the details.
\hfill$\square$

\vspace{15pt}

We now turn to the proof of Theorem \ref{thm:lap0}. 
Note that $s$ in Lemma \ref{lem:lap1} below 
 can be negative. This is due to Lemma \ref{lem:lap0}

\begin{lem}\label{lem:lap1}
For any $s\in \mathbb R$, $A(z)$ is
a ${\bf B}(L^{2,s})$-valued 
continuous function on ${\mathbb C}$.
\end{lem}

{\it Proof.} Since the support of the function $\gamma$
is away from the origin, it is evident that
$\gamma(\xi) |\xi| \in C^{\infty}_0({\mathbb R}^3)$, which 
one can regard as a subset of $S^0_{0, \, 0}$. 
Then it follows 
from Lemma \ref{lem:lap0}
that
$\gamma(D) |D|$ defines a bounded operator
in  $L^{2,s}({\mathbb R}^3)$. This immediately implies the lemma,
because of the fact that 
$A(z) = zI + \gamma(D) |D|$.
\hfill$\square$

\vspace{15pt}

\begin{lem}\label{lem:lap2}
For any $s\ge 0$, $B(z)$ is
a ${\bf B}(L^{2,s}, \, H^{1, -s})$-valued 
continuous function on $D_{ab}$.
\end{lem}

{\it Proof.} In order to decompose the symbol of $B(z)$ 
into a regular part and a singular part, we shall
use the same function 
$\chi\in C_0^{\infty}({\mathbb R}^3)$ as in the 
beginning of Section 2. 
We thus define
\begin{eqnarray*}
B_1(z) &:=& {\cal F}^{-1} 
           \Big[ 
    \frac{\, (\, 1 - \gamma (\xi) \,) \,
  (\, 1 - \chi (\xi) \, ) \, |\xi| \, }{|\xi|^2 - z^2}  
        \Big] 
              {\cal F},  \\
\noalign{\vskip 3pt}
B_2(z) &:=& {\cal F}^{-1} 
           \Big[ 
    \frac{\, (\, 1 - \gamma (\xi) \,) \,
      \chi (\xi)  \, |\xi| \, }{|\xi|^2 - z^2}  
        \Big] 
              {\cal F}. 
\end{eqnarray*}
It is obvious that
\begin{equation}
B(z) = B_1(z) + B_2(z).
\label{eqn:3h}
\end{equation}
Therefore, it is sufficient to show that
both $B_1(z)$ and $B_2(z)$ are
${\bf B}(L^{2,s}, H^{1, -s})$-valued 
continuous functions on $D_{ab}$.

As for $B_1(z)$, we note that the symbol of 
$B_1(z)$ is a $C^{\infty}$-function, and we
shall apply Lemma \ref{lem:lap0}. To this end,
we exploit  the inequalities (\ref{eqn:3cc}) and
(\ref{eqn:3ccc}), and obtain 
\begin{eqnarray}
 \Big| 
\Big( \frac{\partial}{\partial \xi}\Big)^{\alpha}
\Big\{
\frac{\, \big( 1 - \gamma (\xi) \big) \,
  \big( 1 - \chi (\xi) \big) \, |\xi| \, }
   {|\xi|^2 - z^2} 
\Big\} 
\Big|
\le   C_{\alpha} \, \langle\xi\rangle^{-1 -|\alpha|}
\label{eqn:3i} 
\end{eqnarray}
for all $\alpha$, 
where $C_{\alpha}$ is a constant independent of
$z \in D_{ab}$. It then follows from  (\ref{eqn:3i}) and
Lemma \ref{lem:lap0} with $m=1$ that for every 
$s \in \mathbb R$
\begin{equation}
\Vert B_1(z) u \Vert_{H^{1, \, s}} \le 
   C_s \Vert  u \Vert_{L^{2,\, s}},
\qquad  u \in {\cal S}({\mathbb R}^3),
\end{equation}
where $C_s$ is a constant independent of $z \in D_{ab}$.
Therefore, for each $z \in D_{ab}$, $B_1(z)$ can be
extended to 
a bounded operator from $L^{2, \, s}({\mathbb R}^3)$ 
to $H^{1,s}({\mathbb R}^3)$. 
In a similar fashion, we can see that for
$z$, $z^{\prime} \in D_{ab}$
\begin{eqnarray}
 \Big| 
\Big( \frac{\partial}{\partial \xi}\Big)^{\alpha}
\Big\{
\frac{\, (\, 1 - \gamma (\xi) \,) \,
  (\, 1 - \chi (\xi) \, ) \, |\xi| \, }{|\xi|^2 - z^2} 
&-&
\frac{\, (\, 1 - \gamma (\xi) \,) \,
  (\, 1 - \chi (\xi) \, ) \, |\xi| \, }
  {|\xi|^2 - {z^{\prime}}^2} 
\Big\} 
\Big|   \nonumber\\
 \quad 
\le   C_{\alpha} |z - z^{\prime}| \, 
   \langle\xi\rangle^{-3 -|\alpha|}
\label{eqn:3j} 
\end{eqnarray}
for all $\alpha$, 
where the constant $C_{\alpha}$
 is independent of $z$, $z^{\prime} \in D_{ab}$. 
Lemma \ref{lem:lap0} with $m=3$, together with (\ref{eqn:3j}),
gives
\begin{eqnarray*}
\Vert \{ B_1(z) - B_1(z^{\prime}) \}u \Vert_{H^{3, \, s}} 
\le 
C_s |z - z^{\prime}| \, \Vert  u \Vert_{L^{2, \,s}},
\qquad  u \in {\cal S}({\mathbb R}^3), 
\end{eqnarray*}
for every $s\in \mathbb R$, 
where $C_s$ is a constant  being uniform for
$z$, $z^{\prime} \in D_{ab}$.
In particular, $B_1(z)$ is
a ${\bf B}(L^{2,s}, \, H^{1,s})$-valued 
continuous function on $D_{ab}$ for every
 $s\in \mathbb R$. As a result, we can deduce
that $B_1(z)$ is
a ${\bf B}(L^{2,s}, \, H^{1,-s})$-valued 
continuous function on $D_{ab}$ for every
 $s\ge 0$.

As for $B_2(z)$, we exhibit it as a product of 
a pseudodifferential operator with
a smooth symbol and a pseudodifferential operator
with a singular symbol:
\begin{eqnarray*}
B_2(z) &=& {\cal F}^{-1} 
           \Big[ 
    \frac{\, (\, 1 - \gamma (\xi) \,) \,
     }{|\xi|^2 - z^2}  
        \Big] 
              {\cal F} \cdot
{\cal F}^{-1} 
           \Big[ 
 \chi (\xi)  \, |\xi| \,   \Big]   {\cal F}    \\
&=:& B_{2,1}(z) \cdot B_{2,2}.
\end{eqnarray*}
Note that $B_{2,1}(z)$ can be treated in a similar
fashion to $B_1(z)$, and one can deduce that
for every 
$s \in \mathbb R$
\begin{equation}
\Vert B_{2,1}(z) u \Vert_{H^{2, \, s}} 
 \le C_s \Vert  u \Vert_{L^{2, \,s}},
\qquad  u \in {\cal S}({\mathbb R}^3),
\end{equation}
where $C_s$ is a constant independent of $z \in D_{ab}$,
and that
\begin{eqnarray*}
\Vert \{ B_{2,1}(z) - B_{2,1}(z^{\prime}) \}u 
\Vert_{H^{4, \, s}} 
\le 
C_s^{\prime} |z - z^{\prime}| \, 
\Vert  u \Vert_{L^{2, \,s}},
\qquad  u \in {\cal S}({\mathbb R}^3), 
\end{eqnarray*}
for every $s\in \mathbb R$, 
where $C_s^{\prime}$ is a constant  independent of
$z$, $z^{\prime} \in D_{ab}$.
In particular, $B_{2,1}(z)$ is
a ${\bf B}(L^{2,s},\, H^{2,s})$-valued 
continuous function on $D_{ab}$ for every
 $s\in \mathbb R$. Taking into account the fact that
$\chi(\xi)|\xi|$ is a bounded function, we
see that for $s \ge 0$
\begin{eqnarray*}
\Vert B_{2,2} \, u \Vert_{L^{2, \,-s}}  &\le&  
                 \Vert B_{2,2} \, u \Vert_{L^2}     \\
\noalign{\vskip 4pt}
  &\le  & ( \max_{\xi} \chi(\xi)|\xi| \,) \,    
                   \Vert u \Vert_{L^2}         \\
 &\le& ( \max_{\xi} \chi(\xi)|\xi| \,) \, 
  \Vert u \Vert_{L^{2, \,s}}.
\end{eqnarray*}
Hence $B_{2,2} \in {\bf B}(L^{2,s},\, L^{2,-s})$ 
for every $s \ge 0$, which implies that
$B_2(z)= B_{2,1}(z)\,B_{2,2}$  is
a ${\bf B}(L^{2,s},\, H^{2,-s})$-valued 
continuous function on $D_{ab}$ 
for every $s\ge 0$.
Summing up the arguments, we have
completed the proof of the lemma.
\hfill$\square$

\vspace{15pt}
It is clear that we have actually showed the
following assertion in the proof of 
Lemma \ref{lem:lap2}.

\vspace{20pt}
\noindent
{\bf Corollary to Lemma 3.3} {\it 
 There exist 
a ${\bf B}(L^{2,s}, \, H^{1,s})$-valued 
continuous function $B_1(z)$ 
 on $D_{ab}$ for every
 $s\in \mathbb R$ and 
a  ${\bf B}(L^{2,s}, \, H^{2,-s})$-valued 
continuous function $B_2(z)$ 
 on $D_{ab}$ for every
 $s\ge 0$ such that
$B(z)u = B_1(z)u + B_2(z)u$ for all
$u \in L^{2,s}({\mathbb R}^3)$ with $s\ge 0$.
}

\vspace{20pt}

\textbf{Proof of Theorem \ref{thm:lap0}}
\ Assertions (i) and (ii) are  special cases of 
Lemmas \ref{lem:lap1} and  \ref{lem:lap2} respectively,
since we assume $s > 1/2$ in the theorem. To prove 
assertion (iii), we recall a well-known result by
Agmon\cite[Theorem 4.1]{Agmon}: the extended resolvents
${\varGamma}_0^{\pm}(z)$
defined by
\begin{eqnarray}
{\varGamma}_0^{\pm}(z) = \left\{
\begin{array}{ll}
{\varGamma}_0(z)   &  
\mbox{if  $\; z \in {\mathbb C}^{\pm}  $}    \\
&  \\
{\varGamma}_0^{\pm}(\lambda)    
&  \mbox{if  $\;  z = \lambda >0$ }
\end{array}
\right.
\label{eqn:3k}
\end{eqnarray}
are ${\bf B}(L^{2,s},\, H^{2,-s})$-valued 
continuous function on 
${\mathbb C}^{\pm} \cup (0, \, +\infty)$ 
provided that
 $s > 1/2$. In view of assertions (i) and (ii),
the theorem follows from this fact and (\ref{eqn:3g}), 
together with (\ref{eqn:3c}).
\hfill$\square$

\vspace{15pt}

It is worthwhile to improve assertion (ii) of 
Theorem \ref{thm:lap0}.

\begin{thm}\label{thm:lap1}
Under the same assumptions and with the same
notation as in Theorem \ref{thm:lap0},
the operator-valued function 
$B(z)$ has the following
property{\rm:} If  $0 \le s < 5/2$ and  
$t <  s -3/2 $, then $B(z)$ is a 
${\bf B}(L^{2,s}, \; H^{1,t})$-valued 
continuous function on $D_{ab}$.
\end{thm}

{\it Proof.}
We utilize the decomposition (\ref{eqn:3h}) of $B(z)$
made in the proof of Lemma \ref{lem:lap2},
where it was actually shown that
$B_1(z)$ is a ${\bf B}(L^{2,s}, \; H^{1,s})$-valued 
continuous function on $D_{ab}$ for 
any $s \in {\mathbb R}$
(see Corollary to Lemma \ref{lem:lap2}).
 It is therefore sufficient
to prove that
$B_2(z)$ has the property described in the theorem.

We use the same factorization as in the
proof of Lemma \ref{lem:lap2}:
$B_2(z) = B_{2,1}(z)\, B_{2,2}$. Apparently, we 
have shown in the proof of Lemma \ref{lem:lap2} 
that $B_{2,1}(z)$ is a 
 ${\bf B}(L^{2,s}, \; H^{2,s})$-valued 
continuous function on $D_{ab}$ for 
every $s \in {\mathbb R}$. 
Since 
$B_{2,2}$ is equal to 
a pseudodifferential operator 
$\chi(D) \sqrt{-\Delta}$, 
we can apply Umeda\cite[Lemma 5.2]{Umeda2}.
 Thus we see that 
$B_{2,2} \in {\bf B}(L^{2,s}, \; L^{2,t})$
if $0 \le s < n/2+1$ and $t < s - n/2$. 
It then follows that $B_2(z)$ is a
${\bf B}(L^{2,s}, \; H^{2,t})$-valued continuous function
on $D_{ab}$ under the assumption of 
the theorem. 
\hfill$\square$

\newpage
\section{Integral kernels of 
\boldmath $\mathnormal{R_0^{\pm}(\lambda)}$
\unboldmath}

In this section, we shall derive 
the integral kernels 
of the boundary values $R_0^{\pm}(\lambda)$ of 
the resolvent
$R_0(z)$ on the positive half axis
(recall that
 the existence of $R_0^{\pm}(\lambda)$
was assured in the previous section). 
We have to start with
examining the boundary values of 
the complex variable function $\hbox{ci}(-z)$, 
$z \in {\mathbb C} \setminus [0, \, +\infty)$, 
since
the integral kernel $g_z(x)$ of $R_0(z)$ 
contains  the term 
$\hbox{ci}(-z|x|)$ as was shown in
Section 2. 
In connection with 
the integral kernel $g_z(x)$, 
it is worthwhile noting that
all of $\sin (z)$, $\cos (z)$ and
$\hbox{si}(z)$ are entire functions,
but $\hbox{ci}(z)$ is a many-valued function with
a logarithmic branch point at $z=0$;
we shall choose the principal branch
(see Subsection A.1 in Appendix).

By (A.1) in Appendix and the definition of 
the function $h_e(z)$ introduced in Appendix, 
we have 
\begin{eqnarray}
{\rm ci}(-z) =
    - i \, {\rm Arg} \,
(-z) -\gamma - {\rm log} |z| + h_e(z)
\label{eqn:3l}
\end{eqnarray}
for 
$z\in {\mathbb C}\setminus [0, \, +\infty)$.
It follows from (\ref{eqn:3l}) that 
if $\lambda >0$, then
\begin{eqnarray}
{\rm ci}(-(\lambda \pm i \mu) ) \to  
           \pm i \pi + {\rm ci}(\lambda ) 
\quad \hbox{ as } \mu \downarrow 0,
\label{eqn:3m}
\end{eqnarray}
where we have used that fact that $h_e$ is an even function, 
as is remarked
in Appendix.

We now turn to the boundary values of $\ell_z(x)$ on
the positive axis (see (\ref{eqn:2c}) for the
definition of ${\ell}_z(x)$).
Putting $z= \lambda \pm i \mu$ with $\lambda, \, \mu >0$,
we
take the limit of $\ell_z(x)$ as $\mu \downarrow 0$.
 We then see that
\begin{eqnarray}
 \ell_z(x) &\to &
\frac{\, \lambda \,}{2 \, \pi^2 \, |x|}
\big[\, \sin (\lambda |x|) \, 
   \{ \pm i \pi + {\rm ci}(\lambda |x|)   \}    \nonumber\\
 &&   \qquad   \qquad \qquad \quad
- \cos(\lambda |x|) \, 
   \{ - \pi - {\rm si}(\lambda |x|) \} \, \big]
\label{eqn:3n}
\end{eqnarray}
for each $x \not= 0$
as $\mu \downarrow 0$, 
where we have used (\ref{eqn:3m}) and (A.3) in Appendix.
By the fact that 
$e^{\pm i\lambda |x|} = \cos(\lambda |x|) \pm i\sin(\lambda |x|)$,
we get
\begin{eqnarray}
\ell_z(x) \to \frac{\lambda}{\, 2\pi \,}
  \cdot\frac{e^{\pm i \lambda |x|}}{\, |x| \,}
 + m_{\lambda}(x),
\label{eqn:3o}
\end{eqnarray}
for each $x \not= 0$ 
as $\mu \downarrow 0$, where
\begin{eqnarray}
m_{\lambda}(x):=
 \frac{\, \lambda \,}{2 \, \pi^2 \, |x|}
\Big[ \, \sin (\lambda |x|) \, {\rm ci}(\lambda |x|) 
 + \cos(\lambda |x|) \, {\rm si}(\lambda |x|) \, \Big].
\label{eqn:3q}
\end{eqnarray}
In accordance with (\ref{eqn:2d}) in Section 2,
 we define 
\begin{eqnarray}
g_{\lambda}^{\pm}(x) := \frac{\, 1 \,}{2 \pi^2|x|^2}
+ \frac{\lambda}{\, 2\pi \,}
  \cdot\frac{e^{\pm i \lambda |x|}}{\, |x| \,} + 
m_{\lambda}(x).
\label{eqn:3r}
\end{eqnarray}
(Recall that $g_{\lambda + i0}(x)$ in
Introduction, which is exactly the 
same as  $g_{\lambda}^+(x)$
defined above.)
It follows immediately from (\ref{eqn:2d}),
 (\ref{eqn:3o}) and (\ref{eqn:3r}) 
that for $\lambda > 0$
\begin{eqnarray}
g_{\lambda \pm i \mu}(x)  \to g_{\lambda}^{\pm}(x),
\qquad x \not= 0 
\label{eqn:3r+}
\end{eqnarray}
as $\mu \downarrow 0$.
From the view point of the
time independent theory of scattering, it is 
very important that the leading term of (\ref{eqn:3r})
at infinity 
is the second term 
$\lambda e^{\pm i \lambda |x|}/(2\pi |x|)$, 
which is 
the same, up to a constant,
 as the integral kernels
of the boundary values of the
resolvent ${\varGamma}_0(z)$ of 
 $-\Delta$ on ${\mathbb R}^3$.

We finally state a result on the integral 
representations of the boundary values of the 
resolvent $R_0(z)$.

\vspace{15pt}

\begin{thm}\label{thm:R0ker}
Let $s > 1/2$. If $\lambda >0$, then
\begin{eqnarray*}
( R_0^{\pm} (\lambda) u, \, v )_{-s,s}
= \int_{{\mathbb R}^3} \Big\{ 
\int_{{\mathbb R}^3} g_{\lambda}^{\pm}(x-y) \,u(y) \, dy
   \Big\}  \,  \overline{v(x)} \, dx  
\end{eqnarray*}
for all $u$ and $v \in C_0^{\infty}({\mathbb R}^3)$.
\end{thm}

{\it Proof.} It follows from (\ref{eqn:2f}) and 
Theorem \ref{thm:rk} that
\begin{eqnarray}
(R_0(\lambda \pm i\mu)u, \, v)_{L^2} =
\int_{{\mathbb R}^3}
 \Big\{ \int_{{\mathbb R}^3}
 g_{\lambda \pm i \mu}(x-y) \,u(y) \, dy
   \Big\}  \,  \overline{v(x)} \, dx , 
\label{eqn:3t}
\end{eqnarray}
where $\mu >0$. Since $R_0^{\pm}(z)$ defined 
in Theorem \ref{thm:BA-N} are 
${\bf B}(L^{2,s}, \, L^{2,-s})$-valued 
continuous functions on 
${\mathbb C}^{\pm}\cup (0, \, +\infty)$ respectively, 
we see that
\begin{eqnarray}
(R_0(\lambda \pm i\mu)u, \, v)  \to 
 ( R_0^{\pm}(\lambda)u, \, v )_{-s,s}
\label{eqn:3u}
\end{eqnarray}
as $\mu \downarrow 0$. As for the right hand side of
(\ref{eqn:3t}), we shall apply the Lebesgue dominated 
convergence theorem. To this end, we first 
note that $g_z(x)$ is locally
integrable. More precisely,
 in view of (\ref{eqn:2c}), (\ref{eqn:2d}), 
(\ref{eqn:3l}) and the fact that $h_e(z)$ and
$\hbox{\rm si}(z)$ are entire functions, we find 
 that for each pair of $\lambda>0$
and $a >1$, there corresponds a 
positive constant
$C_{\lambda a}$,  independent of $\mu$ 
with $0< \mu <1$, such that
\begin{eqnarray}
| g_{\lambda \pm i \mu} (x) | \le  C_{\lambda a}
  \left\{   \begin{array}{ll}
  1/|x|^2   &  \mbox{ if  $\; |x| \le 1$}     \\
  {}& {}                                   \\
  1   &  \mbox{ if  $\; 1\le |x| \le a$.}  
  \end{array}
  \right. 
\label{eqn:3u+}
\end{eqnarray}
Since $u$ and $v$ lie in
$C_0^{\infty}({\mathbb R}^3)$,
it follows from (\ref{eqn:3u+}) that
\begin{eqnarray}
| g_{\lambda \pm i \mu} (x-y) u(x) \overline{v(y)} | 
\le  C_{\lambda uv}
  \left\{   \begin{array}{ll}
  \displaystyle{\frac{|u(x)v(y)|}{|x-y|^2} }  
&  \mbox{ if  $\; |x -y| \le 1$}     \\
  {}& {}                                   \\
  |u(x)v(y)|   &  \mbox{ otherwise}  
  \end{array}
  \right. 
\label{eqn:3u++}
\end{eqnarray}
where $C_{\lambda uv}>0$ is a constant, 
being dependent on $\lambda$,
$u$ and $v$, but independent 
of $\mu$ 
with $0< \mu <1$. Note that the function on 
the right hand side of (\ref{eqn:3u++}) is integrable on
${\mathbb R}^3 \times {\mathbb R}^3$.
By virtue of (\ref{eqn:3r+}) and (\ref{eqn:3u++}), we 
can apply
the Lebesgue dominated convergence theorem, and see
that
\begin{eqnarray}
\int \!\!\! \int_{{\mathbb R}^6}
 g_{\lambda \pm i \mu} (x-y) u(x) \overline{v(y)} 
  \, dx dy
\to  
\int \!\!\! \int_{{\mathbb R}^6}
 g_{\lambda}^{\pm} (x-y) u(x) \overline{v(y)}  
  \, dx dy
\label{eqn:3v}
\end{eqnarray}
as $\mu \downarrow 0$. Combining (\ref{eqn:3t}) with
 (\ref{eqn:3u}), (\ref{eqn:3v}) gives the theorem.
\hfill$\square$

\vspace{20pt}

It follows from Theorem \ref{thm:R0ker} that
the integral operators defined by
\begin{equation}
G_{\lambda}^{\pm}u (x) :=
\int_{{\mathbb R}^3} 
    g_{\lambda}^{\pm}(x-y) \,u(y) \, dy,
\quad \;\; u \in C_0^{\infty}({\mathbb R}^3) 
\label{eqn:4gg}
\end{equation}
can be extended to bounded  operators
from 
$L^{2, \, s}({\mathbb R}^3)$ to
$H^{1, \, -s}({\mathbb R}^3)$ for $s > 1/2$,
since  
$R_0^{\pm}(\lambda) 
\in {\bf B}(L^{2, \, s}, \, H^{1, \, -s})$
 for $s > 1/2$,
and 
\begin{equation}
R_0^{\pm}(\lambda)u = G_{\lambda}^{\pm}u, 
\quad \;\; u \in C_0^{\infty}({\mathbb R}^3).
\label{eqn:4-last}
\end{equation}
%

\newpage
\section{Estimates on the integral operators}

In this section, we consider 
the Riesz potential $G_0$ (see (\ref{eqn:2e}))
 and the integral
operators $K_{\lambda}^{\pm}$, $M_{\lambda}$ defined by
\begin{eqnarray}
(K_{\lambda}^{\pm}u)(x) &:=& 
   \frac{\lambda}{\, 2\pi \,}\int_{{\mathbb R}^3} 
 \frac{e^{\pm i \lambda |x-y|}}{|x-y|} \, u(y) \, dy,  
 \label{eqn:4a}\\
\noalign{\vskip 6pt}
(M_{\lambda}u )(x) &:=& 
     \int_{{\mathbb R}^3} m_{\lambda}(x-y) \, u(y) \, dy.
\label{eqn:4a+}
\end{eqnarray}
(For the 
definition of  $m_{\lambda}(x)$,
see (\ref{eqn:3q}).)
Our task here is to derive estimates of these operators
in weighted $L^2$-spaces as well as pointwise estimates
of $(G_0u)(x)$, $(K_{\lambda}^{\pm}u)(x)$ and 
$(M_{\lambda}u)(x)$ for 
$u$ belonging to some weighted $L^2$-space or to a suitable
class of functions. We shall apply these estimates 
 in the later sections in order to
examine asymptotic behaviors of the
generalized eigenfunctions of 
$\sqrt{-\Delta} + V(x)$ on ${\mathbb R}^3$.
In connection with this, it is important to notice
that we have formal identities
\begin{eqnarray}
R_0^{\pm}(\lambda) =
G_{\lambda}^{\pm} = 
G_0 + K_{\lambda}^{\pm} +  M_{\lambda},
\label{eqn:5ahlalah}
\end{eqnarray}
which hold at least on
$C^{\infty}_0{({\mathbb R}^3})$;
see (\ref{eqn:3r}), (\ref{eqn:4gg})
and (\ref{eqn:4e+}).

It is well-known (Stein\cite[p.119]{Stein}) that 
the inequality
\begin{eqnarray*}
\Vert G_0 u \Vert_{L^{\infty}} \le  C  
   \Vert  u \Vert_{L^p}
\end{eqnarray*}
cannot hold for any $p \ge 1$.  Furthermore, 
we make a remark that if one defines 
\begin{equation*}
u_0(x) := 
\begin{cases}  
1/|x|  \quad {}\;\;  
    ( |x| \le 1 )  &\\
{}    &\\
0  \quad \;\;\;\;(\hbox{\rm otherwise}), &\\
\end{cases}
\end{equation*}
then $u_0 \in L^{2,s}({\mathbb R}^3)$ for all 
$s \in {\mathbb R}$, and $(G_0 u_0)(0) = +\infty$. 
In spite of these facts, we need to find a class
of functions $u$ for which 
$(G_0 u)(x)$ are bounded functions of $x$. 
Actually, we shall obtain two 
sufficient conditions 
(see Lemmas \ref{lem:ker1+} and \ref{lem:ker1++} below),
either of which is suitable for  showing the
boundedness of generalized
eigenfunctions of $\sqrt{-\Delta} + V(x)$ on 
${\mathbb R}^3$.
It is also well-known (Stein \cite[p. 119]{Stein}) that 
the inequality
\begin{eqnarray*}
\Vert G_0 u \Vert_{L^q} \le  C  
   \Vert  u \Vert_{L^p}
\end{eqnarray*}
holds only if $q^{-1}= p^{-1} -3^{-1}$  in the context of 
the present paper.  When $p=2$, we actually have
\begin{eqnarray}
\Vert G_0 u \Vert_{L^6} \le  C  
   \Vert  u \Vert_{L^2}.
\label{eqn:4b-}
\end{eqnarray}
On the other hand, we are going to show a few 
boundedness 
results on $G_0$ in the framework of weighted
$L^2$-spaces as well as in some other frameworks.

\vspace{15pt}

\begin{lem}\label{lem:ker1}
Let $s> 3/2$. Then 
\begin{description}
\item[\rm (i)] $G_0 \in {\bf B}(L^{2, \, s}, \,L^2)$.
\item[\rm (ii)] $G_0 \in {\bf B}(L^2, \,L^{2, \, -s})$.
\end{description}
\end{lem}

{\it Proof.} Let $u \in L^{2, \, s}({\mathbb R}^3)$.
Since $s > 3/2$, the Schwarz inequality
gives
\begin{eqnarray}
\int_{{\mathbb R}^3} | u(x) | \, dx =
\int_{{\mathbb R}^3}  \langle x \rangle^{-s}  \cdot
     \langle x \rangle^{s} | u(x) | \, dx  \le  
  C_s \Vert u \Vert_{L^{2,s}},
\label{eqn:4b}
\end{eqnarray}
hence $u \in L^1({\mathbb R}^3)$. 
With $B=\{ \, x \; | \; |x| \le 1 \, \}$ 
and $E=\{ \, x \; | \; |x| \ge 1 \, \}$, we decompose
the function $1/(2 \pi^2 |x|^2)$ into two parts:
\begin{eqnarray}
\frac{1}{2 \pi^2 |x|^2}  &=& 
   \frac{1_B(x)}{2 \pi^2 |x|^2} +  
               \frac{1_E(x)}{2 \pi^2 |x|^2} \nonumber\\
&=:&  h_B(x)  +  h_E(x),
\label{eqn:4c}
\end{eqnarray}
where $1_B(x)$ and $1_E(x)$ are the characteristic functions of
the sets $B$ and $E$ respectively.
It is clear that $h_B(x)\in L^1({\mathbb R}^3)$ and  
$h_E(x)\in L^2({\mathbb R}^3)$, and that
\begin{eqnarray}
G_0u = h_B* u +  h_E*u. 
\label{eqn:4c+}
\end{eqnarray}
If we regard $u$ as 
a function belonging to $L^2({\mathbb R}^3)$,  
we can
apply the Young inequality (see Stein\cite[p.271]{Stein})
to $h_B* u$, and obtain
\begin{eqnarray}
\Vert h_B* u \Vert_{L^2} \le  
    \Vert h_B\Vert_{L^1}   \Vert u \Vert_{L^2}  \le
    \Vert h_B\Vert_{L^1}   \Vert u \Vert_{L^{2,s}}.
\label{eqn:4d+}
\end{eqnarray}
If we regard $u$ as 
a function belonging to
 $L^1({\mathbb R}^3)$
 (recall (\ref{eqn:4b})),  
we can also
apply the Young inequality
to $h_E* u$, and obtain
\begin{eqnarray}
\Vert h_E* u \Vert_{L^2} \le  
    \Vert h_E\Vert_{L^2} \Vert u \Vert_{L^1}    \le
    \Vert h_E\Vert_{L^2}   \Vert u \Vert_{L^{2,s}},
\label{eqn:4d++}
\end{eqnarray}
where we have used (\ref{eqn:4b}). Combining  
(\ref{eqn:4c+})--(\ref{eqn:4d++}), we conclude
that assertion (i) is true.

To prove assertion (ii), we note that $G_0$ is
symmetric on $C_0^{\infty}({\mathbb R}^3)$:
\begin{equation*}
( G_0 u, \, v)_{L^2} =
( u, \,  G_0 v)_{L^2}   \quad 
\hbox{for }  
u, \; v \in C_0^{\infty}({\mathbb R}^3),
\end{equation*}
which, together with assertion (i),
 implies 
\begin{equation}
|(  u, \, G_0 v)_{L^2}| \le 
\Vert G_0 u \Vert_{L^2} \Vert v\Vert_{L^2}
\le C \Vert u \Vert_{L^{2, s}} \Vert v\Vert_{L^2}
\label{eqn:4g0}
\end{equation}
for all $u, \; v \in C_0^{\infty}({\mathbb R}^3)$.
We can regard the left hand side of 
(\ref{eqn:4g0})
as the anti-duality
bracket
$(  u, \, G_0 v)_{s, -s}$.
Hence, by the density argument, it follows from
(\ref{eqn:4g0}) that 
\begin{equation*}
\Vert G_0 v \Vert_{L^{2, -s}} \le 
C \Vert v \Vert_{L^2}
\end{equation*}
for all $v \in C_0^{\infty}({\mathbb R}^3)$.
This yields assertion (ii).
\hfill$\square$

\vspace{15pt}

\begin{lem}\label{lem:ker1+}
If $u$ satisfies 
\begin{eqnarray}
|u(x)|  \le C \langle x \rangle^{-\ell},  
\quad \ell >1, \; C>0,
\label{eqn:4de}
\end{eqnarray}
then 
\begin{eqnarray*}
  | G_0 u(x)|                                       
\le C_{\ell} \, 
\Vert \langle \cdot \rangle^{\ell}u \Vert_{ L^{\infty} } 
\times
\left\{  
\begin{array}{ll}
\langle x \rangle^{-(\ell -1)}  
     & \mbox{ if $\; 1  < \ell < 3$},  \\
{}  &  {}             \\
\langle x \rangle^{-2} \log (1 + \langle x \rangle )   
       & \mbox{ if $\;   \ell = 3$},  \\
{}  & {}  \\
\langle x \rangle^{-2}  
       & \mbox{ if $\;   \ell > 3$} . 
\end{array}
\right.
\end{eqnarray*}
\end{lem}

\vspace{6pt}

{\it Proof.} 
It is evident from
the definition (\ref{eqn:2e})  
that we have
\begin{eqnarray}
|G_0u(x)| \le \frac{1}{\, 2 \pi^2 \,} 
   \Vert \langle \cdot \rangle^{\ell} u \Vert_{L^{\infty}}
\int_{{\mathbb R}^3}  \!
\frac{1}{\, |x - y|^2 \langle y \rangle^{\ell} \,} \, dy.
\label{eqn:4de++}
\end{eqnarray}
If we apply Lemma A.1 in Appendix, with
$n=3$, $\beta = 2$ and $\gamma=\ell$, 
to the function
defined by the integral on the
right hand side of (\ref{eqn:4de++}), then
the lemma follows.
\hfill$\square$

\vspace{15pt}

\begin{lem}\label{lem:ker1++}
Suppose that 
\begin{eqnarray}
u \in L^2({\mathbb R}^3) \cap L^q({\mathbb R}^3),  
\quad q >3.
\label{eqn:4de+12}
\end{eqnarray}
Then there exists a constant $C_q$, independent of 
$u$, such that
\begin{eqnarray}
\Vert G_0 u \Vert_{L^{\infty}}  
\le 
C_q (\Vert  u \Vert_{L^2} + \Vert  u \Vert_{L^q} )
\label{eqn:4de+13}
\end{eqnarray}
\end{lem}

{\it Proof.}
We exploit the same 
decomposition of $G_0u$ 
as in (\ref{eqn:4c+}).
If we apply
 the H\"older inequality to $h_B*u$, 
we obtain
\begin{eqnarray}
|h_B*u (x)| \le \frac{1}{\, 2\pi^2 \,} \,
\Big\{ \int_{|x-y| \le 1}
  \frac{1}{\, |x-y|^{2p} \,} \, dy \Big\}^{\!1/p}  \,
  \Vert u  \Vert_{L^q},
\label{eqn:4de+-1}
\end{eqnarray}
where $p^{-1} = 1 - q^{-1}$.  Since $q >3$, it
follows that $2p < 3$. Hence 
the inequality (\ref{eqn:4de+-1}),
 together with the assumption (\ref{eqn:4de+12}), 
implies that $h_B*u(x)$ is a bounded function.
Similarly, if we apply the Schwarz inequality 
to $h_E*u$, we can deduce that 
$h_E*u(x)$ is a bounded function.  Summing up, 
we have shown the inequality (\ref{eqn:4de+13}).
\hfill$\square$

\vspace{20pt}

In order to derive estimates of the operator
$M_{\lambda}$, we need the inequality
\begin{eqnarray}
| \sin (\rho) \, \hbox{ci}(\rho)
  + \cos (\rho) \, \hbox{si}(\rho) |  \le
  \hbox{const.} (1 + \rho)^{-1},  \qquad 
  0 <  \rho  < + \infty,
\label{eqn:4e}
\end{eqnarray}
which follows from the inequalities in the
subsections A.1 and A.2 in Appendix. 
The inequality (\ref{eqn:4e}), together with 
(\ref{eqn:3q}), immediately implies that
for each $\lambda > 0$, there is a 
positive constant $C_{\lambda}$ such that
\begin{eqnarray}
| m_{\lambda}(x) |  \le  C_{\lambda}
   \, |x|^{-1}   \, \langle x \rangle^{-1}.
\label{eqn:4e+}
\end{eqnarray}
It is apparent
 that one can take the constant
$C_{\lambda}$ in (\ref{eqn:4e+}) to be
uniform for $\lambda$ in each compact
interval in $(0, \, + \infty)$.

\vspace{15pt}

\begin{lem}\label{lem:kerm+1}
There exists a positive constant 
$C^{\prime}_{\lambda}$,  being uniform 
for $\lambda$ in each compact
interval in $(0, \, + \infty)$, such that
\begin{equation}
| M_{\lambda} \, u(x) | 
\le C^{\prime}_{\lambda} \, 
   \Vert u \Vert_{L^2}
\label{eqn:kerm+100}
\end{equation}
for all $u \in L^2({\mathbb R}^3)$.
\end{lem}

{\it Proof.} It follows from (\ref{eqn:4e+}) that
$m_{\lambda}\in L^{2}({\mathbb R}^3)$. 
Applying the Schwarz inequality to 
the right hand side of (\ref{eqn:4a+}) 
gives the lemma.
\hfill$\square$

\vspace{10pt}

\begin{lem}\label{lem:ker2}
Let $s> 3/2$.  Then there exists a constant
$C_{s\lambda}$ such that
\begin{eqnarray*}
| M_{\lambda} u(x) |  \le C_{s\lambda}   
( \langle x \rangle^{-2} + \langle x \rangle^{-s} )
  \, \Vert u \Vert_{L^{2,s}}
\end{eqnarray*}
for all $u \in L^{2, s}({\mathbb R}^3)$, 
$C_{s\lambda}$ being uniform 
for $\lambda$ in each compact
interval in $(0, \, + \infty)$.
\end{lem}

{\it Proof.} 
Let $u \in L^{2, s}({\mathbb R}^3)$.
We first note that $M_{\lambda}u(x)$ satisfies
 the inequality (\ref{eqn:kerm+100}),
since we can regard $u$ as an element in 
$L^2({\mathbb R}^3)$. 
Hence, we have
\begin{eqnarray}
| M_{\lambda}u(x) |  \le 
C^{\prime}_{\lambda}
   \Vert u \Vert_{L^{2,s}}.
\label{eqn:4e++}
\end{eqnarray}

We shall next show the inequality
\begin{eqnarray}
| M_{\lambda} u(x) |  \le C_{\lambda} 
\widetilde C_{s}  
  (  | x |^{-2} + \langle x \rangle^{-s} )
   \Vert u \Vert_{L^{2,s}},
\label{eqn:4e+++}
\end{eqnarray}
where $C_{\lambda}$ is the same constant as
in (\ref{eqn:4e+}) and
$\widetilde C_s$ is a constant depending only on $s$.
The inequality (\ref{eqn:4e+++}), together with
the inequality (\ref{eqn:4e++}), gives the lemma. 
In order to show (\ref{eqn:4e+++}), 
we decompose
$M_{\lambda}u(x)$ into three terms:
\begin{eqnarray}
M_{\lambda}u(x) =  
 I(x)  + I\!I(x)  + I\!I\!I(x), 
\label{eqn:4f}
\end{eqnarray}
where
\begin{eqnarray}
I(x) :=   
\int_{|y| \le |x|/2} m_{\lambda}(x-y) \, u(y) \, dy, 
\label{eqn:4f+}
\end{eqnarray}
\begin{eqnarray}
I\!I(x) := 
\int_{\scriptstyle \!\!\! \!\!\! |y| \ge |x|/2  
\atop \scriptstyle |x-y|\ge |x|/2} 
m_{\lambda}(x-y) \, u(y) \, dy, 
\label{eqn:4f++}
\end{eqnarray}
and
\begin{eqnarray}
I\!I\!I(x) := 
\int_{\scriptstyle \!\!\! \!\!\! |y| \ge |x|/2  
\atop \scriptstyle |x-y|\le |x|/2}   
m_{\lambda}(x-y) \, u(y) \, dy. 
\label{eqn:4f+++}
\end{eqnarray}
To deal with $I(x)$, we note 
that $|x - y| \ge |x| - |y| \ge |x|/2$
if
$|y| \le |x|/2$. This fact, together with (\ref{eqn:4e+}), yields
\begin{eqnarray}
|I(x)| &\le&  C_{\lambda}
\int_{|y| \le |x|/2} |x-y|^{-2} \, |u(y)| \, dy  \nonumber\\
\noalign{\vskip 2pt}
&\le& 4 \,C_{\lambda}  \,
|x|^{-2}   \! \int_{|y| \le |x|/2} |u(y)| \, dy \label{eqn:4g}\\
\noalign{\vskip 3pt}
&\le& 4\,C_{\lambda} \,C_s \,
|x|^{-2}  \, \Vert u \Vert_{L^{2,s}} ,  \nonumber
\end{eqnarray}
where we have used (\ref{eqn:4b}) in the last inequality
and the constant $C_s$ is the same one as in (\ref{eqn:4b}).  
It follows  from (\ref{eqn:4e+}) that
\begin{eqnarray}
|I\!I(x)| &\le& 
C_{\lambda} \int_{|x-y|\ge |x|/2} 
|x-y|^{-2} \, |u(y)| \, dy     \nonumber  \\
\noalign{\vskip 4pt}
&\le&  4\, C_{\lambda} \, C_s \, |x|^{-2}   
 \, \Vert u \Vert_{L^{2,s}}.
\label{eqn:4g+}
\end{eqnarray}
To get an estimate of $I\!I\!I(x)$, we should note that
if $|x -y| \le |x| /2$, 
then $|y| \ge |x| - |x -y| \ge |x|/2$, hence
$\langle y \rangle  \ge \langle x \rangle /2$. 
By using this fact and 
(\ref{eqn:4e+}), we have
\begin{eqnarray}
|I\!I\!I(x)| &\le&  C_{\lambda} \int_{|x-y|\le |x|/2} 
|x-y|^{-1} \langle x-y \rangle^{-1} \, |u(y)| \, dy     \nonumber\\
&\le&  C_{\lambda} 
\, \Big\{
\int_{|x-y|\le |x|/2} 
\frac{\langle y \rangle^{-2s} }
  {|x-y|^{2}\langle x-y \rangle^{2}} \, dy   \Big\}^{1/2}  
                    \Vert u \Vert_{L^{2,s}}               \nonumber\\
\noalign{\vskip 2pt}
&\le& 2^s \, C_{\lambda} \, 
   \langle x \rangle^{-s} \Vert u \Vert_{L^{2,s}},   
\label{eqn:4gg+} 
\end{eqnarray}
where we have used the Schwarz inequality in the second
inequality. Finally we deduce from
(\ref{eqn:4f}) -- (\ref{eqn:4gg+})  that  
(\ref{eqn:4e+++}) is verified.
\hfill$\square$

\vspace{15pt}

As an immediate corollary to Lemma \ref{lem:ker2}, 
we obtain a boundedness result
on the operator $M_{\lambda}.$

\begin{lem}\label{lem:ker3}
If $s> 3/2$, then $M_{\lambda} \in 
{\bf B}(L^{2, \, s}, \,L^2)$.
Moreover, the operator norm of $M_{\lambda}$
is bounded by a constant $C_{s\lambda}$,
which is uniform for $\lambda$ in
each compact interval in $(0, \, + \infty)$.
\end{lem}

\vspace{15pt}

We shall close this section with estimates of the operator
$K_{\lambda}^{\pm}$.

\begin{lem}\label{lem:ker4}
Let $s>1/2$. Then
there exists  a positive constant $C_{s}$ such
that
\begin{eqnarray*}
  | K_{\lambda}^{\pm} u(x)|                                       
\le C_{s} \, \lambda
\, \Vert u \Vert_{ L^{2, \, s} } 
\left\{  
\begin{array}{ll}
\langle x \rangle^{-(s -1/2)}  & \mbox{ if $\;  1/2  < s < 3/2$},  \\
{}  &  {}             \\
\langle x \rangle^{-1} \{ \log (1 + \langle x \rangle )\}^{1/2}   
       & \mbox{ if $\;   s = 3/2$},  \\
{}  & {}  \\
\langle x \rangle^{-1}  
       & \mbox{ if $\;   s > 3/2$}  
\end{array}
\right.
\end{eqnarray*}
for all $u \in L^{2, \, s}({\mathbb R}^3)$.
\end{lem}

{\it Proof.} Let $u \in L^{2, \, s}({\mathbb R}^3)$. 
Then applying the Schwarz inequality
to (\ref{eqn:4a}), 
we have
\begin{eqnarray}
| K_{\lambda}^{\pm}u(x)| 
\le \frac{\lambda}{\, 2 \pi \,}
 \Big\{ 
 \int_{{\mathbb R}^3}  
\frac{1}{\,|x-y|^2 \langle y \rangle^{2s} \,} \, dy 
   \Big\}^{\! 1/2} \,
   \Vert u \Vert_{L^{2,s}}.   \label{eqn:4h}
\end{eqnarray}  
We now apply Lemma A.1 in Appendix with
$n=3$, $\beta=2$ and $\gamma = 2 s >1$, and obtain the lemma.
\hfill$\square$

\vspace{15pt}
As an immediate consequence of Lemma \ref{lem:ker4}, 
we obtain a boundedness result on the operators
$K_{\lambda}^{\pm}$.

\begin{lem}\label{lem:ker5}
If $s>1$, then  $K_{\lambda}^{\pm} 
   \in {\bf B}(L^{2, \, s}, \,L^{2, \, -s})$.
Moreover, the operator norms of
$K_{\lambda}^{\pm}$ are bounded by 
$C_s \lambda$, where $C_s$ is a constant
depending only on $s$.
\end{lem}

\vspace{15pt}
Summing up all the results of Lemma \ref{lem:ker1}(ii)
and
Lemmas \ref{lem:kerm+1} and \ref{lem:ker4},
we see that (\ref{eqn:5ahlalah}) hold
on $L^{2, \, s}({\mathbb R}^3)$, $s > 1/2$, i.e.,
\begin{eqnarray}
R_0^{\pm}(\lambda) u=
G_{\lambda}^{\pm} u= 
G_0 u + K_{\lambda}^{\pm}u +  M_{\lambda}u
\label{eqn:5ahlalahAdd}
\end{eqnarray}
for all $\lambda >0$ and all 
$u \in L^{2, \, s}({\mathbb R}^3)$ 
with $s > 1/2$.

\newpage
\section{Radiation conditions for
\boldmath$\mathnormal{\sqrt{-\Delta}}$
\unboldmath}

This section is devoted to discussing radiation conditions 
for
$\sqrt{-\Delta}$ on ${\mathbb R}^3$.
The main result in this section is
Theorem \ref{thm:FRC}.

It is well-known that the radiation condition
$$
\Big(\frac{\partial}{\partial r} - 
  i \lambda \Big) u = O(r^{-2})   
  \; \mbox{ as $r=|x| \to \infty$}
$$
 was first introduced in order to
single out an outgoing solution of the Helmholtz equation
$(-\Delta - \lambda^2)u = f$ in ${\mathbb R}^3$, 
where $\lambda >0$.
The outgoing solution is the one which behaves as 
$e^{i \lambda r}/r$ at infinity.
In the present paper
 we shall exploit the Ikebe-Sait\={o}'s
formulation
of the radiation conditions for the Helmholtz equation,
which we regard as a special case of the 
time-independent Schr\"odinger
equations investigated in  
Ikebe-Sait\={o}\cite[Theorems 1.4, 1.5 and 
Remark 1.6]{IkebeSaito}. 
See also Sait\={o}\cite{Saito1}, \cite{Saito2}
and Pladdy-Sait\={o}-Umeda\cite{Pladdy}
for the formulation of the 
radiation conditions.

\vspace{15pt}

\begin{thm}[Ikebe-Sait\={o}]\label{thm:IS} 
Let $1/2 < s < 1$. 
\begin{description}
\item[\rm (i)] 
Suppose that 
$u$ belongs 
$L^{2, -s}({\mathbb R}^3)\cap 
H^2_{\rm loc}({\mathbb R}^3)$ 
and satisfies the equation
\begin{eqnarray}
(-\Delta - \lambda^2)u = 0  ,  
\qquad \lambda >0 ,
\label{eqn:5-0}
\end{eqnarray}
and, in addition, that $u$ satisfies either 
the outgoing radiation condition
\begin{eqnarray}
\Big( \frac{\partial}{\partial x_j} - 
i \lambda  \omega_j \Big) u \in 
L^{2, \, s-1}({\mathbb R}^3),
  \quad  j= 1, \, 2, \, 3 ,
\label{eqn:5-0a}
\end{eqnarray}
or the incoming radiation condition
\begin{eqnarray}
\Big( \frac{\partial}{\partial x_j} + 
i \lambda  \omega_j \Big) u 
\in L^{2, \, s-1}({\mathbb R}^3),
  \quad  j= 1, \, 2, \, 3 ,
\label{eqn:5-0b}
\end{eqnarray}
where  $\omega= x / |x|$.
Then $u$ vanishes identically.
\item[\rm (ii)] Suppose that 
$f \in L^{2, \, s}({\mathbb R}^3)$ 
and $\lambda >0$. Then 
$v^+(\lambda, \, f):= {\varGamma}^+_0(\lambda^2)f$ and
$v^-(\lambda, \, f):= {\varGamma}^-_0(\lambda^2)f$ satisfy 
the equation
\begin{eqnarray}
(-\Delta - \lambda^2)u = f  ,  
\qquad \lambda >0 
\label{eqn:5-0f}
\end{eqnarray}
 with 
the outgoing radiation condition {\rm (\ref{eqn:5-0a})} 
and 
the incoming radiation condition {\rm (\ref{eqn:5-0b})}
respectively. 
{\rm(}For the definition
of ${\varGamma}_0^{\pm}(z)$, 
see {\rm(\ref{eqn:3f})} and {\rm(\ref{eqn:3k})} {\rm).}
\end{description}
\end{thm}

\vspace{15pt}
It is not difficult to find radiation conditions
for $\sqrt{-\Delta}$ in a formal manner,
because it is easy to see that
$\sqrt{-\Delta}\,(\!\sqrt{-\Delta}u) = -\Delta u$
is formally valid. Actually a difficulty
arises if one tries to make sense of 
$\sqrt{-\Delta}\,(\!\sqrt{-\Delta}u)$ 
for $u \in L^{2, -s}({\mathbb R}^3)$ with $s<0$. 
The difficulty comes from the fact that
the symbol $|\xi|$ is singular at the 
origin $\xi = 0$ (see
 Lieb-Loss\cite[\S 7.15]{LiebLoss}).
In order to overcome the difficulty, we
need to clarify the function spaces to which
$\sqrt{-\Delta}u$ belongs when $u$ belongs
to $L^{2, -s}({\mathbb R}^3)$ with $s<0$. 
The clarification was  made in
Umeda\cite{Umeda2}.
We reproduce \cite[Theorem 5.8]{Umeda2}
for the reader's convenience.

\vspace{15pt}

\begin{thm}[Umeda]\label{thm:U3-5} 
Let  $\ell \in \mathbb R$. 
If $s$ and $t$ satisfy either
\begin{eqnarray}
 s \ge 0,  \;\;\; t < \min \{1, \, s -3/2\}
\end{eqnarray}
or
\begin{eqnarray}
-5/2 < s <0, \;\;\;  t < s -3/2,
\end{eqnarray}
 then
$\sqrt{-\Delta}$ is a bounded operator from
$H^{\ell, \, s}({\mathbb R}^3)$ to 
$H^{\ell -1, \, t}({\mathbb R}^3)$.
\end{thm}

\vspace{15pt}

With the aid of Theorem \ref{thm:U3-5}, we shall 
first 
make sense  of $\sqrt{-\Delta}\,(\!\sqrt{-\Delta}u)$ 
for $u \in {\cal S}({\mathbb R}^3)$.

\vspace{15pt}

\begin{lem}\label{lem:5-1}
If $\varphi \in  {\cal S}({\mathbb R}^3)$, then 
$\sqrt{-\Delta}\,(\!\sqrt{-\Delta}\varphi)\in
H^{-1, t}({\mathbb R}^3)$ for
all $t < 1$, and
\begin{eqnarray}
\sqrt{-\Delta}\,(\!\sqrt{-\Delta} \, \varphi) 
= -\Delta \varphi  
  \quad    in  \;\; {\cal S}^{\prime}({\mathbb R}^3).
\label{eqn:5-1}
\end{eqnarray}
\end{lem}

{\it Proof.} Let $\varphi \in  {\cal S}({\mathbb R}^3)$.
By virtue of \cite[Theorem 4.4]{Umeda2},
 we find 
that $\sqrt{-\Delta} \, \varphi 
\in L^{2, s}({\mathbb R}^3)$
for any $s < 5/2$. It follows from
Theorem \ref{thm:U3-5} that 
$\sqrt{-\Delta}\,(\!\sqrt{-\Delta}\varphi)$
makes sense, and that 
$\sqrt{-\Delta}\,(\!\sqrt{-\Delta}\varphi)$ belongs to
$H^{-1, t}({\mathbb R}^3)$ for all $t < 1$.
It follows, in particular, that
$\sqrt{-\Delta}\,(\!\sqrt{-\Delta}\varphi) 
\in {\cal S}^{\prime}({\mathbb R}^3)$.

To prove (\ref{eqn:5-1}), we take a test function
$\psi \in {\cal S}({\mathbb R}^3)$. By definition of 
the action
of $\sqrt{-\Delta}$ on distributions we have
\begin{eqnarray}
\langle 
\sqrt{-\Delta}\,(\!\sqrt{-\Delta}\varphi), \, 
\psi \rangle = 
( \sqrt{-\Delta}\varphi, \, 
{ \sqrt{-\Delta} \; 
\overline{\psi}} \, )_{-s, s}
\label{eqn:5-1a}
\end{eqnarray}
if $-5/2 < s <5/2$.
(It follows from
\cite[Theorem 4.4]{Umeda2} that the mapping
$\psi \mapsto  
(\sqrt{-\Delta}\varphi, \, 
{ \sqrt{-\Delta} \; 
\overline{\psi}} \, )_{-s, s}$ 
 is a continuous linear functional 
on ${\cal S}({\mathbb R}^3)$,
because one can regard $\sqrt{-\Delta}\varphi$ as 
a function belonging to $L^{2, -s}({\mathbb R}^3)$ for 
any  $s > -5/2$, and because one  finds that
$$
|(
\sqrt{-\Delta}\varphi, \, 
{ \sqrt{-\Delta} \; 
\overline{\psi}} \, )_{-s, s} |\le 
\Vert \sqrt{-\Delta}\varphi \Vert_{L^{2, -s}}
\Vert \sqrt{-\Delta} \,
\overline{\psi} \Vert_{L^{2, s}}
$$
for any $s$ with $-5/2 < s <5/2$.)
It is clear that we can regard the right hand side
of (\ref{eqn:5-1a}) as the inner product in 
$L^2({\mathbb R}^3)$, and we have
\begin{eqnarray*}
\langle 
\sqrt{-\Delta}\,(\!\sqrt{-\Delta}\varphi), \, 
\psi \rangle 
&=& ( \sqrt{-\Delta}\varphi, \, 
{ \sqrt{-\Delta} \; \overline{\psi}} )_{L^2}   \\
&=&
( |\xi| \, {\cal F}[\,\varphi\,], \, 
|\xi| \,{\cal F}[\,{\overline{\psi}} \,])_{L^2_{\xi}} \\
&=&
(-\Delta \varphi, \, \overline{\psi})_{L^2}   \\
&=&
\langle 
-\Delta \varphi, \, 
\psi \rangle ,    
\end{eqnarray*}
where we have used the Plancherel formula
twice. This proves (\ref{eqn:5-1}).
\hfill$\square$

\vspace{15pt}

\begin{lem}\label{lem:5-2}
Suppose that $\ell \in \mathbb R$ and $0< s < 1$.
If $u \in H^{\ell, -s}({\mathbb R}^3)$, then
\begin{eqnarray}
\sqrt{-\Delta}\,(\!\sqrt{-\Delta} \, u) 
= -\Delta u  
  \quad    in  \;\; {\cal S}^{\prime}({\mathbb R}^3).
\label{eqn:5-2}
\end{eqnarray}
\end{lem}

{\it Proof.}
Let $u$ be in $H^{\ell, -s}({\mathbb R}^3)$. 
Since ${\cal S}({\mathbb R}^3)$ is dense in 
$H^{\ell, -s}({\mathbb R}^3)$, we can
choose a sequence $\{ \varphi_j \} 
\subset {\cal S}({\mathbb R}^3)$ so that
$\varphi_j \to u$ in $H^{\ell, -s}({\mathbb R}^3)$ 
as $j \to \infty$. 
It follows from Theorem \ref{thm:U3-5}  that
\begin{eqnarray}
\sqrt{-\Delta} \varphi_j  \to \sqrt{-\Delta} u
 \quad \mbox{ in }\;  H^{\ell-1, \,t}({\mathbb R}^3)
\label{eqn:5-2a}
\end{eqnarray}
for any $t < -s -3/2$. In view of the hypothesis 
that $0 < s <1$,
we can find that (\ref{eqn:5-2a}) holds for any $t$ 
satisfying 
$-5/2  < t < -s - 3/2$.
Therefore, it follows from Theorem \ref{thm:U3-5} again
that
\begin{eqnarray}
\sqrt{-\Delta}\,(\!\sqrt{-\Delta} \varphi_j )
 \to \sqrt{-\Delta}\,(\!\sqrt{-\Delta} u)
   \quad \mbox{ in }\;  H^{\ell-2, \,t}({\mathbb R}^3)
\label{eqn:5-2b}
\end{eqnarray}
for any $t < -s -3$. In particular, we have
\begin{eqnarray}
-\Delta \varphi_j 
 \to \sqrt{-\Delta}\,(\!\sqrt{-\Delta} u)
   \quad \mbox{ in } \; {\cal S}^{\prime}({\mathbb R}^3),
\label{eqn:5-2c}
\end{eqnarray}
where we have used Lemma \ref{lem:5-1}.
On the other hand, by using the fact that
$\varphi_j \to u$ in $H^{\ell, -s}({\mathbb R}^3)$,
we obtain 
\begin{eqnarray}
-\Delta \varphi_j 
 \to {-\Delta} u
 \quad \mbox{ in } \; {\cal S}^{\prime}({\mathbb R}^3).
\label{eqn:5-2d}
\end{eqnarray}
Combining (\ref{eqn:5-2c}) with (\ref{eqn:5-2d}) gives
(\ref{eqn:5-2}).
\hfill$\square$

\vspace{20pt}

We shall now establish the
radiation conditions for $\sqrt{-\Delta}$ in the 
same formulation as in Theorem \ref{thm:IS}.

\vspace{10pt}

\begin{thm} \label{thm:FRC}
Let $1/2 < s < 1$.
\begin{description}
\item[\rm (i)] 
Suppose that $u$ belongs to
$L^{2, -s}({\mathbb R}^3)\cap H^1_{\rm loc}({\mathbb R}^3)$ 
and satisfies the equation
\begin{eqnarray}
(\sqrt{-\Delta} - \lambda) u = 0 
           \quad in  \;\; {\cal S}^{\prime}({\mathbb R}^3), 
  \qquad \lambda > 0,
\label{eqn:Tfrc}
\end{eqnarray}
and, in addition, that $u$ satisfies either of
the outgoing radiation condition
{\rm(\ref{eqn:5-0a})} 
or 
the incoming radiation condition 
{\rm(\ref{eqn:5-0b})}.
Then $u$ vanishes identically.
\item[\rm (ii)]
Suppose that 
$f \in L^{2, \, s}({\mathbb R}^3)$ 
and $\lambda >0$. Then 
$u^+_0(\lambda, \, f):= R^+_0(\lambda)f$ and
$u^-_0(\lambda, \, f):= R^-_0(\lambda)f$ satisfy 
the equation 
\begin{eqnarray}
(\sqrt{-\Delta} - \lambda ) 
  u = f 
\quad \hbox{ in } \; 
   {\cal S}^{\prime}({\mathbb R}^3)
\label{eqn:6-RC+}
\end{eqnarray}
with 
the outgoing radiation condition {\rm (\ref{eqn:5-0a})} 
and 
the incoming radiation condition {\rm (\ref{eqn:5-0b})}
respectively. 
\end{description}
\end{thm}

\vspace{15pt}
A very important consequence of Theorem \ref{thm:FRC}
is the fact that the 
radiation conditions (\ref{eqn:5-0a}) and 
(\ref{eqn:5-0b})  characterize the boundary 
values $R^+_0(\lambda)$ and
$R^-_0(\lambda)$ respectively.

In order to prove Theorem \ref{thm:FRC}, we need to
prepare two lemmas.
One might regard
the equality (\ref{eqn:3-4a}) below
as straightforward.
Unfortunately, this is not the case.
Indeed, there exists a diffuculty 
to make sense of
$\sqrt{-\Delta} R^{\pm}_0(\lambda)f$.
The reason for this difficulty is the same
as the ones mentioned before Theorem \ref{thm:U3-5}, 
namely, the fact that
$R^{\pm}_0(\lambda)f$ merely belong to
$L^{2,-s}({\mathbb R}^3)$
with $s > 1/2$.
Nevertheless we can prove, with the 
aid of  theorems in Umeda\cite{Umeda2}, that
(\ref{eqn:3-4a}) is true.

\vspace{15pt}

\begin{lem}\label{lem:3-4}
Suppose that
 $\lambda >0$ and
$f \in L^{2,s}({\mathbb R}^3)$, $s >1/2$. 
Then
\begin{eqnarray}
(\sqrt{-\Delta} - \lambda ) 
  R^{\pm}_0 (\lambda) f = f 
\quad \hbox{ in } \; 
   {\cal S}^{\prime}({\mathbb R}^3).
\label{eqn:3-4a}
\end{eqnarray}
\end{lem}

{\it Proof.} We can assume, without 
loss of generality, that
$1/2 < s < 5/2$. It then follows from
Theorem \ref{thm:U3-5}
(cf. \cite[Theorem 4.6]{Umeda2}) that 
$\sqrt{-\Delta} R^{\pm}_0(\lambda)f$ make
sense.
In order to show (\ref{eqn:3-4a}), we take a test
function 
$\psi \in {\cal S}({\mathbb R}^3)$.
We then have
\begin{eqnarray}
\langle \, (\sqrt{-\Delta} - \lambda \mp i \mu) 
  R_0 (\lambda \pm i \mu) f , \, \psi \rangle
 =  \langle f, \, \psi \rangle
\label{eqn:3-4b}
\end{eqnarray}
for all $\mu >0$, since 
$R_0 (\lambda \pm i \mu) f$ belong to
$H^1({\mathbb R}^3)$, the domain of the selfadjoint
operator $H_0$,
and since
\begin{eqnarray*}
(\sqrt{-\Delta} - \lambda \mp i \mu) 
  R_0 (\lambda \pm i \mu) f & =&
(H_0 - (\lambda \pm i \mu) )
  R_0 (\lambda \pm i \mu) f   \\
&=& f.
\end{eqnarray*}
By definition of the action of $\sqrt{-\Delta}$
on $L^{2, \, -s}({\mathbb R}^3)$, the left hand side
of (\ref{eqn:3-4b}) becomes
\begin{eqnarray}
(   R_0 (\lambda \pm i \mu) f , \, 
 {\sqrt{-\Delta} \, \overline{\psi} }
\, )_{-s, s}
 - ( R_0 (\lambda \pm i \mu)f, 
  \,  (\lambda \mp i \mu) 
\overline{\psi} )_{-s, s}.
\label{eqn:3-4c}
\end{eqnarray}
(Note that 
$\sqrt{-\Delta}\,\overline{\psi} 
\in L^{2, \, t}({\mathbb R}^3)$ for any 
$t < 5/2$; see \cite[Theorem 4.4]{Umeda2}.)
It follows from Theorem \ref{thm:BA-N} that
\begin{eqnarray}
\lim_{\mu \downarrow 0} \, 
(  R_0 (\lambda \pm i \mu) f , \, 
 {\sqrt{-\Delta} \, \overline{\psi} } \,
)_{-s, s}
=
(  R_0^{\pm} (\lambda) f , \, 
 {\sqrt{-\Delta} \, \overline{\psi} } \,
)_{-s, s}
\label{eqn:3-4d}
\end{eqnarray}
Combining (\ref{eqn:3-4c}), (\ref{eqn:3-4d})
with (\ref{eqn:3-4b}), we conclude that
\begin{eqnarray*}
(  R_0^{\pm} (\lambda) f , \, 
 {\sqrt{-\Delta} \, \overline{\psi} } \,
)_{-s, s}
-
(  R_0^{\pm} (\lambda) f , \, 
   \lambda \overline{\psi}  
)_{-s, s}
= \langle   
   f , \, \psi  
\rangle
\end{eqnarray*}
for any test function 
$\psi \in {\cal S}({\mathbb R}^3)$. 
This completes the proof.
\hfill$\square$

\vspace{15pt}

\begin{lem}\label{lem:5-3}
Suppose that $1/2 < s < 1$ and 
$\lambda >0$.
If $u$ belongs to\linebreak 
$\hbox{\rm Ran}\big(R_0^+(\lambda)\big)$, 
then $u$
satisfies the the outgoing radiation condition
{\rm(\ref{eqn:5-0a})}. Similarly, if $u$ belongs to
$\hbox{\rm Ran}\big(R_0^-(\lambda)\big)$, 
then $u$ satisfies 
the incoming radiation condition
{\rm(\ref{eqn:5-0b})}.
\end{lem}

{\it Proof.}
We only give the proof 
for $u \in\hbox{\rm Ran}\big(R_0^+(\lambda)\big)$.
The proof for  
$u \in\hbox{\rm Ran}\big(R_0^-(\lambda)\big)$
is similar.

By assumption, one can 
find an $f \in L^{2, \,s}({\mathbb R}^3)$ such that
$u =R_0^+(\lambda) f$. It follows from
Theorem \ref{thm:lap0}, together with
Corollary to Lemma 3.3, that
there exist $A(\lambda) \in {\bf B}(L^{2,s})$,
$B_1(\lambda) \in {\bf B}(L^{2,s}, \; H^{1,s})$ 
and
$B_2(\lambda) \in {\bf B}(L^2, \; H^2)$ 
such that
\begin{eqnarray}
u = {\varGamma}_0^+(\lambda^2) \, A(\lambda)f 
  + B_1(\lambda)f   +  B_2(\lambda)f.
\label{eqn:rc+1}
\end{eqnarray}
By Theorem \ref{thm:IS}(ii), 
the first term on the right hand side
of (\ref{eqn:rc+1}) satisfies 
the outgoing radiation condition
{\rm(\ref{eqn:5-0a})}.
Since $B_1(\lambda)f \in H^{1, \, s}({\mathbb R}^3)$,
it is straightforward to see that
\begin{eqnarray*}
\Big( \frac{\partial}{\partial x_j} - 
i \lambda  \omega_j \Big) B_1(\lambda) f \in 
L^{2, \,s}({\mathbb R}^3) 
\subset L^{2, \, s-1}({\mathbb R}^3),
  \quad  j= 1, \, 2, \, 3,
\end{eqnarray*}
that is, the second term on the right hand side
of (\ref{eqn:rc+1}) satisfies 
{\rm(\ref{eqn:5-0a})}.
Finally, since $B_2(\lambda)f \in H^2({\mathbb R}^3)$, 
it 
follows that
\begin{eqnarray*}
\Big( \frac{\partial}{\partial x_j} - 
i \lambda  \omega_j \Big) B_2(\lambda) f \in 
H^1({\mathbb R}^3) \subset L^{2, \, s-1}({\mathbb R}^3),
  \quad  j= 1, \, 2, \, 3,
\end{eqnarray*}
where we have used the assumption that $s<1$.
Hence the last term on the right hand side
of (\ref{eqn:rc+1}) satisfies 
{\rm(\ref{eqn:5-0a})}. 
\hfill$\square$

\vspace{20pt}

{\bf Proof of Theorem \ref{thm:FRC}} 
\ It follows from (\ref{eqn:Tfrc}) 
that $\sqrt{-\Delta}u = \lambda u$,
hence $\sqrt{-\Delta}u$ belongs 
to $L^{2, -s}({\mathbb R}^3)$ with
$1/2 <s < 1$.
By Lemma \ref{lem:5-2}, 
it makes sense to consider
$\sqrt{-\Delta}\,(\!\sqrt{-\Delta}u)$, and we see that
$u$ satisfies
\begin{eqnarray*}
( -\Delta - \lambda^2 )u = 0 
     \quad    in  \;\; 
    {\cal S}^{\prime}({\mathbb R}^3),
\end{eqnarray*}
which implies that
$-\Delta u = \lambda^2 u$ belongs to 
$L^2_{\rm loc}({\mathbb R}^3)$. Therefore,
 we find that 
$u \in H^2_{\rm loc}({\mathbb R}^3)$. 
It is  evident
that we can apply Theorem \ref{thm:IS}(i)
 and obtain assertion (i) of the theorem.

Assertion (ii) of the theorem is an
 immediate consequence of Lemmas \ref{lem:3-4}
and \ref{lem:5-3}.
\hfill$\square$

\newpage
\section{Radiation conditions for
\boldmath$\mathnormal{\sqrt{-\Delta}+ V}$
\unboldmath}

This section is devoted to discussing radiation 
conditions 
for
$\sqrt{-\Delta} + V$ on ${\mathbb R}^3$.
As mentioned in Introduction, we assume that
$V(x)$ is a real-valued measurable function on
${\mathbb R}^3$ satisfying (\ref{eqn:Vdecay}).
 Under this assumption,
it is obvious that $V= V(x) \times$ is a bounded
selfadjoint operator in $L^2({\mathbb R}^3)$, and that
$H:= H_0 + V$ defines a selfadjoint operator
in $L^2({\mathbb R}^3)$, of which domain is 
$H^1({\mathbb R}^3).$
For $z \in \rho(H)$,  we write
\begin{eqnarray*}
R(z) = (H-z)^{-1}.
\end{eqnarray*}
It is clear that $H$ is essentially selfadjoint
on $C_0^{\infty}({\mathbb R}^3)$, 
since $H$ is a bounded  selfadjoint
perturbation of $H_0$, which is 
essentially selfadjoint
on $C_0^{\infty}({\mathbb R}^3)$ (see Section 2).
Since $V$ is relatively 
compact with respect to $H_0$, it follows
from 
Reed-Simon\cite[p.113, Corollary 2]{ReedSimon4}
that
\begin{eqnarray*}
\sigma_{\rm ess}(H) = \sigma_{\rm ess}(H_0) = 
 [0, \, +\infty).
\end{eqnarray*}

Before establishing the radiation conditions 
for $\sqrt{-\Delta}+V(x)$,
we need to remark that 
$\sigma_{\rm p}(H)\cap(0, \, + \infty)$ is a discrete set.
 This fact 
was first proved by Simon\cite[Theorem 2.1]{Simon} 
in a general setting,
and later recovered by 
Ben-Artzi and Nemirovsky\cite[Theorem 4A]{BenArtzi-Nem} 
also in a general setting.
Moreover, Simon\cite[Theorem 2.1]{Simon}
proved that
each eigenvalue in the set
$\sigma_{\rm p}(H)\cap(0, \, + \infty)$ 
has finite multiplicity.

To formulate the main theorem in this section, we exploit 
a result, which is a special case of
 Ben-Artzi and Nemirovsky\cite[Theorem 4A]{BenArtzi-Nem}.

\vspace{15pt}

\begin{thm}[Ben-Artzi and Nemirovski] \label{thm:BA-N2}
Let $\sigma > 1$ and $s > 1/2$. Then 
\begin{description}
\item[\rm (i)] The continuous spectrum 
$\sigma_{\rm c}(H)=[0, \, +\infty)$ is absolutely continuous, 
except possibly for a discrete set of embedded eigenvalues
$\sigma_{\rm p}(H) \cap (0, \, + \infty)$, which can 
accumulate only at $0$ and $+\infty$.

\item[\rm (ii)] For any 
$\lambda \in (0, \, +\infty) \setminus 
\sigma_{\rm p}(H)$, there
exist the limits
$$
R^{\pm}(\lambda) = \lim_{\mu \downarrow 0} 
   R(\lambda \pm i \mu)   \quad
  {\it in } \; {\bf B}(L^{2,\,s}, \; H^{1,\,-s}).
$$
\item[\rm (iii)] The operator-valued 
functions $R^{\pm}(z)$ defined by
\begin{equation*}
R^{\pm}(z) = 
\begin{cases}  
R(z) &\text{\it if $\;  z \in {\mathbb C}^{\pm}$}  \\
{}    &\\
R^{\pm}(\lambda)  &\text{\it if $\;  
z=\lambda \in (0, \, +\infty) \setminus
\sigma_{\rm p}(H) $}
\end{cases}
\end{equation*}
are ${\bf B}(L^{2,\,s}, \; H^{1,\,-s})$-valued 
continuous functions.
\end{description}
\end{thm}

 We now state the main result in this section, which
establishes
the radiation conditions 
for $\sqrt{-\Delta}+V(x)$.

\vspace{5pt}

\begin{thm} \label{thm:RC}     
Let $\sigma >1$
 and $1/2 < s < \min\, (\sigma/2, \, 1)$. 
\begin{description}
\item[\rm (i)] 
Suppose that $u$ belongs to
$L^{2, -s}({\mathbb R}^3)\cap 
 H^1_{\rm loc}({\mathbb R}^3)$ 
and satisfies the equation
\begin{eqnarray}
(\sqrt{-\Delta} + V(x) - \lambda) u = 0 
           \quad in  \;\; {\cal S}^{\prime}({\mathbb R}^3), 
 \; \quad \lambda \in (0, \, + \infty) 
\setminus \sigma_{\rm p}(H)
\label{eqn:rc+}
\end{eqnarray}
and, in addition, that $u$ satisfies either of
the outgoing radiation condition
{\rm(\ref{eqn:5-0a})} 
or 
the incoming radiation condition 
{\rm(\ref{eqn:5-0b})}.
Then $u$ vanishes identically.
\item[\rm (ii)]
Suppose that 
$f \in L^{2, \, s}({\mathbb R}^3)$ 
and $\lambda \in (0, \, + \infty) 
\setminus \sigma_{\rm p}(H)$. Then 
$u^+(\lambda, \, f):= R^+(\lambda)f$ and
$u^-(\lambda, \, f):= R^-(\lambda)f$ satisfy 
the equation 
\begin{eqnarray}
(\sqrt{-\Delta} + V(x) - \lambda ) 
  u = f 
\quad \hbox{ in } \; 
   {\cal S}^{\prime}({\mathbb R}^3)
\label{eqn:6-RC+v}
\end{eqnarray}
with 
the outgoing radiation condition {\rm (\ref{eqn:5-0a})} 
and 
the incoming radiation condition {\rm (\ref{eqn:5-0b})}
respectively. 
\end{description}
\end{thm}

\vspace{15pt}

The same remark after Theorem \ref{thm:FRC}
applies to Theorem \ref{thm:RC}, namely,
Theorem \ref{thm:RC} gives the characterization
of the boundary 
values $R^+(\lambda)$ and
$R^-(\lambda)$ in terms of the 
radiation conditions (\ref{eqn:5-0a}) and 
(\ref{eqn:5-0b}) respectively. 

We shall give a proof of Theorem \ref{thm:RC}
by means of a series of lemmas, but
only for $u$ satisfying
the outgoing radiation condition
{\rm(\ref{eqn:5-0a})}.
The proof for $u$ satisfying
the incoming radiation condition
{\rm(\ref{eqn:5-0b})} is similar.

\vspace{15pt}
\begin{lem}\label{lem:5-4}
Let $\sigma >1$, and
suppose that $1/2 < s < \sigma/2$. Then
\begin{eqnarray}
R^{\pm}(z) \big( I + V  R^{\pm}_0(z) \big) 
&=& R^{\pm}_0(z)  
    \quad \hbox{\rm on } L^{2, \, s}({\mathbb R}^3),
\label{eqn:5-4a}   \\
\noalign{\vskip 4pt}
R^{\pm}_0(z) \big( I - V  R^{\pm}(z) \big)  
&=& R^{\pm}(z)  
   \quad \hbox{\rm on } L^{2, \, s}({\mathbb R}^3)
\label{eqn:5-4b}
\end{eqnarray}
for all 
$z \in {\mathbb C}^{\pm}\cup 
   \{ (0, \, +\infty) \setminus \sigma_p(H) \}$.
\end{lem}

{\it Proof.} We shall give the proof only
in the case where the superscripts are \lq\lq +," 
the plus sign.
If $z \in {\mathbb C}^+$, it is apparent
that
\begin{eqnarray*}
( H - z) R_0(z) &=& I + VR_0(z)  
    \quad \hbox{\rm on } L^{2}({\mathbb R}^3),   \\
( H_0 - z) R(z) &=& I - VR(z)  
    \quad \hbox{\rm on } L^{2}({\mathbb R}^3), 
\end{eqnarray*}
from which it follows that
\begin{eqnarray}
R(z) \big( I + V  R_0(z) \big) 
&=& R_0(z)  \quad \hbox{\rm on } L^2({\mathbb R}^3),
\label{eqn:5-4c}   \\
\noalign{\vskip 4pt}
R_0(z) \big( I - V  R(z) \big)  
&=& R(z)  \quad \hbox{\rm on } L^2({\mathbb R}^3).
\label{eqn:5-4d}
\end{eqnarray}

In order to proceed to the extended resolvents,
we now regard that
$R^+_0(z)$ and $R^+(z)$ are
${\bf B}(L^{2,\,s}, L^{2,\,-s} )$-valued 
continuous functions
on
${\mathbb C}^+ \cup (0, \, +\infty)$
and 
${\mathbb C}^+\cup 
   \{ (0, \, +\infty) \setminus \sigma_p(H) \}$
respectively. 
By (\ref{eqn:Vdecay}), and by
the assumption that $1/2 < s < \sigma/2$, we see 
that $V \in {\bf B}(L^{2, -s}, L^{2, s})$, and hence
$VR^+_0(z)$ and $VR^+(z)$ are
${\bf B}(L^{2,\,s})$-valued continuous functions
on
${\mathbb C}^+ \cup (0, \, +\infty)$
and 
${\mathbb C}^+\cup 
   \{ (0, \, +\infty) \setminus \sigma_p(H) \}$
respectively. 
Therefore, 
we conclude from (\ref{eqn:5-4c}) and 
(\ref{eqn:5-4d}) that  the assertion 
of the lemma is valid.
\hfill$\square$

\vspace{20pt}

As a corollary to Lemma \ref{lem:5-4}, we
obtain the following result.

\vspace{15pt}

\begin{lem}\label{lem:5-5}
Let $\sigma >1$, and
suppose that $1/2 < s < \sigma/2$. Then
\begin{eqnarray*}
{\rm Ran}\big( R^{\pm}(z) \big) = 
   {\rm Ran}\big( R_0^{\pm}(z) \big)
\end{eqnarray*}
for every $z \in {\mathbb C}^{\pm}\cup 
 \{ (0, \, +\infty) \setminus \sigma_p(H) \}$.
\end{lem}

\vspace{15pt}
\begin{lem}\label{lem:5-6}
Let $\sigma >1$, and
suppose that $1/2 < s < \sigma/2$. 
Then
\begin{eqnarray}
\big(I - R^{\pm}(z) V \big) 
\big( I +  R^{\pm}_0(z) V  \big) 
&=& I 
  \quad \hbox{\rm on } L^{2, \, -s}({\mathbb R}^3),
\label{eqn:5-6a}   \\
\noalign{\vskip 4pt}
\big(I + R^{\pm}_0(z)V \big) 
 \big( I - R^{\pm}(z) V   \big)  
&=& I 
  \quad \hbox{\rm on } L^{2, \, -s}({\mathbb R}^3)
\label{eqn:5-6b}
\end{eqnarray}
for every 
$z \in {\mathbb C}^{\pm}\cup 
   \{ (0, \, +\infty) \setminus \sigma_p(H) \}$.
\end{lem}

{\it Proof.} We shall only give the proof 
of (\ref{eqn:5-6a})
in the case where the superscripts are \lq\lq +." 
The proof 
of (\ref{eqn:5-6a})
in the other case and the proof of
(\ref{eqn:5-6b}) are similar.

We first show that for every $z \in {\mathbb C}^+$
\begin{eqnarray}
\big( I - R(z) V \big) \, 
 \big( I + R_0(z) V \big) = I 
   \quad \hbox{\rm on }  L^2({\mathbb R}^3).
\label{eqn:5-6c}
\end{eqnarray}
In fact, if 
$u$ belongs to 
$H^1({\mathbb R}^3)$,
 we then have
\begin{eqnarray*}
( I - R(z) V) u &=& R(z) (H - z)u - R(z) Vu   \\
  &=&   R(z) (H_0 - z)u
\end{eqnarray*}
and
\begin{eqnarray*}
( I + R_0(z) V) u &=& R_0(z) (H_0 - z)u + R_0(z) Vu   \\
  &=&   R_0(z) (H - z)u
\end{eqnarray*}
(recall that  
$\hbox{\rm Dom}(H) = \hbox{\rm Dom}(H_0)
= H^1({\mathbb R}^3)$).
Hence we get
\begin{eqnarray*}
\big( I - R(z) V \big)  
\big( I + R_0(z) V \big) u
&=& R(z) (H_0 - z)  R_0(z) (H - z)u   \\ 
&=& u
\end{eqnarray*}
for all $u \in H^1({\mathbb R}^3)$, where
we have used the fact that
$\big( I + R_0(z) V \big) u \in H^1({\mathbb R}^3)$ 
when $u \in H^1({\mathbb R}^3)$.
Since $H^1({\mathbb R}^3)$ is dense in 
  $L^2({\mathbb R}^3)$, we can deduce that
(\ref{eqn:5-6c}) is true.

We next work in the weighted $L^2$-spaces.
As mentioned in the proof of 
Lemma \ref{lem:5-4}, we have
$V \in {\bf B}(L^{2, \, - s}, \,L^{2, \, s})$.
Also, as mentioned in the second half
 of the proof of
Lemma \ref{lem:5-4}, 
we can regard that
$R^+_0(z)$ and $R^+(z)$ are
${\bf B}(L^{2,\,s}, L^{2,\,-s} )$-valued 
continuous functions
on
${\mathbb C}^+ \cup (0, \, +\infty)$
and 
${\mathbb C}^+\cup 
   \{ (0, \, +\infty) \setminus \sigma_p(H) \}$
respectively.
Therefore 
$R^+_0(z)V$ and $R^+(z)V$ are
${\bf B}(L^{2, \, -s})$-valued 
continuous functions
on
${\mathbb C}^+ \cup (0, \, +\infty)$
and 
${\mathbb C}^+\cup 
   \{ (0, \, +\infty) \setminus \sigma_p(H) \}$
respectively.
Thus, we can conclude from (\ref{eqn:5-6c}) 
that
(\ref{eqn:5-6a}) in the case where
the superscripts are
 the plus sign is true. 
\hfill$\square$

\vspace{15pt}

{\bf Proof of Theorem \ref{thm:RC}}\ \ 
We first prove assertion (i) of the theorem.
Let $u$ belong to
$L^{2, -s}({\mathbb R}^3)\cap H^1_{\rm loc}({\mathbb R}^3)$ 
and satisfy the equation (\ref{eqn:rc+})
together with
the outgoing radiation condition
{\rm(\ref{eqn:5-0a})}.
By (\ref{eqn:rc+}), we have 
\begin{eqnarray}
(\sqrt{-\Delta} - \lambda ) u = -Vu
\quad \hbox{ in } \;{\cal S}^{\prime}({\mathbb R}^3).
\label{eqn:7p1}
\end{eqnarray}
Since 
$Vu$ belongs to $L^{2, \, s}({\mathbb R}^3)$ 
by the fact that 
$V \in {\bf B}(L^{2, \, - s}, \,L^{2, \, s})$,
it  follows from Lemma \ref{lem:3-4}
that
\begin{eqnarray}
(\sqrt{-\Delta} - \lambda ) 
  R^+_0 (\lambda) Vu = Vu 
\quad \hbox{ in } \; 
   {\cal S}^{\prime}({\mathbb R}^3).
\label{eqn:7p2}
\end{eqnarray}
Combining (\ref{eqn:7p1}) with (\ref{eqn:7p2}) gives
\begin{eqnarray}
(\sqrt{-\Delta} - \lambda ) 
  (u + R^+_0 (\lambda)Vu )= 0
\quad \hbox{ in } \; 
   {\cal S}^{\prime}({\mathbb R}^3).
\label{eqn:7p3}
\end{eqnarray}
By virtue of Lemma \ref{lem:5-3}  and the 
fact that 
$R^+_0(\lambda)Vu \in H^{1, -s}({\mathbb R}^3)$,
it follows that  
$u + R^+_0 (\lambda)Vu$ belongs to
$L^{2, -s}({\mathbb R}^3)\cap H^1_{\rm loc}({\mathbb R}^3)$ 
and satisfies the outgoing radiation condition 
(\ref{eqn:5-0a}).
Hence we can apply Theorem \ref{thm:FRC} and conclude 
that
\begin{eqnarray}
u+ R^+_0 (\lambda)Vu = 0.
\label{eqn:7p4}
\end{eqnarray}
Since $u$ belongs to
$L^{2, -s}({\mathbb R}^3)$, it follows
from (\ref{eqn:7p4}) and Lemma \ref{lem:5-6}
that $u$ vanishes identically.

We next prove assertion (ii). 
It follows from Lemmas \ref{lem:5-5} and \ref{lem:5-3}
that $u^+(\lambda, \, f)$ satisfies the outgoing 
radiation condition (\ref{eqn:5-0a}). 
In order to show that $u^+(\lambda, \, f)$ is 
a solution to the
equation (\ref{eqn:6-RC+v}),
we follow the idea exploited in the proof of 
Lemma \ref{lem:3-4}.
Thus we start with
\begin{eqnarray*}
(\sqrt{-\Delta} + V - \lambda - i \mu) 
  R (\lambda + i \mu) f  = f,  \quad  \forall \mu >0,
\end{eqnarray*}
which implies that
\begin{eqnarray}
\langle \, 
(\sqrt{-\Delta} + V - \lambda - i \mu) 
  R (\lambda + i \mu) f , \, \psi \rangle
 =  \langle f, \, \psi \rangle
\label{eqn:7p5}
\end{eqnarray}
for any  test
function 
$\psi \in {\cal S}({\mathbb R}^3)$.
By definition of the action of $\sqrt{-\Delta}$
on $L^{2, \, -s}({\mathbb R}^3)$,
 the left hand side
of (\ref{eqn:7p5}) becomes
\begin{eqnarray}
(  R(\lambda + i \mu) f , \, 
 {\sqrt{-\Delta} \, \overline{\psi} }
\, )_{-s, s} &+&                    
(   V R(\lambda + i \mu) f , \, 
  \, \overline{\psi} 
\, )_{-s, s}                \nonumber \\              
&-& ( R (\lambda + i \mu)f, 
  \,  (\lambda - i \mu) 
\overline{\psi} \, )_{-s, s}.  \quad
\label{eqn:7p6}
\end{eqnarray}
(Note again that 
$\sqrt{-\Delta}\,\overline{\psi} 
\in L^{2, \, t}({\mathbb R}^3)$ for any 
$t < 5/2$.)
It follows from Theorem \ref{thm:BA-N2} that
\begin{eqnarray}
\lim_{\mu \downarrow 0} \, 
(  R (\lambda + i \mu) f , \, 
{\sqrt{-\Delta} \, \overline{\psi} } \,
)_{-s, s}
=
(  R^+ (\lambda) f , \, 
 {\sqrt{-\Delta} \, \overline{\psi} } \,
)_{-s, s}.
\label{eqn:7p7}
\end{eqnarray}
Similarly, we have
\begin{eqnarray}
\lim_{\mu \downarrow 0} \,  \{ 
(  V R (\lambda + i \mu) f , \, 
{\overline{\psi} } \,
)_{-s, s}
- ( R (\lambda + i \mu)f, 
  \,  (\lambda - i \mu) 
\overline{\psi} )_{-s, s}   \}   \\
=
( V R^+ (\lambda) f , \, 
 {\overline{\psi} } \,
)_{-s, s} 
- (  R^+ (\lambda) f , \, 
  \lambda {\overline{\psi} } \,
)_{-s, s}.
\label{eqn:7p8}
\end{eqnarray}
Combining (\ref{eqn:7p5}) with
 (\ref{eqn:7p6}) -- (\ref{eqn:7p8}) yields
\begin{eqnarray*}
\langle \, 
(\sqrt{-\Delta} + V - \lambda ) 
  R^+(\lambda) f , \, \psi \rangle
 =  \langle f, \, \psi \rangle
\end{eqnarray*}
for any test function 
$\psi \in {\cal S}({\mathbb R}^3)$. Thus we have
shown that 
$u^+(\lambda, \, f) = R^+(\lambda) f$
satisfies the equation (\ref{eqn:6-RC+v}).
\hfill$\square$


\newpage
\section{Generalized eigenfunctions}

Two tasks are set in this section.
One of them is to construct generalized 
eigenfunctions
of $\sqrt{-\Delta} + V(x)$ on ${\mathbb R}^3$,
 which are 
the superposition of plane waves and solutions of 
the equation
(\ref{eqn:6-RC+}), for some $\lambda$ and $f$, 
satisfying 
the outgoing or the incoming radiation condition. 
To this end,
we shall adopt the idea in 
Agmon\cite{Agmon}
(cf. Kato and Kuroda\cite{KatoKuroda}).
The other task is to 
show that the generalized 
eigenfunctions to be constructed
 are characterized as the unique
solutions to  
integral equations, which
we shall call the
modified Lippmann-Schwinger equations.

We shall write the plane wave  
$e^{ix\cdot k}$ as 
$\varphi_0 (x, \, k)$:
\begin{equation}
\varphi_0 (x, \, k) := e^{ix\cdot k}.
\label{eqn:8-0a}
\end{equation}
It should be noted that 
one can easily sees that
\begin{eqnarray*}
-\Delta_x \varphi_0 (x, \, k) = |k|^2 \varphi_0 (x, \, k), 
\end{eqnarray*}
which is a starting point when one discusses the 
generalized eigenfunction expansion for
the Schr\"odinger operator
$-\Delta + V(x)$.
On the contrary, it is not trivial to justify
\begin{equation}
\sqrt{-\Delta_x}
 \varphi_0 (x, \, k) = |k| \varphi_0 (x, \, k)
\quad \hbox{ in } \; 
{\cal S}^{\prime}({\mathbb R}^3_x),
\label{eqn:8-0b}
\end{equation}
which is formally obvious though.
The reason why (\ref{eqn:8-0b}) is nontrivial is
that $\varphi_0 (x, \, k)$ does not belong to
the Sobolev space $H^{\ell}({\mathbb R}^3_x)$
for any $\ell \in  \mathbb R$.
In fact, the Fourier transform of
$\varphi_0 (x, \, k)$ with 
respect to the variable $x$ is
a Delta-function 
$(2\pi)^{3/2} \delta (\xi - k)$,
which is obviously not a function in
$L^1_{\rm loc}({\mathbb R}^3_{\xi})$, 
whereas we have
\begin{eqnarray*}
H^{\ell}({\mathbb R}^3) =
 \{ \, f \; | \; 
 \langle \xi \rangle^{\ell} 
\hat f  \in L^2({\mathbb R}^3_{\xi})  \; \}
\end{eqnarray*}
by definition.

By virtue of some results in Umeda\cite{Umeda2}
 we shall be able to make sense of 
$\sqrt{-\Delta_x}\varphi_0 (x, \, k)$ and
prove that (\ref{eqn:8-0b}) is valid.

\vspace{15pt}
\begin{lem}\label{lem:8-1}
For every $k \in {\mathbb R}^3$, 
$\varphi_0(x, \, k)$ satisfies
the pseudodifferential equation 
{\rm (\ref{eqn:8-0b})}.
\end{lem}

{\it Proof.}
It is straightforward to see that
$\varphi_0 (x, \, k)$ belongs to
$L^{2, \, s}({\mathbb R}^3_x)$ for 
every $s < -3/2$.
This fact, together with Theorem \ref{thm:U3-5},
implies that  $\sqrt{-\Delta_x}\varphi_0 (x, \, k)$
makes sense. Taking a test function 
$\psi \in {\cal S}({\mathbb R}^3)$, we get
\begin{eqnarray}
\langle 
\sqrt{-\Delta_x}\varphi_0 (\cdot, \, k), \, \psi 
\rangle 
= ( \varphi_0 (\cdot, \, k), \, 
\sqrt{-\Delta_x} \, \overline{\psi} )_{s, -s}
\label{eqn:8-1a}
\end{eqnarray}
for all $s$ with $-5/2 < s < -3/2$, 
where we have used the fact 
that 
$\sqrt{-\Delta_x} \, \overline{\psi} 
\in L^{2, \, t}({\mathbb R}^3)$
for any $t < 5/2$.
The right hand side of (\ref{eqn:8-1a}) equals
\begin{eqnarray*}
\int e^{ix\cdot k} \,
\overline{\sqrt{-\Delta} \, \overline{\psi (x)} }
\, dx
&=&
(2\pi)^{3/2} \,
\overline{
{\cal F}[\sqrt{-\Delta} \, \overline{\psi}\,](k) }   \\
&=& (2\pi)^{3/2} \, |k| \,
\overline{
{\cal F}[\, \overline{\psi}\,]  (k) }.
\end{eqnarray*}
Noting that 
\begin{eqnarray*}
\overline{
{\cal F}[\, \overline{\psi}\,]  (k) } 
=(2\pi)^{-3/2} 
\int \varphi_0(x, \, k) \, \psi(x) \, dx ,
\end{eqnarray*}
we obtain
\begin{eqnarray}
 ( \varphi_0 (\cdot, \, k), \, 
\sqrt{-\Delta_x} \, \overline{\psi} )_{s, -s}
&=&
\int |k| \, \varphi_0(x, \, k) \, 
\psi(x) \, dx                       \nonumber\\
\noalign{\vskip 4pt}
&=&
\langle |k| \, 
\varphi_0( \cdot, \, k), \, \psi 
\rangle.
\label{eqn:8-1b}
\end{eqnarray}
Combining (\ref{eqn:8-1b}) with (\ref{eqn:8-1a})
gives the lemma.
\hfill$\square$

\vspace{15pt}

Following 
Agmon\cite{Agmon}, we define two families
of generalized eigenfunctions of
$\sqrt{-\Delta} + V(x)$ on ${\mathbb R}^3$ 
by
\begin{equation}
\varphi^{\pm}(x, \, k):= 
   \varphi_0(x, \, k) - R^{\mp}(|k|)
 \{ V(\cdot) \varphi_0(\cdot, \, k) \}(x)
\label{eqn:8-2}
\end{equation}
for $k$ with 
$|k| \in (0, \, +\infty) \setminus \sigma_p(H)$. 
Note that the second terms on the
right hand side of (\ref{eqn:8-2}) 
make sense, provided that
$|V(x)| \le C \langle x \rangle^{-\sigma}$,
$\sigma >2$. In fact,
$V(\cdot) \varphi(\cdot, \, k) \in L^{2, s}({\mathbb R}^3)$
for all $s$ with $1/2 < s < \sigma - 3/2$.

\vspace{15pt}

\begin{thm}\label{thm8-1}
Let $\sigma >2$. 
If $|k| \in (0, \, +\infty) \setminus \sigma_p(H)$,
then both $\varphi^{\pm}(x, \, k)$ satisfy
the equation
\begin{eqnarray}
(\sqrt{-\Delta_x} + V(x) ) u = |k| u  \quad
\hbox{ in } \;{\cal S}^{\prime}({\mathbb R}_x^3).
\label{eqn:8-3a}
\end{eqnarray}
\end{thm}

{\it Proof.}
As remarked just before the theorem, we
see that
$V(\cdot) \varphi_0(\cdot, \, k)$ 
belongs to $L^{2, s}({\mathbb R}^3)$ for all
$s$ with
$1/2 < s < \sigma - 3/2$. Hence, by
Theorem \ref{thm:RC}(ii), we get
\begin{eqnarray}
(\sqrt{-\Delta_x} + V(x) -|k|) 
\big[ \, R^{\mp}(|k|) \{ 
  V(\cdot) \varphi_0(\cdot, \, k) \} \, \big](x) 
                       \qquad \qquad       \nonumber\\
 = V(\cdot) \varphi_0(\cdot, \, k)       \quad
\hbox{ in } \;{\cal S}^{\prime}({\mathbb R}_x^3),
\label{eqn:8-3b}
\end{eqnarray}
which, together with Lemma \ref{lem:8-1},
implies that
\begin{align*}
(\sqrt{-\Delta_x} &+ V(x) ) 
\varphi^{\pm}(x, \, k)              \\
\noalign{\vskip 4pt}
 &=
(\sqrt{-\Delta_x} + V(x) ) \varphi_0(x, \, k) \\
& \qquad - (\sqrt{-\Delta_x} + V(x) )
\big[ R^{\mp}(|k|) \{ 
  V(\cdot) \varphi_0(\cdot, \, k)  \} \big](x)   \\
\noalign{\vskip 4pt}
&= 
|k|\varphi_0(x, \, k) + V(x)\varphi_0(x, \, k) \\
& \qquad -V(x)  \varphi_0(x, \, k)  
- |k| \, \big[ R^{\mp}(|k|) \{ 
  V(\cdot) \varphi_0(\cdot, \, k) \} \big](x) \\
\noalign{\vskip 4pt}
&=
|k| \big[ \,
\varphi_0(x, \, k) -
R^{\mp}(|k|) \{ 
V(\cdot) \varphi_0(\cdot, \,k)  \}(x)  \, \big] .  
\end{align*}
By the definition (\ref{eqn:8-2}), this
gives the theorem.
\hfill$\square$

\vspace{15pt}
{\it Remark.} For each $k$ with
 $|k| \in (0, \, +\infty) \setminus \sigma_p(H)$, the 
generalized eigenfunctions $\varphi^{\pm}(x, \, k)$
are unique in the following sense:
If $\tilde{\varphi}^+(x, \, k)$ 
(resp. $\tilde{\varphi}^-(x, \, k)$) satisfies 
the equation (\ref{eqn:8-3a}), and
in addition, 
$\tilde{\varphi}^+(x, \, k) - \varphi_0(x, \, k)$
(resp. 
$\tilde{\varphi}^-(x, \, k) - \varphi_0(x, \, k)$)
belongs to
$L^{2, -s}({\mathbb R}^3)
  \cap H^1_{\rm loc}({\mathbb R}^3)$,
$1/2 < s < \min ( \, \sigma/2, 1)$, 
and satisfies the incoming radiation condition 
(\ref{eqn:5-0b}) (resp. the outgoing radiation 
condition (\ref{eqn:5-0a})), then
$\tilde{\varphi}^+(x, \, k) =\varphi^+(x, \, k)$
(resp.\linebreak
$\tilde{\varphi}^-(x, \, k) =\varphi^-(x, \, k)$).
This is a direct consequence of assertion (i) of
Theorem \ref{thm:RC}.

\vspace{15pt}

We are in a position to
 show that the generalized 
eigenfunctions $\varphi^+ (x,\, k)$ and
$\varphi^- (x,\, k)$, defined by (\ref{eqn:8-2}), 
are characterized as the unique
solutions to the 
integral equations
\begin{equation}
\varphi (x) = \varphi_0(x, \, k) - 
  \int_{{\mathbb R}^3} g_{|k|}^-(x-y) \, 
V(y) \, \varphi (y) \,dy
\label{eqn:9-1}
\end{equation}
and
\begin{equation}
\varphi (x) = \varphi_0(x, \, k) - 
  \int_{{\mathbb R}^3} g_{|k|}^+(x-y) \, 
V(y) \, \varphi (y) \,dy
\label{eqn:9-2}
\end{equation}
respectively. (Recall that $g_{\lambda}^{\pm}(x -y)$ are
the integral kernels of the 
boundary values $R_0^{\pm}(\lambda)$.
See Theorem \ref{thm:R0ker}.)
We call (\ref{eqn:9-1}) and (\ref{eqn:9-2}) the
modified Lippmann-Schwinger equations,
because the leading terms of $g_{\lambda}^{\pm}(x-y)$ are
the same, up to a constant, as the  integral
kernels of the Lippmann-Schwinger equations, namely,
\begin{equation*}
g_{\lambda}^{\pm}(x-y) 
= \frac{\lambda}{\, 2\pi \,}
\cdot\frac{e^{\pm i \lambda |x-y|}}{\, |x-y| \,}
 + O( |x - y|^{-2} )   
\quad \hbox{as } |x - y| \to  + \infty.
\end{equation*}
(Recall (\ref{eqn:3r}) and (\ref{eqn:4e+}).)

Our generalized eigenfunctions 
$\varphi^{\pm}(x, \, k)$ are expected to behave
like the plane wave $\varphi_0(x, \, k)$,
which belongs to $L^{2, -s}({\mathbb R}^3)$ 
only for 
$s > 3/2$. 
Thus it is natural to take  
$L^{2, -s}({\mathbb R}^3)$, 
with $s > 3/2$, to be
the space of functions in which
we deal with the integral equations
(\ref{eqn:9-1}) and (\ref{eqn:9-2}).
It is evident from Theorem \ref{thm:R0ker}
that 
(\ref{eqn:9-1}) and (\ref{eqn:9-2}) can be
formally rewritten in the forms
$\big(I + R_0^-(|k|)V\big)\varphi = 
\varphi_0(\cdot, \, k)$ and
$\big(I + R_0^+(|k|)V\big)\varphi = 
\varphi_0(\cdot, \, k)$ 
respectively.
For these reasons, we prepare the following 
lemma, which is a variant of Lemma \ref{lem:5-6}.
The only difference between Lemmas \ref{lem:5-6} 
and \ref{lem:9-1} lies in their assumptions.
In Lemma \ref{lem:9-1}, $s$ is allowed to be 
greater than $3/2$.

\vspace{15pt}
\begin{lem}\label{lem:9-1}
Let $\sigma >2$, and
suppose that $1/2 < s < \sigma - 1/2$. 
Then the conclusions of Lemma \ref{lem:5-6} 
hold.
\end{lem}

{\it Proof.} We  only give the proof 
of (\ref{eqn:5-6a})
in the case where the superscripts are 
the plus sign. 
The proof 
of (\ref{eqn:5-6a})
in the other case and the proof of
(\ref{eqn:5-6b}) are similar.

It is obvious that we shall 
follow the line of the proof
of Lemma  \ref{lem:5-6}.
By assumption, we can 
choose $t$ so that
\begin{equation}
1/2 < t 
     < \min\,(s, \, \sigma - s).
\label{eqn:9-4}
\end{equation}
We note that $R_0^+(z)$ and $R^+(z)$ can
be regarded as 
${\bf B}(L^{2, \, t}, L^{2, -t})$-valued
continuous functions on
${\mathbb C}^+ \cup (0, \, +\infty)$
and 
${\mathbb C}^+\cup 
  \{ (0, \, +\infty) \setminus \sigma_p(H) \}$
respectively, as
mentioned in the proof of Lemma \ref{lem:5-4}.
From this fact, we can deduce that
$R_0^+(z)$ and $R^+(z)$ are
${\bf B}(L^{2, \, t}, L^{2, -s})$-valued
continuous functions on
${\mathbb C}^+ \cup (0, \, +\infty)$
and 
${\mathbb C}^+\cup 
  \{ (0, \, +\infty) \setminus \sigma_p(H) \}$
respectively, since $-s < -t$.
In view of (\ref{eqn:9-4}) we have
$V \in {\bf B}(L^{2, -s}, L^{2,\,t})$.
Therefore, $R_0^+(z)V$ and $R^+(z)V$ are
${\bf B}(L^{2, -s})$-valued
continuous functions on
${\mathbb C}^+ \cup (0, \, +\infty)$
and 
${\mathbb C}^+\cup 
  \{ (0, \, +\infty) \setminus \sigma_p(H) \}$
respectively. 
Recalling (\ref{eqn:5-6c}), which  
was shown to be valid for all 
$z \in {\mathbb C}^+$, we conclude that
(\ref{eqn:5-6a})
in the case where the superscripts are \lq\lq +" 
holds. 
\hfill$\square$

\vspace{15pt}

\begin{thm}\label{thm:8-2}
Let $\sigma >2$, and suppose that 
$3/2 < s < \sigma - 1/2$.
If $|k| \in (0, \, +\infty) \setminus \sigma_p(H)$,
then  $\varphi^+(x, \, k)$  and
 $\varphi^-(x, \, k)$ are the unique solution of 
the
modified Lippmann-Schwinger equations
{\rm (\ref{eqn:9-1})} and 
{\rm (\ref{eqn:9-2})} in $L^{2, -s}({\mathbb R}^3_x)$ 
respectively.
\end{thm}

{\it Proof.} We shall give the proof only 
for $\varphi^+(x, \, k)$.

It follows from the definition (\ref{eqn:8-2}) 
that
\begin{equation}
\varphi^+(\cdot, \, k)= 
   \big(I- R^-(|k|)V\big)\varphi_0(\cdot, \, k),
\label{eqn:9-5} 
\end{equation}
where we regard $\varphi_0(\cdot, \, k)$
as a function belonging to 
$L^{2, -s}({\mathbb R}^3_x)$.
Combining (\ref{eqn:9-5}) with (\ref{eqn:5-6b}),
 we have
\begin{equation}
\big(I+ R^-_0(|k|)V\big)\varphi^+(\cdot, \, k)= 
   \varphi_0(\cdot, \, k),
\label{eqn:9-6} 
\end{equation}
from which we obtain
\begin{equation}
\varphi^+(\cdot, \, k)= 
   \varphi_0(\cdot, \, k) -
R^-_0(|k|)V\varphi^+(\cdot, \, k).
\label{eqn:9-6a} 
\end{equation}
Since the integral kernel of $R^-_0(|k|)$ is
given by
$g_{|k|}^-(x - y)$, we conclude from (\ref{eqn:9-6a})
that  $\varphi^+(x, \, k)$ satisfies the
modified Lippmann-Schwinger equation (\ref{eqn:9-1}).
Uniqueness follows from (\ref{eqn:9-6}) and 
(\ref{eqn:5-6a}).
\hfill$\square$

\newpage
\section{Continuity of the generalized eigenfunctions}

The aim of this section is to prove the following
result.

\vspace{15pt}

\begin{thm} \label{thm:gefcont}     
Let $\sigma >2$. Then the generalized eigenfunctions
$\varphi^{\pm}(x, \, k)$ defined by 
{\rm (\ref{eqn:8-2})} have the following
properties:
\begin{description}
\item[\rm (i)] 
For each interval
$[a, \, b] \subset (0, \, +\infty) \setminus \sigma_p (H)$, 
there exists
a constant $C_{ab}$, depending on $a$ and $b$, such that
\begin{equation}
|\varphi^{\pm}(x, \, k) | \le  C_{ab}
\label{eqn:gef-1}
\end{equation}
for all $(x, \, k) \in 
{\mathbb R}^3 \times 
\{ \, k \, | \, a \le |k| \le b \, \}$. 
\item[\rm (ii)] 
$\varphi^{\pm}(x, \, k)$ are continuous
functions      
on ${\mathbb R}^3_x \times \big\{ \; k \; | \; |k| 
\in (0, \, +\infty) \setminus
   \sigma_p(H) \; \big\}$.
%
%
\end{description}
\end{thm}

\vspace{15pt}

We shall give a proof of Theorem \ref{thm:gefcont} 
by means of a series of lemmas.
Hence, throughout the present section
we shall assume that
$$
\sigma >2
$$
without saying so every time.
We shall first prepare a few lemmas 
and then prove assertion (i) of 
Theorem \ref{thm:gefcont}. 
We shall next show a few lemmas, 
of  which combination directly gives
a proof of assertion (ii) of Theorem \ref{thm:gefcont}. 
The estimate (\ref{eqn:gef-1}) will be
useful in the discussions for the proof 
of assertion (ii).

\vspace{15pt}

\begin{lem}\label{lem:9-a}
If $s > 3/2$, then $\varphi^{\pm}(\cdot, \, k)$
are $L^{2, -s}({\mathbb R}^3_x)$-valued continuous
functions on
$\, \big\{ \; k \; | \; |k| 
\in (0, \, +\infty) \setminus
   \sigma_p(H) \; \big\}$.
\end{lem}

{\it Proof.}
We  note that $\varphi_0(\cdot, \, k)$ is
$L^{2, -s}({\mathbb R}^3_x)$-valued continuous
function on
$\, {\mathbb R}^3_k$.
On the other hand,
for any $t$ with
$1/2 < t < \sigma -3/2$,
$V(\cdot)\varphi_0(\cdot, \, k)$ is
$L^{2, t}({\mathbb R}^3_x)$-valued continuous
function on
$\, {\mathbb R}^3_k$ (see 
the assumption (\ref{eqn:Vdecay})).
This fact, together with 
 Theorem \ref{thm:BA-N2} (iii), implies that
$R^{\mp}(|k|)\{ 
V(\cdot)\varphi_0(\cdot, \, k) \}$ 
are $L^{2, -t}({\mathbb R}^3_x)$-valued continuous
functions on
$\, \big\{ \; k \; | \; |k| 
\in (0, \, +\infty) \setminus
   \sigma_p(H) \; \big\}$.
Since $t$ can be taken to be less than $s$,
it follows that 
$R^{\mp}(|k|)\{ 
V(\cdot)\varphi_0(\cdot, \, k) \}$ 
are $L^{2, -s}({\mathbb R}^3_x)$-valued continuous
functions on
$\, \big\{ \; k \; | \; |k| 
\in (0, \, +\infty) \setminus
   \sigma_p(H) \; \big\}$.
In view of 
 the definition (\ref{eqn:8-2}),
 we have proved the lemma. 
 \hfill$\square$

\vspace{15pt}

\begin{lem}\label{lem:9-b}
If $s > 3/2$, then 
$V(\cdot)\varphi^{\pm}(\cdot, \, k)$
are $L^{2, \, \sigma -s}({\mathbb R}^3_x)$-valued 
continuous functions on
$\, \big\{ \; k \; | \; |k| 
\in (0, \, +\infty) \setminus
   \sigma_p(H) \; \big\}$.
\end{lem}

{\it Proof.} The lemma is a direct consequence
of Lemma \ref{lem:9-a} and 
the assumption (\ref{eqn:Vdecay}).
\hfill$\square$

\vspace{15pt}

In the rest of this section, we assume that
$s$ satisfies the inequalities
\begin{equation}
\frac{3}{\, 2 \,} < s <
\sigma - \frac{1}{\, 2 \,}
\label{eqn:9-1-before}
\end{equation}
In order to prove assertion (i) of 
Theorem \ref{thm:gefcont}, we need  intermediate
estimates, which only assure that, for each $k$,
$\varphi^{\pm}(x, \, k)$ are sums of
bounded functions of $x$ and 
functions of $x$ belonging to
$L^6({\mathbb R}^3)\cap L^{2, \, -t}({\mathbb R}^3)$
for all $t > 3/2$.
To derive the intermediate estimates mentioned above,
 we appeal to
Theorem \ref{thm:8-2}; 
assuming that
$|k| \in (0, \, +\infty) \setminus \sigma_p(H)$,
 we have
\begin{equation}
\varphi^{\pm}(x, \, k) =
\varphi_0(x, \, k) - 
G_{|k|}^{\mp} \big( \, V(\cdot) 
\varphi^{\pm}(\cdot, \, k) \big) (x) ,
\label{eqn:9-1-change}
\end{equation}
(see (\ref{eqn:4gg}) and
(\ref{eqn:4-last}) for the notation 
$G_{|k|}^{\mp}$).
According to the identities (\ref{eqn:5ahlalah}),
we then decompose $\varphi^{\pm}(x, \, k)$ 
into two parts:
\begin{equation}
\varphi^{\pm}(x, \, k) =
\psi_0^{\pm}(x, \, k) + \psi_1^{\pm}(x, \, k),
\label{eqn:9-2-change}
\end{equation}
where
\begin{eqnarray}
\psi_0^{\pm}(x, \, k) &:=& \varphi_0(x, \, k) -
K_{|k|}^{\mp} \big( \, V(\cdot)   
\varphi^{\pm}(\cdot, \, k) \big) (x)  \qquad  \nonumber\\
&& \qquad \qquad \qquad \qquad 
-M_{|k|} \big( \, V(\cdot) 
\varphi^{\pm}(\cdot, \, k) \big) (x) ,   
\label{eqn:9-2-0} \\
\noalign{\vskip 4pt}
\psi_1^{\pm}(x, \, k) &:=& - 
G_0 \big( \, V(\cdot) 
\varphi^{\pm}(\cdot, \, k) \big) (x) .
\label{eqn:9-2-1}
\end{eqnarray}

\vspace{15pt}

\begin{lem}\label{lem:9-c}
Suppose that 
$[a, \, b] \subset (0, \, +\infty) 
\setminus \sigma_p(H)$.
Then there exists a constant $C_{ab}$, depending
on $a$ and $b$, such that
\begin{equation}
|\psi_0^{\pm}(x, \, k) | \le  C_{ab}
\label{eqn:psi-0}
\end{equation}
for all 
$(x, \, k) \in {\mathbb R}^3 \times
\{ \, k \, | \, a \le |k| \le b \, \}$.
\end{lem}

{\it Proof.}
Let $k$ satisfy $a \le |k| \le b$. 
Appealing to the definition (\ref{eqn:9-2-0}), 
we have
\begin{eqnarray}
|\psi_0^{\pm}(x, \, k)| &\le&  1 +  \,
|\, K_{|k|}^{\mp} \big( \, V(\cdot)   
\varphi^{\pm}(\cdot, \, k) \big) (x)|   \nonumber \\
&& \qquad \qquad \qquad 
+ \, |\, M_{|k|} \big( \, V(\cdot)   
\varphi^{\pm}(\cdot, \, k) \big) (x)|   \nonumber \\
&\le& 1 +  \, C\, b \, 
\Vert  V(\cdot)   
   \varphi^{\pm}(\cdot, \, k)
\Vert_{L^{2, \, \sigma -s}}    \nonumber \\
&& \qquad \qquad \qquad \quad
+ \, 
C_{ab} \,
\Vert
 V(\cdot)   
\varphi^{\pm}(\cdot, \, k)
\Vert_{L^2},
\label{eqn:psi-0-1}
\end{eqnarray}
where we have used 
Lemmas \ref{lem:ker4} and 
 \ref{lem:kerm+1}
(note that $\sigma - s > 1/2$ by
(\ref{eqn:9-1-before})). 
Here we note that
the constants $C$  and 
$C_{ab}$ in (\ref{eqn:psi-0-1})
are independent of $k$ with
$a \le |k| \le b$.
Lemma \ref{lem:9-b}, together with 
(\ref{eqn:psi-0-1}), implies the lemma.
\hfill$\square$

\vspace{15pt}
\begin{lem}\label{lem:9-d}
Let 
$|k|\in (0, \, +\infty) \setminus \sigma_p(H)$. 
Then we have
\begin{equation*}
\psi_1^{\pm}(\cdot, \, k) \in 
L^6({\mathbb R}^3)
\cap 
L^{2, \, -t}({\mathbb R}^3)
\end{equation*}
for every $t > 3/2$.
Moreover, for each compact 
interval
$[a, \, b] \subset (0, \, +\infty) 
\setminus \sigma_p(H)$ 
and each $t >3/2$,
there corresponds a positive
constant $C_{tab}$ such that
\begin{equation*}
\Vert \psi_1^{\pm}(\cdot, \, k) 
  \Vert_{L^6} 
+
\Vert \psi_1^{\pm}(\cdot, \, k) 
 \Vert_{L^{2, -t}} 
\le C_{tab}
\end{equation*}
for all $k$ with
$a \le |k| \le b$.
\end{lem}

{\it Proof.}
Since $\sigma -s >1/2$ by 
(\ref{eqn:9-1-before}), 
it follows from Lemma \ref{lem:9-b} that
$V(\cdot)\varphi^{\pm}(\cdot, \, k) 
\in L^2({\mathbb R}^3)$. 
Then the definition of $\psi_1^{\pm}$ and
the inequality (\ref{eqn:4b-}) show that
\begin{equation}
\Vert
\psi_1^{\pm}(\cdot, \, k) 
 \Vert_{L^6}
\le
C \, 
\Vert
V(\cdot) \varphi^{\pm}(\cdot, \, k) 
 \Vert_{L^2} ,
\label{eqn:psi-1-a}
\end{equation}
where $C$ is a constant
independent of $k$.
Similarly, 
the definition of $\psi_1^{\pm}$ and
Lemma \ref{lem:ker1}(ii)  give
\begin{equation}
\Vert
\psi_1^{\pm}(\cdot, \, k) 
 \Vert_{L^{2, -t}}
\le
C_t \, 
\Vert
V(\cdot) \varphi^{\pm}(\cdot, \, k) 
 \Vert_{L^2} 
\label{eqn:psi-1-b}
\end{equation}
for every $t>3/2$, where
the constant $C_t$ is dependent
on $t$ but independent of $k$. 
The assertions of the lemma
now follow from (\ref{eqn:psi-1-a}),
(\ref{eqn:psi-1-b}) and Lemma \ref{lem:9-b}.
\hfill$\square$

\vspace{15pt}

{\bf Proof of assertion(i) of 
Theorem \ref{thm:gefcont} }
In view of (\ref{eqn:9-2-change}) and
Lemma \ref{lem:9-c}, it is
sufficient to show that
there exists a constant
$C_{ab}$ such that
\begin{equation}
|\psi_1^{\pm}(x, \, k) | \le C_{ab}
\label{eqn:psi-1-c}
\end{equation}
for all
$(x, \, k) \in {\mathbb R}^3 \times
\{ \, k \, | \, a \le |k| \le b \, \}$.

It follows from
(\ref{eqn:9-2-change})
and 
(\ref{eqn:9-2-1}) that
\begin{equation}
\psi_1^{\pm}(x, \, k) 
= 
-G_0 \big( V(\cdot) 
\psi_0^{\pm}(\cdot, \, k) \big) (x) 
-G_0 \big( V(\cdot) 
\psi_1^{\pm}(\cdot, \, k) \big) (x) .
\label{eqn:psi-1-d}
\end{equation}
We apply Lemma \ref{lem:9-c} to the
first term on the
right hand side  of (\ref{eqn:psi-1-d})
and appeal to the definition (\ref{eqn:2e}) 
of $G_0$, and obtain
\begin{eqnarray}
\big| 
G_0 \big( V(\cdot) 
\psi_0^{\pm}(\cdot, \, k) \big) (x)  
\big|
 \le
\frac{\, \Vert  
V(\cdot)  \langle \cdot \rangle^{\sigma} 
\Vert_{L^{\infty}} }{ 2 \pi^2 }
 \int_{{\mathbb R}^3}  
\frac{C_{ab}}
{\,|x-y|^2 \langle y \rangle^{\sigma} \,} \, dy ,
\label{eqn:psi-1-e}
\end{eqnarray}  
where the constant $C_{ab}$ is the same 
as in (\ref{eqn:psi-0}), and is 
independent of $k$ with
$a \le |k| \le b$. By virtue
of Lemma A.1 in Appendix,
the function of $x$ defined by 
the integral on the right hand side
of (\ref{eqn:psi-1-e}) is bounded 
on ${\mathbb R}^3$.
Thus the
first term on the
right hand side  of (\ref{eqn:psi-1-d}) possesses
the desired estimate.
To handle the second term on the
right hand side  of (\ref{eqn:psi-1-d}), 
we decompose it into 
two parts:
\begin{eqnarray}
G_0 \big( V(\cdot) 
\psi_1^{\pm}(\cdot, \, k) \big) (x)  
&=&
\frac{1}{\, 2 \pi^2 \,}
\int_{|x -y| \le 1}  
\frac{\, V(y) \psi_1^{\pm}(y, \, k) \,}
{|x - y|^2 }   
\, dy      \nonumber\\
\noalign{\vskip 4pt}
&& \;\; + \;
\frac{1}{\, 2 \pi^2 \,}
\int_{|x -y| \ge 1}
\frac{\, V(y) \psi_1^{\pm}(y, \, k) \,}
{|x - y|^2 }   
\, dy      \nonumber\\
\noalign{\vskip 4pt}
&=:& I^{\pm}(x, \, k) +
I\!I^{\pm}(x, \, k).
\label{eqn:psi-1-f}
\end{eqnarray}
We apply the H\"older inequality to
$I^{\pm}(x, \, k)$, and get
\begin{eqnarray}
|I^{\pm}(x, \, k)|  
&\le&
\frac{1}{\, 2 \pi^2 \,}
\Big\{
\int_{|x - y| \le 1}
\big( \frac{1}{\, |x - y| ^2\,}
\big)^{6/5} \, dy
\Big\}^{5/6}     \nonumber \\
\noalign{\vskip 4pt}
&& \;\; \quad \times  \;
\Big\{
\int_{|x - y| \le 1}
\big|  V(y) \psi_1^{\pm}(y, \, k)
\big|^6 \, dy
\Big\}^{1/6}        \nonumber\\
\noalign{\vskip 4pt}
&\le&
\frac{1}{\, 2 \pi^2 \,} 
\Big\{
\int_{|y| \le 1}
 |y|^{-12 / 5}\, dy
\Big\}^{5/6}
\Vert V \Vert_{L^{\infty}}
\Vert 
\psi_1^{\pm}(\cdot, \, k)   
\Vert_{L^6} .  \quad \;\;
\label{eqn:psi-1-g}
\end{eqnarray}
Since $-12/5 > -3$, Lemma \ref{lem:9-d}
and (\ref{eqn:psi-1-g}) imply that
$I^{\pm}(x, \, k)$ satisfy the 
desired estimate. Similarly, we apply 
the Schwarz inequality to
$I\!I^{\pm}(x, \, k)$, and we obtain
\begin{eqnarray}
|I\!I^{\pm}(x, \, k)|  
&\le&
\frac{1}{\, 2 \pi^2 \,}
\Big\{
\int_{|x - y| \ge 1}
\big( \frac{1}{\, |x - y| ^2\,}
\big)^2 \, dy
\Big\}^{1/2}     \nonumber \\
\noalign{\vskip 4pt}
&& \;\; \quad \times  \;
\Big\{
\int_{|x - y| \ge 1}
\big|  V(y) \psi_1^{\pm}(y, \, k)
\big|^2 \, dy
\Big\}^{1/2}        \nonumber\\
\noalign{\vskip 4pt}
&\le&
\frac{1}{\, 2 \pi^2 \,} 
\Big\{
\int_{|y| \ge 1}
 |y|^{-4}\, dy
\Big\}^{2}         \nonumber \\
\noalign{\vskip 4pt}
&& \quad \quad
\times  \;\Vert V(\cdot) 
\langle \cdot \rangle^{\sigma} 
\Vert_{L^{\infty}}
\Vert 
\psi_1^{\pm}(\cdot, \, k)   
\Vert_{L^{2, -\sigma}} .  
\label{eqn:psi-1-h}
\end{eqnarray}
Since $-4 < -3$ and
$\sigma > 3/2$, Lemma \ref{lem:9-d},
together with (\ref{eqn:psi-1-h}), 
implies that
$I\!I^{\pm}(x, \, k)$ have the 
desired estimate.
Summing up, we have shown that 
(\ref{eqn:psi-1-c})
holds for all
$(x, \, k)$ in ${\mathbb R}^3 \times
\{ \, k \, | \, a \le |k| \le b \, \}$.
\hfill$\square$

\vspace{15pt}

In order to prepare lemmas, of which combination
will directly
give the proof of  assertion (ii)
of Theorem \ref{thm:gefcont}, it is convenient to
write
\begin{eqnarray}
\psi_{0\kappa}^{\pm}(x, \, k) &:=& 
 -K_{|k|}^{\mp} \big( \, V(\cdot)   
\varphi^{\pm}(\cdot, \, k) \big) (x) ,   
\label{eqn:psi-0K} \\
\psi_{0\mu}^{\pm}(x, \, k) &:=& 
- M_{|k|} \big( \, V(\cdot) 
\varphi^{\pm}(\cdot, \, k) \big) (x) .
\label{eqn:psi-0M}
\end{eqnarray}
According to 
(\ref{eqn:9-2-change})--(\ref{eqn:9-2-1}),
we then have
\begin{equation}
\varphi^{\pm}(x, \, k)
= \varphi_0(x, \, k) +
\psi_{0\kappa}^{\pm}(x, \, k) + 
\psi_{0\mu}^{\pm}(x, \, k) + 
\psi_1^{\pm}(x, \, k).
\label{eqn:9-3-1}
\end{equation}

\vspace{15pt}
\begin{lem}\label{lem:9-e}
 $\psi_{0\kappa}^{\pm}(x, \, k)$ are
continuous  on
${\mathbb R}^3_x \times \big\{ \; k \; | \; |k| 
\in (0, \, +\infty) \setminus
   \sigma_p(H) \; \big\}$.
\end{lem}

{\it Proof.}
Let $(x_0, \, k_0)$ be an arbitrary
point in 
${\mathbb R}^3_x \times \big\{ \; k \; | \; |k| 
\in (0, \, +\infty) \setminus
   \sigma_p(H) \; \big\}$.
We shall show that
\begin{equation}
\psi_{0\kappa}^{\pm}(x, \, k) 
 \to 
\psi_{0\kappa}^{\pm}(x_0, \, k_0)
\quad 
\hbox{\rm as }{(x, \, k)\to (x_0, \, k_0)}.
\label{eqn:9-e-1}
\end{equation}

Let $\varepsilon > 0$ be given.
One can then  choose $r>0$ so that
\begin{equation}
\frac{|k_0|}{\, 2\pi \,} \,
\Vert V \Vert_{L^{\infty}} 
\Big\{\!\!
\sup_{y \in {\mathbb R}^3 
    \atop
 |k - k_0| \le r} 
\!\!\!
 | \varphi^{\pm}(y, \, k) |
\Big\}
\int_{|y| \le 2r} 
\frac{1}{\, |y| \,} \,dy < \varepsilon.
\label{eqn:9-e-2}
\end{equation}
Note that, by virtue of assertion (i) of
Theorem  \ref{thm:gefcont}, the supremum in
(\ref{eqn:9-e-2}) is finite.
To show (\ref{eqn:9-e-1}), we write
\begin{eqnarray}
&&\psi_{0\kappa}^{\pm}(x, \, k) 
- 
\psi_{0\kappa}^{\pm}(x_0, \, k_0)  \nonumber\\
\noalign{\vskip 3pt}
&=& \{ \psi_{0\kappa}^{\pm}(x, \, k) 
- 
\psi_{0\kappa}^{\pm}(x_0, \, k) \}   
+
\{ \psi_{0\kappa}^{\pm}(x_0, \, k) 
- 
\psi_{0\kappa}^{\pm}(x_0, \, k_0) \}  \qquad  
\label{eqn:9-e-3}\\
\noalign{\vskip 3pt}
&=:&
 I^{\pm}_{0\kappa}(x, \, k)
+ I\!I^{\pm}_{0\kappa}(k).    \nonumber
\end{eqnarray}
If $|x - x_0| \le r$ and $|k - k_0| \le r$,
we then have, 
appealing to the definition (\ref{eqn:4a}),
\begin{eqnarray}
|I^{\pm}_{0\kappa}(x, \, k)| 
\le
2 \frac{|k|}{\, |k_0| \,} \varepsilon  
+ 
\frac{|k|}{\, 2 \pi \,}
\Big\{\!\!
\sup_{y \in {\mathbb R}^3 
    \atop
 |k - k_0| \le r} 
\!\!\!
 | \varphi^{\pm}(y, \, k) |
\Big\} \times    \qquad \qquad \qquad \nonumber \\
\times \int_{{\mathbb R}^3} \Big|  \,  
\Big\{
1_{E(x, \, 2r)}(y) 
\frac{\, e^{\pm i|k| |x - y|}\, }{|x - y|} 
-
1_{E(x_0, \, 2r)}(y) 
\frac{\, e^{\pm i|k| |x_0 - y|}\, }{|x_0 - y|} 
\Big\} V(y)
\Big| \, dy,
\label{eqn:9-e-4}
\end{eqnarray}
where 
$E(x, \, 2r)= \{ \, y \, | \, |x - y| > 2r \, \}$ 
 and we have used (\ref{eqn:9-e-2}).
We note here that
\begin{equation}
1_{E(x, \, 2r)}(y) 
\frac{1}{\,|x - y|\,}  \le
\frac{3}{\, 2 \,} \times
1_{E(x_0, \, r)}(y) 
\frac{1}{\,|x_0 - y|\,}
\label{eqn:9-e-4-1}
\end{equation}
whenever $|x - x_0| \le r$.
Hence,
the integrand in (\ref{eqn:9-e-4}) is bounded,
for all $(x, \, k)$ with $|x - x_0| \le r$,
 by
the function
\begin{equation}
\frac{\, 5 \,}{2} \times 
1_{E(x_0, \, r)}(y)
\frac{1}{\, |x_0 - y| \,} \, |V(y)|
\label{eqn:9-e-5}
\end{equation}
which is in $L^1({\mathbb R}^3_y)$
(recall that we made
the assumption (\ref{eqn:Vdecay})
with $\sigma >2$).
Hence we can apply the Lebesgue dominated convergence
theorem to the integral in (\ref{eqn:9-e-4}), and 
deduce that
\begin{equation}
\limsup_{(x,k) \to (x_0, k_0)}
| I_{0\kappa}^{\pm}(x, \, k) | \le 2 \varepsilon .
\label{eqn:9-e-6}
\end{equation}
In a similar fashion to (\ref{eqn:9-e-4}), 
if $|k - k_0| \le r$,
we have
\begin{eqnarray}
|I\!I^{\pm}_{0\kappa}(k)| 
\le | K_{|k|}^{\mp}\big( \, V(\cdot)   
\varphi^{\pm}(\cdot, \, k) \big) (x_0)
- K_{|k_0|}^{\mp}\big( \, V(\cdot)   
\varphi^{\pm}(\cdot, \, k) \big) (x_0) | \quad \nonumber\\
\noalign{\vskip 3pt}
 + | K_{|k_0|}^{\mp}\big( \, V(\cdot)   
\varphi^{\pm}(\cdot, \, k) \big) (x_0)
- K_{|k_0|}^{\mp}\big( \, V(\cdot)   
\varphi^{\pm}(\cdot, \, k_0) \big) (x_0) |   \nonumber\\
\noalign{\vskip 3pt}
\le
\Big( 
\frac{\,|k| \,}{\, |k_0| \,} + 3
\Big)
\varepsilon  
+ 
\frac{1}{\, 2 \pi \,}
\int_{{\mathbb R}^3} \Big|  \,  
1_{E(x_0, \, 2r)}(y) 
\frac{1}{\,|x_0 - y|\,} \times
     \qquad \qquad \qquad    \nonumber \\
\noalign{\vskip 4pt}
\Big| \, 
|k|e^{\mp i|k| |x_0 - y|}
-
 |k_0|e^{\mp i|k_0| |x_0 - y|} \Big| \,
\big| V(y) \varphi^{\pm}(y, \,k) \big|
 \, dy  \qquad                         
\label{eqn:9-e-7}\\
\noalign{\vskip 4pt}
+ \; \frac{\,|k_0|\,}{\, 2\pi \,}
\int_{{\mathbb R}^3} \Big|  \,  
1_{E(x_0, \, 2r)}(y) 
\frac{e^{\mp i|k_0| |x_0 - y|}}{\,|x_0 - y|\,} 
V(y)  \Big| 
\times
     \qquad \qquad \qquad    \nonumber \\
\noalign{\vskip 4pt}
|\varphi^{\pm}(y, \,k) - 
\varphi^{\pm}(y, \,k_0) \big|
 \, dy.    \qquad \qquad
\label{eqn:9-e-8}
\end{eqnarray}
The integral in (\ref{eqn:9-e-7}) is 
estimated by
\begin{eqnarray}
\Big\{\, 
\int_{{\mathbb R}^3} \Big|  \,  
1_{E(x_0, \, 2r)}(y) \,
\frac{\,|V(y)|^2 \langle y \rangle^{2s}}{\,|x_0 - y|^2\,} 
\times        
 \qquad \qquad\qquad\qquad
               \quad\qquad \qquad \nonumber \\
\noalign{\vskip 4pt}
\Big| \, 
|k|e^{\mp i|k| |x_0 - y|}
- |k_0|e^{\mp i|k_0| |x_0 - y|} 
\,\Big|^2   dy  \Big\}^{1/2}
\, \Vert 
\varphi^{\pm}(\cdot, \,k)  
\Vert_{L^{2, -s}}
\label{eqn:9-e-9}
\end{eqnarray}
In view of (\ref{eqn:9-1-before}), it
follows that
\begin{eqnarray*}
1_{E(x_0, \, 2r)}(y) \,
\frac{\,|V(y)|^2 \langle y\rangle^{2s}}{\,|x_0 - y|^2\,} 
\in L^1({\mathbb R}^3_y).
\end{eqnarray*}
Therefore, applying the Lebesgue dominated convergence
theorem to the integral in 
(\ref{eqn:9-e-9}) and appealing to 
Lemma \ref{lem:9-a}, we see 
that the integral in 
(\ref{eqn:9-e-7}) tends to $0$ as $k$ 
approaches $k_0$.
Also, the integral in (\ref{eqn:9-e-8}) is 
estimated by
\begin{eqnarray*}
\Big\{ 
\int_{{\mathbb R}^3} \Big|  \,  
1_{E(x_0, \, 2r)}(y) \,
\frac{\,|V(y)|^2 \langle y \rangle^{2s}}{\,|x_0 - y|^2\,} 
\,\Big|^2   dy  \Big\}^{1/2}    \times  
                                \qquad\qquad\nonumber\\
\noalign{\vskip 4pt}
\, \Vert 
\varphi^{\pm}(\cdot, \,k)
-    \varphi^{\pm}(\cdot, \,k_0)   
\Vert_{L^{2, -s}},
\end{eqnarray*}
which tends to $0$, 
by Lemma \ref{lem:9-a}, as
$k$ approaches $k_0$. 
Thus, we have shown that   
\begin{equation}
\limsup_{(x,k) \to (x_0, k_0)}
| I\!I_{0\kappa}^{\pm}(k) | 
\le 4 \varepsilon .
\label{eqn:9-e-10}
\end{equation}
Combining (\ref{eqn:9-e-3}), (\ref{eqn:9-e-6})
and (\ref{eqn:9-e-10}), we deduce that
\begin{equation}
\limsup_{(x,k) \to (x_0, k_0)}
| \psi_{0\kappa}^{\pm}(x, \, k)
- \psi_{0\kappa}^{\pm}(x_0, \, k_0) | 
\le 6 \varepsilon .
\label{eqn:9-e-11}
\end{equation}
Since $\varepsilon$ is arbitrary,
(\ref{eqn:9-e-11}) implies (\ref{eqn:9-e-1}).
\hfill$\square$

\vspace{15pt}
\begin{lem}\label{lem:9-f}
 $\psi_{0\mu}^{\pm}(x, \, k)$ are
continuous  on
${\mathbb R}^3_x \times \big\{ \; k \; | \; |k| 
\in (0, \, +\infty) \setminus
   \sigma_p(H) \; \big\}$.
\end{lem}

{\it Proof.}
The proof is similar to that of Lemma \ref{lem:9-e}.

Let $(x_0, \, k_0)$ be an arbitrary
point in 
${\mathbb R}^3_x \times \big\{ \; k \; | \; |k| 
\in (0, \, +\infty) \setminus
   \sigma_p(H) \; \big\}$.
We shall show that
\begin{equation}
\psi_{0\mu}^{\pm}(x, \, k) 
 \to 
\psi_{0\mu}^{\pm}(x_0, \, k_0)
\quad 
\hbox{\rm as }{(x, \, k)\to (x_0, \, k_0)}.
\label{eqn:9-f-1}
\end{equation}
To show this, we first need to appeal to
the definition (\ref{eqn:3q}) of $m_{\lambda}(x)$
and the inequality (\ref{eqn:4e}). 
We then have
\begin{equation}
|m_{\lambda}(x)| \le 
\frac{\, \hbox{const.}\,}{2 \pi^2}
\cdot \frac{\lambda}{\, |x| \,} 
(1 + \lambda |x|)^{-1}
\le
\hbox{const.}^{\!\prime} \frac{\lambda}{\, |x| \,}, 
\label{eqn:9-f-2} 
\end{equation}
where const. is the same as in (\ref{eqn:4e}) and
$\hbox{const.}^{\!\prime} := \hbox{const.}/2\pi^2$.

Let $\varepsilon > 0$ be given.
We choose $r>0$ so that
\begin{equation}
\hbox{const.}^{\!\prime} \, |k_0| \,
\Vert V \Vert_{L^{\infty}} 
\Big\{\!\!
\sup_{y \in {\mathbb R}^3 
    \atop
 |k - k_0| \le r} 
\!\!\!
 | \varphi^{\pm}(y, \, k) |
\Big\}
\int_{|y| \le 2r} 
\frac{1}{\, |y| \,} \,dy < \varepsilon.
\label{eqn:9-f-3}
\end{equation}
Similarly to (\ref{eqn:9-e-3}), we write
\begin{eqnarray}
&&\psi_{0\mu}^{\pm}(x, \, k) 
- 
\psi_{0\mu}^{\pm}(x_0, \, k_0)  \nonumber\\
\noalign{\vskip 3pt}
&=& \{ \psi_{0\mu}^{\pm}(x, \, k) 
- 
\psi_{0\mu}^{\pm}(x_0, \, k) \}   
+
\{ \psi_{0\mu}^{\pm}(x_0, \, k) 
- 
\psi_{0\mu}^{\pm}(x_0, \, k_0) \}  
             \label{eqn:9-f-4}\qquad  \\
\noalign{\vskip 3pt}
&=:&
 I^{\pm}_{0\mu}(x, \, k)
+ I\!I^{\pm}_{0\mu}(k).    \nonumber
\end{eqnarray}
If $|x - x_0| \le r$ and $|k - k_0| \le r$,
then it follows from
 the definition (\ref{eqn:4a+})  and
(\ref{eqn:9-f-2}),
(\ref{eqn:9-f-3}) that
\begin{eqnarray}
|I^{\pm}_{0\mu}(x, \, k)| 
\le
2 \frac{|k|}{\, |k_0| \,} \varepsilon  
+ 
\Big\{\!\!
\sup_{y \in {\mathbb R}^3 
    \atop
 |k - k_0| \le r} 
\!\!\!
 | \varphi^{\pm}(y, \, k) |
\Big\} \times    \qquad \qquad \qquad\qquad \nonumber \\
\times \int_{{\mathbb R}^3} \Big|    
\Big\{
1_{E(x, \, 2r)}(y) \,
m_{|k|}(x - y) 
-
1_{E(x_0, \, 2r)}(y) \,
m_{|k|}(x_0 - y)
\Big\} V(y)
\Big| \, dy.
\label{eqn:9-f-5}
\end{eqnarray}
Noting (\ref{eqn:9-f-2}) and (\ref{eqn:9-e-4-1}),
we find that
the integrand in (\ref{eqn:9-f-5}) is 
bounded by the function
\begin{equation}
\frac{\, 5 \,}{2} \times 
1_{E(x_0, \, r)}(y)
\frac{\,\hbox{const.}^{\!\prime}(|k_0| + r)\,}
{\, |x_0 - y| \,} \, |V(y)|
\in L^1({\mathbb R}^3_y)
\label{eqn:9-f-6}
\end{equation}
for all $(x, \, k)$ with
$|x - x_0| \le r$, $|k - k_0| \le r$.
Therefore, the Lebesgue dominated convergence
theorem applied to the integral in 
(\ref{eqn:9-f-5}) gives
\begin{equation}
\limsup_{(x,k) \to (x_0, k_0)}
| I_{0\mu}^{\pm}(x, \, k) | \le 2 \varepsilon .
\label{eqn:9-f-7}
\end{equation}
Here we have used the fact that
$m_{|k|}(x)$ is continuous on
$\big\{ \,{\mathbb R}^3_x \setminus \{ 0 \} \}
\times {\mathbb R}^3_k$.
In a similar manner to (\ref{eqn:9-f-5}),
if $|x - x_0| \le r$ and $|k - k_0| \le r$,
 then we have
\begin{eqnarray}
|I\!I^{\pm}_{0\mu}(k)| 
\le
\Big(
\frac{\,|k|\,}{\, |k_0| \,} + 3 
\Big)
\varepsilon     
\qquad\qquad \qquad\qquad 
    \qquad \qquad \qquad\qquad \qquad \nonumber \\
+ 
\int_{{\mathbb R}^3} \Big|  \,  
1_{E(x_0, \, 2r)}(y) 
\big( 
  m_{|k|}(x_0 -y ) -  m_{|k_0|}(x_0 -y )
\big)
      \Big| \,
\big| V(y) \varphi^{\pm}(y, \,k) \big|
 \, dy  \qquad                         
                  \label{eqn:9-f-8}\\
\noalign{\vskip 4pt}
+ \; 
\int_{{\mathbb R}^3} \Big|  \,  
1_{E(x_0, \, 2r)}(y)   \,
m_{|k_0|}(x_0 - y)
V(y)  \Big| 
|\varphi^{\pm}(y, \,k) - 
\varphi^{\pm}(y, \,k_0) \big|
 \, dy.   \qquad
\label{eqn:9-f-9}
\end{eqnarray}
The integral in (\ref{eqn:9-f-8}) is 
estimated by
\begin{eqnarray}
\Big\{
\int_{{\mathbb R}^3}     
1_{E(x_0, \, 2r)}(y)   \,
\big| 
 \big( 
m_{|k|}(x_0 -y ) -  m_{|k_0|}(x_0 -y ) 
\big)
 \,V(y) 
\langle y \rangle^{s} \, \big|^2
\, dy \Big\}^{1/2}  \times  \nonumber\\
      \times \,
\Vert  \varphi^{\pm}(\cdot, \,k) \Vert_{L^{2, -s}}.
\qquad\qquad
\label{eqn:9-f-10}
\end{eqnarray}
In view of  the inequality (\ref{eqn:4e+}) and 
the continuity of $m_{|k|}(x)$, as
mentioned after (\ref{eqn:9-f-7}), 
we can apply the Lebesgue dominated convergence
theorem to the integral in (\ref{eqn:9-f-10}),
and deduce that
the integral in (\ref{eqn:9-f-8}) tends to
 $0$ as $k$ approaches $k_0$.
Also, the integral in (\ref{eqn:9-f-9}) is 
estimated by
\begin{eqnarray}
\Big\{
\int_{{\mathbb R}^3}     
1_{E(x_0, \, 2r)}(y)   \,
\big| 
   m_{|k_0|}(x_0 -y )
\, V(y) 
\langle y \rangle^{s}
\big|^2 \, dy \Big\}^{1/2}  \times  
                       \qquad\nonumber\\
      \times  \,
\Vert  
\varphi^{\pm}(\cdot, \,k) 
- \varphi^{\pm}(\cdot, \,k_0) 
\Vert_{L^{2, -s}},
\label{eqn:9-f-11}
\end{eqnarray}
which tends to $0$, by Lemma \ref{lem:9-a} and
the inequality (\ref{eqn:4e+}), 
as $k$ approaches $k_0$. Thus we have
shown that 
\begin{equation*}
\limsup_{(x,k) \to (x_0, k_0)}
| I\!I_{0\mu}^{\pm}(x, \, k) | 
\le 4 \varepsilon .
\end{equation*}
By the same arguments as in the end of the
proof of Lemma \ref{lem:9-e}, 
we conclude that (\ref{eqn:9-f-1}) is verified.
\hfill$\square$

\vspace{15pt}
\begin{lem}\label{lem:9-g}
 $\psi_1^{\pm}(x, \, k)$ are
continuous  on
${\mathbb R}^3_x \times \big\{ \; k \; | \; |k| 
\in (0, \, +\infty) \setminus
   \sigma_p(H) \; \big\}$.
\end{lem}

The proof of Lemma \ref{lem:9-g} is 
similar to those of
Lemmas \ref{lem:9-e} and \ref{lem:9-f}.
Actually it is much easier because
the integral kernel of the
operator $G_0$ is independent of the 
variable $k$ (recall the definitions
 (\ref{eqn:2e})  and 
(\ref{eqn:9-2-1})).
For this reason, we omit the proof of 
Lemma \ref{lem:9-g}.

\vspace{15pt}

{\bf Proof of assertion(ii) of 
Theorem \ref{thm:gefcont} } Assertion(ii) is a direct
consequence of (\ref{eqn:9-3-1}) and 
Lemmas \ref{lem:9-e}, \ref{lem:9-f} and \ref{lem:9-g}.
\hfill$\square$

\newpage
\section{Asymptotic behaviors of the generalized eigenfunctions}

We shall first show that the generalized
eigenfunctions $\varphi^{\pm}(x, \, k)$,
defined by (\ref{eqn:8-2}), are
distorted plane waves, and give 
estimates of the differences between
 $\varphi^{\pm}(x, \, k)$ and the
plane wave 
$\varphi_0(x, \, k)=e^{ix\cdot k}$ 
(Theorem \ref{thm:dpw}).
We shall next prove that 
$\varphi^{\pm}(x, \, k)$ are asymptotically
equal to the sums of the plane wave and
the spherical waves $e^{\mp i |x||k|}/|x|$
under the  assumption that 
$\sigma >3$, and shall 
give estimates of the differences 
between $\varphi^{\pm}(x, \, k)$ and
the sums mentioned above (Theorem \ref{thm:sw}).

In view of the definition (\ref{eqn:8-2}) and
Theorem \ref{thm:RC}(ii), it
is clear that $\varphi^-(x, \, k)$
(resp. $\varphi^+(x, \, k)$)
 is the sum of the plane wave $e^{ix\cdot k}$
and the solution of the equation 
(\ref{eqn:6-RC+v}) with the outgoing
radiation condition (\ref{eqn:5-0a})
 \ (resp. the incoming radiation 
condition (\ref{eqn:5-0b})). 
However, the radiation conditions (\ref{eqn:5-0a})
and (\ref{eqn:5-0b}) are generalizations of 
the radiation condition mentioned in the 
beginning of Section 6, and this generalization
makes it unclear that
$$
R^{\mp}(|k|)
 \{ V(\cdot) \varphi_0(\cdot, \, k) \}(x)
$$
behave as $e^{\mp i |x||k|}/|x|$ at
infinity.
Theorem \ref{thm:sw} shows that
this is indeed the case if $\sigma >3$.

\vspace{15pt}
\begin{thm} \label{thm:dpw}
Let $\sigma >2$. If 
$|k| \in (0, \, +\infty) \setminus \sigma_p (H)$,  then
\begin{eqnarray}
\big| 
\varphi^{\pm}(x, \, k) - 
e^{x\cdot k}
\big| 
 \le C_k 
\begin{cases}  
\langle x \rangle^{-(\sigma -2)}  
&\text{\it if $\,\;  2 < \sigma < 3$,}  \\
{}    &\\
\langle x\rangle^{-1} 
\log( 1 + \langle x\rangle)  
       &\text{\it if $\,\;   \sigma = 3$,}  \\
{}    &\\
\langle x\rangle^{-1}   
       &\text{\it if $\,\;   \sigma > 3$,} 
\end{cases}
\label{eqn:10-1}
\end{eqnarray}
where the constant $C_k$ is uniform for $k$ 
in  any compact subset of 
$\{ \, k \, | \, |k| 
\in (0, \, +\infty )\setminus 
\sigma_p(H) \, \}$.  
\end{thm}

{\it Proof.}
In view of (\ref{eqn:9-3-1}), it is
sufficient to show that
all of $\psi_{0\kappa}^{\pm}(x, \, k)$,
$\psi_{0\mu}^{\pm}(x, \, k)$ and 
$\psi_1^{\pm}(x, \, k)$ satisfy the
estimates (\ref{eqn:10-1}).

By  assertion(i) of Theorem \ref{thm:gefcont} and
the definitions (\ref{eqn:psi-0K})
and (\ref{eqn:4a}), we have
\begin{equation}
|\psi_{0\kappa}^{\pm}(x, \, k)|
\le 
\frac{|k|}{\, 2 \pi \,}
\Vert 
\langle \cdot \rangle^{\sigma}
V(\cdot) \varphi^{\pm}(\cdot, \, k)
\Vert_{L^{\infty}} 
\int_{{\mathbb R}^3}  
\frac{1}{\, |x - y|\langle y \rangle^{\sigma} \,}
\, dy.
\label{eqn:10-2}
\end{equation}
If we apply Lemma A.1 in Appendix, with
$n=3$, $\beta=1$ and 
$\gamma = \sigma$, to the integral on the right
hand side of (\ref{eqn:10-2}), we can deduce 
from 
assertion(i) of Theorem \ref{thm:gefcont} 
and(\ref{eqn:10-2}) that
$\psi_{0\kappa}^{\pm}(x, \, k)$ satisfy the
desired estimates.
By the definitions (\ref{eqn:psi-0M}),
(\ref{eqn:4a+}) and
the inequality (\ref{eqn:4e+}), we get
\begin{equation}
|\psi_{0\mu}^{\pm}(x, \, k)|
\le 
C_{|k|}
\Vert 
\langle \cdot \rangle^{\sigma}
V(\cdot) \varphi^{\pm}(\cdot, \, k)
\Vert_{L^{\infty}} 
\int_{{\mathbb R}^3}  
\frac{1}{\, |x - y|^2 \langle y \rangle^{\sigma} \,}
\, dy,
\label{eqn:10-3}
\end{equation}
where the constant $C_{|k|}$ is the 
one specified in (\ref{eqn:4e+}).
Similarly, by the definition 
(\ref{eqn:2e}), we obtain
\begin{equation}
|\psi_1^{\pm}(x, \, k)|
\le 
\frac{1}{\, 2 \pi^2 \,}
\Vert 
\langle \cdot \rangle^{\sigma}
V(\cdot) \varphi^{\pm}(\cdot, \, k)
\Vert_{L^{\infty}} 
\int_{{\mathbb R}^3}  
\frac{1}{\, |x - y|^2 
\langle y \rangle^{\sigma} \,}
\, dy.
\label{eqn:10-4}
\end{equation}
Lemma A.1 with $n=3$, $\beta=2$ and 
$\gamma = \sigma$ now gives
\begin{eqnarray}
|\psi_{0\mu}^{\pm}(x, \, k)|+
| \psi_1^{\pm}(x, \, k)| 
 \le C_k ^{\prime}
\begin{cases}  
\langle x \rangle^{-(\sigma -1)}  
&\text{ if $\,\;  2 < \sigma < 3$,}  \\
{}    &\\
\langle x\rangle^{-2} 
\log( 1 + \langle x\rangle)  
       &\text{ if $\,\;   \sigma = 3$,}  \\
{}    &\\
\langle x\rangle^{-2}   
       &\text{ if $\,\;   \sigma > 3$,} 
\end{cases}
\label{eqn:10-5}
\end{eqnarray}
where the constant $C_k^{\prime}$ is 
uniform for $k$ 
in  any compact subset of 
$\{ \, k \, | \, |k| 
\in (0, \, +\infty )\setminus 
\sigma_p(H) \, \}$.  
\hfill$\square$

\vspace{20pt}

\begin{thm} \label{thm:sw}
Let $\sigma >3$, and suppose that 
$|k| \in (0, \, +\infty) \setminus \sigma_p (H)$.  
Then for $|x| \ge 1$ we have
\begin{eqnarray}
\lefteqn{   \Big| \varphi^{\pm}(x, \, k) - 
\Big( e^{ix\cdot k}  + 
   \frac{e^{\mp i |k| \, |x|}}{\,|x| \,}\, 
f^{\pm}(|k|, \, \omega_x, \, \omega_k) \, \Big)
\Big|                                      }   \nonumber\\
\noalign{\vskip 5pt}
  & & 
\le C_k  \left\{  
\begin{array}{ll}
|x|^{-(\sigma -1)/2}  & {\it if }\; 
 3 < \sigma < 5,              \\
{}  &  {} \\
|x|^{-2} \log (1 + |x| )   
       & {\it if }\;   \sigma = 5,   \\
{}  & {}  \\
|x|^{-2}   
       & {\it if }\;   \sigma > 5,     
\end{array}\right.      \label{eqn:10-6}
\end{eqnarray}
where $\omega_x = x/|x|$, $\omega_k = k/|k|$,
\begin{equation}
f^{\pm}(\lambda , \, \omega_x, \, \omega_k) = 
 - \frac{\lambda}{\, 2 \pi \,}  
\int_{{\mathbb R}^3} 
e^{\pm i \lambda \omega_x \cdot y} \, V(y) \,
  \varphi^{\pm}(y, \, \lambda \omega_k) \, dy,
\label{eqn:10-6-1}
\end{equation}
and the constant $C_k$ is uniform 
for $k$ 
in  any compact subset of 
$\{ \, k \, | \, |k| 
\in (0, \, +\infty )\setminus 
\sigma_p(H) \, \}$.
\end{thm}

\vspace{15pt}

We shall give a proof of Theorem \ref{thm:sw}
by means of a series of lemmas.

\vspace{15pt}

\begin{lem}  \label{lem:10-1}
Let $\sigma >3$. Then
$$
\big| \, \varphi^{\pm}(x, \, k) -  
\big( 
 e^{ix\cdot k} + 
\psi_{0\kappa}^{\pm}(x, \, k)   
 \big)  \, \big|
  \le C_k \langle x \rangle^{-2},
$$
where $C_k$ is a constant uniform 
for $k$ 
in  any compact subset of 
$\{ \, k \, | \, |k| 
\in (0, \, +\infty )\setminus 
\sigma_p(H) \, \}$.
\end{lem}

{\it Proof.}
The lemma is a direct consequence of (\ref{eqn:10-5}) 
and (\ref{eqn:9-3-1}).
\hfill$\square$

\vspace{15pt}

In view of Lemma \ref{lem:10-1}, it is apparent that
we need to evaluate the differences

\begin{eqnarray}
\psi_{0\kappa}^{\pm}(x, \, k) 
- \frac{e^{\mp i |k| \, |x|}}{\,|x| \,}\, 
f^{\pm}(|k|, \, \omega_x, \, \omega_k),
\label{eqn:10-7}
\end{eqnarray}
which are equal to
\begin{eqnarray}
\frac{\, |k| \,}{2\pi}
\int_{{\mathbb R}^3} \Big\{
\frac{\, e^{\mp i |k|(|x| - \omega_x \cdot y)}\, }{|x|} 
-
\frac{\, e^{\mp i |k||x-y|}\, }{|x-y|}
\Big\}
\, V(y) \,
  \varphi^{\pm}(y, \, |k| \omega_k) \, dy
\label{eqn:10-8}
\end{eqnarray}
by (\ref{eqn:psi-0K}), (\ref{eqn:4a}) and 
(\ref{eqn:10-6-1}). Thus we are led to
consider the following integrals:
\begin{eqnarray}
\frac{1}{\, |x| \,}
  \int_{{\mathbb R}^3} 
e^{ia(|x| - \omega_x \cdot y)}
  u(y) \, dy ,    \label{eqn:10-9-a}\\
\noalign{\vskip 4pt}
\int_{{\mathbb R}^3}   
\frac{\, e^{ia|x- y|} \,}{|x - y|} \,
  u(y) \, dy ,         \label{eqn:10-9-b}
\end{eqnarray}
 and 
their difference.
The same integrals as in (\ref{eqn:10-9-a}) and
(\ref{eqn:10-9-b}) were discussed in
Ikebe\cite[\S 3]{Ikebe}, though our 
arguments below are slightly different
from those of \cite{Ikebe}, and 
our estimates are slight refinements of 
those of \cite{Ikebe}. 

\vspace{15pt}

\begin{lem}  \label{lem:10-2}
Let $a \in \mathbb R$ and
let $u$ satisfy 
\begin{equation}
| u(x)| \le C 
     \langle x \rangle^{-\sigma},
\quad \sigma > 3.
\label{eqn:lem10-2-1}
\end{equation}
Then for $|x| \ge 1$ we have
\begin{eqnarray}
\Big|
  \int_{|y| \ge \sqrt{|x|}} 
e^{ia(|x| - \omega_x \cdot y)}
  u(y) \, dy   
   \Big|     
\le C_1 \, 
\Vert 
  \langle \cdot \rangle^{\sigma}u 
\Vert_{L^{\infty}}
\, |x|^{-(\sigma -3)/2} ,  
            \label{eqn:lem10-2-2}\\
\noalign{\vskip 4pt}
\Big|
\int_{|y| \ge \sqrt{|x|}}   
\frac{\, e^{ia|x- y|} \,}{|x - y|} \,
  u(y) \, dy   \Big|    
\le C_2  \, 
\Vert 
  \langle \cdot \rangle^{\sigma}u 
\Vert_{L^{\infty}}
\,
|x|^{-(\sigma -1)/2},   
               \label{eqn:lem10-2-3}
\end{eqnarray}
where the constants $C_1$  and $C_2$
 are independent
of $\,a$.
\end{lem}

{\it Proof.}
It follows that
\begin{eqnarray}
\Big|
  \int_{|y| \ge \sqrt{|x|}} 
e^{ia(|x| - \omega_x \cdot y)}
  u(y) \, dy   
   \Big|   
&\le& 
\int_{|y| \ge \sqrt{|x|}}  
\Vert 
  \langle \cdot \rangle^{\sigma}u 
\Vert_{L^{\infty}}
\, \, |y|^{-\sigma} \, dy    \nonumber\\
\noalign{\vskip 3pt}
&\le&
 C_1 
\Vert 
  \langle \cdot \rangle^{\sigma}u 
\Vert_{L^{\infty}}
\,
\, |x|^{-(\sigma -3)/2}.
\label{eqn:lem10-2-4}
\end{eqnarray}

To show (\ref{eqn:lem10-2-3}), we decompose the
integral in (\ref{eqn:lem10-2-3})
into two parts:
\begin{eqnarray}
\lefteqn{\int_{|y| \ge \sqrt{|x|}}   
\frac{\, e^{ia|x- y|} \,}{|x - y|} \,
  u(y) \, dy
=} \qquad\qquad   \nonumber \\
\noalign{\vskip 4pt}
&&\Big\{
\int_{F_0(x)} +  \int_{F_1(x)}
\Big\}
\frac{\, e^{ia|x- y|} \,}{|x - y|} \,
  u(y) \, dy,
\label{eqn:lem10-2-5}
\end{eqnarray}
where
\begin{eqnarray*}
F_0(x) := \{ \, y \in {\mathbb R}^3 \,\, |
 \,\; |y| \ge \sqrt{|x|}, \;\,
   |x -y| \le \frac{\, |x| \,}{2} \, \},  \\
F_1(x) := \{ \, y \in {\mathbb R}^3 \,\, |
 \,\; |y| \ge \sqrt{|x|}, \;\,
   |x -y| \ge \frac{\, |x| \,}{2} \, \}.
\end{eqnarray*}
If $y \in F_0(x)$, then
\begin{eqnarray*}
|y|=|x - (x - y)| \ge |x| - |x - y| \ge 
    \frac{\, |x| \,}{2},
\end{eqnarray*}
hence we have
\begin{eqnarray}
\lefteqn{ \Big| 
\int_{F_0(x)} 
\frac{\, e^{ia|x- y|} \,}{|x - y|} \,
           u(y) \, dy  \Big|  } \qquad\nonumber\\
&\le& 
\int_{F_0(x)}  \frac{1}{\, |x - y| \,}
 \, 
\Vert 
  \langle \cdot \rangle^{\sigma}u 
\Vert_{L^{\infty}}
\,
|y|^{-\sigma}   \, dy   \nonumber \\
\noalign{\vskip 4pt}
&\le&
\Vert 
  \langle \cdot \rangle^{\sigma}u 
\Vert_{L^{\infty}}
\,
 2^{\sigma}|x|^{-\sigma} \!
\int_{|x-y| \le |x|/2}  
  \frac{1}{\, |x - y| \,} \, dy   \nonumber \\
\noalign{\vskip 3pt}
 &=& C^{\prime} \, 
\Vert 
  \langle \cdot \rangle^{\sigma}u 
\Vert_{L^{\infty}}
\,
|x|^{-(\sigma -2)}.
\label{eqn:lem10-2-6}
\end{eqnarray}
If $y \in F_1(x)$, then $|x-y| \ge |x|/2$, therefore 
we get
\begin{eqnarray}
\lefteqn{\Big| 
\int_{F_1(x)} 
\frac{\, e^{ia|x- y|} \,}{|x - y|} \,
  u(y) \, dy  \Big| }  \qquad  \nonumber\\
&\le& 
\int_{|y| \ge \sqrt{|x|}  } 
\frac{2}{\, |x| \,} \,
\Vert 
  \langle \cdot \rangle^{\sigma}u 
\Vert_{L^{\infty}}
\, 
  |y|^{-\sigma} \, dy     \nonumber \\
\noalign{\vskip 3pt}
&\le&  C^{\prime\prime} \,
\Vert 
  \langle \cdot \rangle^{\sigma}u 
\Vert_{L^{\infty}}
\,
|x|^{-(\sigma -1)/2}. 
\label{eqn:lem10-2-7}
\end{eqnarray}
Since $\sigma -2 > (\sigma -1)/2$, 
we conclude from 
(\ref{eqn:lem10-2-5})--(\ref{eqn:lem10-2-7})
that the inequality (\ref{eqn:lem10-2-3}) holds.
\hfill$\square$

\vspace{20pt}

In view of (\ref{eqn:10-9-a}), (\ref{eqn:10-9-b}) 
and Lemma \ref{lem:10-2}, we now need to consider 
the integral
\begin{equation}
  \int_{|y| \le \sqrt{|x|}} 
\Big(
\frac{1}{\, |x| \,} \,
e^{ia(|x| - \omega_x \cdot y)}
-
\frac{\, e^{ia|x- y|} \,}{|x - y|} \,
\Big)
  u(y) \, dy .
\label{eqn:10-d-in}
\end{equation}
To  get an estimate on the integral 
(\ref{eqn:10-d-in}), 
we split it into two parts:
\begin{eqnarray}
\frac{1}{\, |x| \,} 
\int_{|y| \le \sqrt{|x|}} 
\Big(
e^{ia(|x| - \omega_x \cdot y)}
-
e^{ia|x-y|}
\Big) \, u(y) \, dy     
\label{eqn:10-d-in-1}   \qquad\qquad  \\
\noalign{\vskip 4pt}
+ \int_{|y| \le \sqrt{|x|}} 
\Big(
\frac{1}{\, |x| \,}
-
\frac{1}{\, |x-y| \,}
\Big) 
e^{ia|x-y|}
\, u(y) \, dy, 
\label{eqn:10-d-in-2}
\end{eqnarray}
and evaluate these two integrals separately.

\vspace{15pt}

\begin{lem}  \label{lem:10-3}
Under the same assumptions 
as in Lemma \ref{lem:10-2},
 we have
\begin{eqnarray}
\lefteqn{\Big|
\frac{1}{\, |x| \,}
  \int_{|y| \le \sqrt{|x|}} 
\Big(
e^{ia(|x| - \omega_x \cdot y)}
-
e^{ia|x-y|}
\Big) \, u(y) \, dy 
\Big|  }   \nonumber\\
\noalign{\vskip 5pt}
  & & 
\le C_3\,|a|  \, 
\Vert 
  \langle \cdot \rangle^{\sigma}u 
\Vert_{L^{\infty}}
\,
\left\{  
\begin{array}{ll}
|x|^{-(\sigma -1)/2}  & {\it if }\; 
 3 < \sigma < 5,              \\
{}  &  {} \\
|x|^{-2} \log (1 + |x| )   
       & {\it if }\;   \sigma = 5,   \\
{}  & {}  \\
|x|^{-2}   
       & {\it if }\;   \sigma > 5     
\end{array}\right.      \label{eqn:10-3-1}
\end{eqnarray}
for $|x| \ge 1$,
where the constant $C_3$  is 
independent of  $a$.
\end{lem}

{\it Proof.}
We start with  simple remarks that
\begin{equation}
|x - y| = |x| \,
\Big(
1 -2 \, \frac{\, \omega_x \cdot y \,}{|x|}
 + \frac{|y|^2}{\, |x|^2 \,}
\Big)^{\!1/2}
\label{eqn:10-3-2}
\end{equation}
and
\begin{equation}
\big| (1+ \rho)^{1/2} - 
(1 + \frac{\rho}{\,2\,}) \big|
\le \frac{\, \sqrt{2} \,}{2} \, \rho^2,
\qquad  \, \rho \ge -\frac{1}{\, 2 \,}.
\label{eqn:10-3-3}
\end{equation}
It is easy to see that 
\begin{equation}
 \big| 
-2 \, \frac{\, \omega_x \cdot y \,}{|x|}
 + \frac{|y|^2}{\, |x|^2 \,}
\big| \le \frac{1}{\, 2 \,}
\label{eqn:10-3-4}
\end{equation}
if $\sqrt{|x|} \ge 5$ and
$|y| \le \sqrt{|x|}$. 
Hence, it follows
from 
(\ref{eqn:10-3-2})--(\ref{eqn:10-3-4}) that
\begin{equation}
 \big| 
|x-y| -
( |x| -  \omega_x \cdot y )
 \big| 
\le 3\sqrt{2} \, \frac{\; |y|^2 }{|x|}
\label{eqn:10-3-5}
\end{equation}
when $\sqrt{|x|} \ge 5$ and
$|y| \le \sqrt{|x|}$. 
Using the inequality
\begin{eqnarray*}
| e^{i \alpha} -  e^{i \alpha} |
\le | \alpha - \beta |,  
\qquad \, \alpha, \, \beta  \in {\mathbb R},
\end{eqnarray*}
we have
\begin{eqnarray}
\Big|
  \int_{|y| \le \sqrt{|x|}} 
\big(
e^{ia(|x| - \omega_x \cdot y)}
-
e^{ia|x-y|} \,
\big) \, u(y) \, dy 
\Big|     \qquad\qquad\qquad  \nonumber\\
\noalign{\vskip 4pt}
\le 
  \int_{|y| \le \sqrt{|x|}} 
\big| a(|x| - \omega_x \cdot y) 
- a|x-y| \big|
\, 
\Vert 
  \langle \cdot \rangle^{\sigma}u 
\Vert_{L^{\infty}}
\,
 \langle y\rangle^{-\sigma} 
\, dy \,   
\label{eqn:10-3-6}      \\
\noalign{\vskip 4pt}
\le  3 \sqrt{2} \, 
\, |a| \, 
\Vert 
  \langle \cdot \rangle^{\sigma}u 
\Vert_{L^{\infty}}
\,
\frac{1}{\, |x| \,}  \,
\int_{|y| \le \sqrt{|x|}} \, |y|^2 \,
\langle y\rangle^{-\sigma} \, dy
\qquad\qquad
\label{eqn:10-3-7}
\end{eqnarray}
when $\sqrt{|x|}\ge 5$.
Here  
we have used (\ref{eqn:10-3-5}) 
in the second inequality (\ref{eqn:10-3-7}). 
Now we have
\begin{eqnarray}
\int_{|y| \le \sqrt{|x|}}  \, |y|^2 \,
\langle y\rangle^{-\sigma} \, dy 
\qquad\qquad\qquad \qquad\qquad 
\qquad\qquad \qquad\qquad        \nonumber \\
\noalign{\vskip 4pt}
\le
2^{\sigma /2} \! \int_{|y| \le \sqrt{|x|}} 
\, (1 + |y| )^{2 -\sigma} \, dy  
 \qquad\qquad 
(\because \langle y \rangle 
  \ge 
\frac{1}{\sqrt{2} \,} \, (1+ |y|)\, )  \nonumber \\
=
2^{\sigma /2} \,
\times 4\pi  
\int_0^{\sqrt{|x|}}  \, 
(1 + r)^{-\sigma +4} \, dr
\qquad\qquad \qquad\qquad\quad
             \label{eqn:10-3-8}\\
\noalign{\vskip 6pt}
\le
 2^{(\sigma +4)/2} \pi \times
\left\{  
\begin{array}{ll}
\displaystyle{
\frac{|x|^{-(\sigma -5)/2}}{5-\sigma}    } 
    & \;\; {\rm if }\;\; 
 3 < \sigma < 5,              \\
{}  &  {} \\
\log (1 + |x| )   
       & \;\;{\rm if }\;\;   \sigma = 5,   
\label{eqn:10-3-9}\\
{}  & {}  \\
\displaystyle{\frac{1}{\, \sigma -5 \,} } 
  & \;\;{\rm if }\;\;   \sigma > 5,     
\end{array}\right.     \quad
 \label{eqn:10-3-10}
\end{eqnarray}
where we have used spherical polar coordinates 
in (\ref{eqn:10-3-8}). Combining (\ref{eqn:10-3-10})
with (\ref{eqn:10-3-7}) yields the desired inequalities.
\hfill$\square$

\vspace{15pt}

\begin{lem}  \label{lem:10-4}
Under the same assumptions 
as in Lemma \ref{lem:10-2},
 we have
\begin{eqnarray}
\lefteqn{\Big|
 \int_{|y| \le \sqrt{|x|}} 
\Big(
\frac{1}{\, |x| \,}
-
\frac{1}{\, |x-y| \,}
\Big) 
e^{ia|x-y|}
\, u(y) \, dy 
\Big|  }   \nonumber\\
\noalign{\vskip 5pt}
  & & 
\le C_4\, 
\Vert 
  \langle \cdot \rangle^{\sigma}u 
\Vert_{L^{\infty}}
\,
\left\{  
\begin{array}{ll}
|x|^{-\sigma/2}  & {\it if }\; 
 3 < \sigma < 4,              \\
{}  &  {} \\
|x|^{-2} \log (1 + |x| )   
       & {\it if }\;   \sigma = 4,   \\
{}  & {}  \\
|x|^{-2}   
       & {\it if }\;   \sigma > 4     
\end{array}\right.      \label{eqn:10-4-1}
\end{eqnarray}
for $|x| \ge 1$, 
where the constant $C_4$  is 
independent of  $a$.
\end{lem}

{\it Proof.}
If $\sqrt{|x|}\ge 5$ and
$|y| \le \sqrt{|x|}$, then
the inequality (\ref{eqn:10-3-5}) implies
\begin{eqnarray*}
\big| \, |x - y| - |x| \, \big|
\le |y| + 3\sqrt{2} 
\frac{\; |y|^2 \,}{|x|}.
\end{eqnarray*}
Also, if $\sqrt{|x|}\ge 5$ and
$|y| \le \sqrt{|x|}$, we then have
\begin{eqnarray*}
|x-y| \ge |x| - |y| \ge 
|x| - \frac{\, |x| \,}{5} 
=\frac{\, 4 \,}{5} \,|x| .
\end{eqnarray*}
Using these two inequalities, we
arrive at
\begin{eqnarray}
\Big|
 \int_{|y| \le \sqrt{|x|}} 
\Big(
\frac{1}{\, |x| \,}
-
\frac{1}{\, |x-y| \,}
\Big) 
e^{ia|x-y|}
\, u(y) \, dy 
\Big|        \qquad\quad \nonumber \\
\noalign{\vskip 4pt}
\le
\int_{|y| \le \sqrt{|x|}} 
\frac{\, 5 \,}{4}\cdot
   \frac{1}{\, |x|^2 }\,
\big( |y| +  3\sqrt{2} 
\frac{\; |y|^2 \,}{|x|} \, \big)
\, 
\Vert 
  \langle \cdot \rangle^{\sigma}u 
\Vert_{L^{\infty}}
\, 
\langle y \rangle^{-\sigma} \, dy
     \quad  \nonumber\\
\noalign{\vskip 4pt}
=
\frac{\, 5 \,
\Vert 
  \langle \cdot \rangle^{\sigma}u 
\Vert_{L^{\infty}}
\, 2^{\sigma/2}\,}{4}
\Big(
   \frac{1}{\, |x|^2 }
\int_{|y| \le \sqrt{|x|}} 
( 1 + |y|)^{1 - \sigma} \, dy  
\qquad\qquad \nonumber \\
\noalign{\vskip 4pt}
+ \; 3 \sqrt{2}
\frac{1}{\, |x|^3 }
\int_{|y| \le \sqrt{|x|}} 
( 1 + |y|)^{2 - \sigma} \, dy
\Big)
\label{eqn:10-4-2}
\end{eqnarray}
provided that $\sqrt{|x|}\ge 5$.
By introducing spherical polar
coordinates, we obtain
\begin{eqnarray}
\int_{|y| \le \sqrt{|x|}} 
\, (1 + |y| )^{1 -\sigma} \, dy   
\qquad\qquad\qquad\qquad\quad \nonumber \\
\noalign{\vskip 4pt}
\le
4\pi 
\left\{  
\begin{array}{ll}
\displaystyle{
\frac{|x|^{-(\sigma -4)/2}}{4-\sigma}    } 
    & \;\; {\rm if }\;\; 
 3 < \sigma < 4,              \\
{}  &  {} \\
\log (1 + |x| )   
       & \;\;{\rm if }\;\;   \sigma = 4,   \\
{}  & {}  \\
\displaystyle{\frac{1}{\, \sigma - 4 \,} } 
  & \;\;{\rm if }\;\;   \sigma > 4.     
\end{array}\right.     \quad
 \label{eqn:10-4-3}
\end{eqnarray}
Combining (\ref{eqn:10-4-2}) with
(\ref{eqn:10-4-3}) and (\ref{eqn:10-3-10}),
we conclude that
the desired inequalities are verified.
\hfill$\square$

\vspace{15pt}

{\bf Proof of Theorem \ref{thm:sw} }\  
We write
\begin{eqnarray*}
\varphi^{\pm}(x, \, k) - 
\Big( e^{ix\cdot k}  + 
   \frac{e^{\mp i |k| \, |x|}}{\,|x| \,}\, 
f^{\pm}(|k|, \, \omega_x, \, \omega_k) \, \Big) 
        \qquad\qquad\qquad\qquad\qquad\\
=
\varphi^{\pm}(x, \, k) -  
\big( 
 e^{ix\cdot k} + 
\psi_{0\kappa}^{\pm}(x, \, k)   \big) 
\qquad\qquad\qquad\qquad\qquad\qquad\qquad \\ 
+
\;\big(
\psi_{0\kappa}^{\pm}(x, \, k)  
 - \frac{e^{\mp i |k| \, |x|}}{\,|x| \,}\, 
f^{\pm}(|k|, \, \omega_x, \, \omega_k) \, \big)  
            \qquad\qquad\qquad\qquad\\
=
\varphi^{\pm}(x, \, k) -  
\big( 
 e^{ix\cdot k} + 
\psi_{0\kappa}^{\pm}(x, \, k)   \big) 
\qquad\qquad\qquad\qquad\qquad\qquad\qquad \\
\noalign{\vskip 4pt}
+ \; \frac{\, |k| \,}{2\pi}  \cdot
\frac{1}{\, |x| \,}
\int_{|y|\ge\sqrt{|x|} }   
    e^{\mp i |k|(|x| - \omega_x \cdot y)}         
\, V(y) \,
  \varphi^{\pm}(y, \, |k| \omega_k) \, dy   
                                \qquad\quad\\
\noalign{\vskip 4pt}
- \;  \frac{\, |k| \,}{2\pi} 
\int_{|y|\ge\sqrt{|x|} } 
\frac{\, e^{\mp i |k||x-y|}\, }{|x-y|}
V(y) \,
  \varphi^{\pm}(y, \, |k| \omega_k) \, dy    
                    \qquad\qquad\qquad\quad \\
\noalign{\vskip 4pt}
+ \; \frac{\, |k| \,}{2\pi}  \cdot
\frac{1}{\, |x| \,}
\int_{|y|\le\sqrt{|x|} }   
   \big(
 e^{\mp i |k|(|x| - \omega_x \cdot y)}  
 -
   e^{\mp i |k||x-y|}   
   \big)       
\, V(y) \,
  \varphi^{\pm}(y, \, |k| \omega_k) \, dy   \\
\noalign{\vskip 4pt}
+ \;
\frac{\, |k| \,}{2\pi}
\int_{|y|\le\sqrt{|x|}}
   \Big\{
\frac{1}{\, |x| \,} 
-
\frac{1}{\, |x-y| \,}
\Big\}
\, e^{\mp i |k||x-y|} \,V(y) \,
  \varphi^{\pm}(y, \, |k| \omega_k) \, dy,  \quad
\end{eqnarray*}
where we have used the fact 
that (\ref{eqn:10-7}) equals 
(\ref{eqn:10-8}), and decomposed 
the integral in (\ref{eqn:10-8}) 
into four parts.
Now the conclusion of the theorem
follows from assertion(i) of 
Theorem \ref{thm:gefcont}
and Lemmas 
\ref{lem:10-1}--\ref{lem:10-4}.
\hfill$\square$


\newpage
\section{Appendix}

\vspace{12pt}

In this appendix we shall derive a few 
formulae and estimates concerning 
the cosine integral and 
the sine integral functions for the
reader's convenience, 
the formulae and estimates
which seem not to be found in the literature.
We begin with the definitions of these
functions and some basic facts (cf.
\cite{Erdelyi1} and \cite{Erdelyi2}). 

\vspace{20pt}

{\it A.1. The cosine integral function.} The definition is
$$
\hbox{ci}(\rho) = -\hbox{Ci}(\rho) = 
    \int_\rho^{+\infty}  \frac{\, \cos t \,}{t} \, dt , 
\quad \rho >0
$$
(cf. \cite[p. 386]{Erdelyi1}).  We have

\begin{equation*}
|\hbox{ci}(\rho) | \le  \hbox{const.}
\begin{cases}  
{\rho}^{-1}  &\text{  if $\; \rho \ge 1$,}    \\
{}    &\\
1+|\log \rho |  &\text{  if $\; 0 < \rho \le 1$.}
\end{cases}
\end{equation*}
\vspace{8pt}\noindent
The estimate for $\rho \ge 1$ follows from 
$$
\hbox{ci}(\rho) =  - \frac{\sin \rho}{\, \rho \,} +
 \frac{\cos \rho}{\, \rho^2 \,} 
 -2 \int_\rho^{+\infty}  \frac{\, \cos t \,}{t^3} \, dt, 
$$
which can be shown by repeated use of 
integration by parts. 
The estimate for $0 < \rho \le 1$ follows from 
\cite[p. 145, Formula(6)]{Erdelyi2}.

The cosine integral function $\hbox{ci}(\rho)$ has 
an analytic continuation
$\hbox{ci}(z)$, which is a many-valued function with a 
logarithmic branch-point
at $z=0$ (see \cite[p.145]{Erdelyi2} ). In this paper, 
we choose the principal branch:
$$
\hbox{ci}(z)= -\gamma -{\rm Log} \, z -
  \sum_{m=1}^{\infty} \frac{(-1)^m}{(2m)!2m} z^{2m},
  \quad  
z \in {\bf C }\setminus(-\infty, \, 0],      \eqno(\rm A.1)
$$
where $\gamma$ is Euler's constant and 
$|{\rm Im}\, {\rm Log} \, z| < \pi$.
Note that the power series 
$$
h_e(z) := \sum_{m=1}^{\infty} \frac{(-1)^m}{(2m)!2m} z^{2m}
$$ 
on the right hand side
of (A.1)    is 
an entire function and satisfies that $h_e(-z)=h_e(z)$, i.e.,
$h_e(z)$ is an even function.

\vspace{15pt}

{\it A.2. The sine integral function.} The definition is
$$
\hbox{si}(\rho) = 
    - \int_\rho^{+\infty}  \frac{\, \sin t \,}{t} \, dt , 
\quad \rho >0 
$$
(cf. \cite[p. 386]{Erdelyi1}).
Since
$$
\hbox{si}(\rho) = - \frac{\pi}{\, 2 \,} + 
     \int_0^{\rho}  \frac{\, \sin t \,}{t} \, dt ,
$$
we can show, by integration by parts, that
$$
|\hbox{si}(\rho) | \le  \hbox{const.}
 ( 1 + |\rho| \, )^{-1} .
$$
Moreover, we see that
$\hbox{si}(\rho)$ has an analytic continuation
$\hbox{si}(z)$:
$$
\hbox{si}(z)=- \frac{\pi}{\, 2 \,}
    + \sum_{m=0}^{\infty} \frac{(-1)^m}{(2m+1)! \, 
(2m+1)}\, z^{2m+1}.                          \eqno(\rm A.2)
$$ 
It follows from (A.2) that
$\hbox{si}(z)$ is an entire
function and satisfies that
$$
\hbox{si}(-z)= -\pi - \hbox{si}(z)        \eqno(\rm A.3)                  
$$
(cf. \cite[p.145]{Erdelyi2}).

\vspace{20pt}

{\it A.3. Laplace transforms.} In computing the
resolvent kernel of $\sqrt{-\Delta}$ in Section 2, 
we applied the following formula
$$
\int_0^{+\infty} \! e^{-pt} \frac{\, 1\, }{t^2 + a^2} \, dt
 = - \frac{\, 1 \,}{a} \{\, 
        \hbox{ci}(ap) \sin (ap) + 
   \hbox{si}(ap) \cos(ap)  \, \},           \eqno(\rm A.4)
$$
where $\hbox{Re }p >0, \;a>0$ 
(cf. \cite[p. 269, Formula(46)]{Erdelyi1}).

\vspace{10pt}

For the purpose of applications in Section 2,
it is convenient to replace $p$ in (A.4) 
 with $-z$.  We thus have the
function
\begin{eqnarray*}
\lefteqn{ 
 - \{ \, {\rm ci}(-z)\sin (-z) + {\rm si}(-z) \cos(-z) \} }  \\
&\qquad =& \sin (z) \, {\rm ci}(-z)  - \cos(z) \, {\rm si}(-z)  ,
\end{eqnarray*}
which is holomorphic in ${\mathbb C}\setminus [0, \, + \infty)$.

\vspace{20pt}

{\it A.4. Estimates of a convolution.} We have often encountered
the convolution of the form
$$
\Phi(x):= \int_{{\mathbb R}^n} 
\frac{1}{\, |x-y|^{\beta} 
\langle y \rangle^{\gamma} \,} \, dy
$$
in the previous sections, and 
used Lemma A.1 below several times. Although the
results exhibited in Lemma A.1 are well-known, it 
appears neither
 in a convenient form for our purpose
(see Ikebe \cite{Ikebe}) nor in an
 accessible form
 (see Kuroda \cite{Kuroda} which is written 
in Japanese) in the literature. 
For this reason, we reproduce the results here 
for the reader's convenience.

\vspace{20pt}
Lemma A.1. {\it If $0 < \beta < n$ and $\beta + \gamma >n$,
then
$\Phi (x)$ is a bounded continuous function satisfying
\begin{equation*}
|\Phi (x)| \le  C_{\beta\gamma n}
\begin{cases}  
{\langle x \rangle}^{-(\beta + \gamma -n )}  
      &\text{  if $\; 0 < \gamma < n$,}   \\
{}    &\\
{\langle x \rangle}^{-\beta} 
   \log(1 + \langle x \rangle)
    &\text{  if $\; \gamma = n$, }  \\
{}    &\\
{\langle x \rangle}^{-\beta} 
  &\text{  if $ \; \gamma >n$,  }
\end{cases}
\end{equation*}
where $C_{\beta\gamma n}$ is a constant depending
on $\beta$, $\gamma$ and $n$.}

\vspace{15pt}

We shall divide the proof into
four steps.

\vspace{15pt}

{\it Step 1. \ $\Phi (x)$ is a continuous function 
on ${\mathbb R}^n$.}

\vspace{5pt}

{\it Proof.}
Let $x_0$ be an arbitrary point in ${\mathbb R}^n$,
and let $\varepsilon >0$ be given.
Since $0 < \beta <n$, we can choose 
$r>0$ so that%
$$
 \int_{|y| \le 2r}
\frac{1}{\, |y|^{\beta} \,} \, dy
 < \varepsilon .           \eqno(\rm A.5)
$$
We then decompose $\Phi(x)$ into two
parts:
$$
\Phi(x)
= \Big( 
\int_{B(x, \,2r)} 
+ \int_{E(x,\,2r)}  
\Big)
\frac{1}{\, |x-y|^{\beta} 
\langle y \rangle^{\gamma} \,} \, dy   
=: \Phi_b(x) + \Phi_e(x),     \eqno(\rm A.6)
$$
where 
$B(x, \, 2r) = \{ \, y \, | \, |x - y| \le 2r \, \}$
and $E(x, \, 2r)$ is the same as 
in the proof of Lemma \ref{lem:9-e}. 
By (A.5), we get
$$
0 < \Phi_b(x) < \varepsilon  \eqno(\rm A.7)
$$
for all $x \in {\mathbb R}^n$.
It follows from the 
definition of $\Phi_e(x)$ that
$$
\displaylines{
\quad \Phi_e(x) - \Phi_e(x_0) \hfill  \cr
\noalign{\vskip 4pt}
\hfill
= \int_{{\mathbb R}^n} \!
\Big\{
1_{E(x, \, 2r)}(y) 
\frac{1}{|x - y|^{\beta} 
\langle y \rangle^{\gamma}} 
-
1_{E(x_0, \, 2r)}(y) 
\frac{1}{|x_0 - y|^{\beta} 
\langle y \rangle^{\gamma}} 
\Big\} \, dy.
\qquad\qquad
\llap{(A.8)}\cr}
$$
Note that 
the inequality (\ref{eqn:9-e-4-1}) implies that
$$
1_{E(x, \, 2r)}(y) 
\frac{1}{\,|x - y|^{\beta} \,}
\le
\big( 
\frac{3}{\, 2 \,}
\big)^{\beta} \times
1_{E(x_0, \, 2r)}(y) 
\frac{1}{|x_0 - y|^{\beta} }
$$
whenever $|x - x_0| \le r$.
Hence, the integrand in (A.8) is 
bounded by
$$
\Big\{
\big( 
\frac{3}{\, 2 \,}
\big)^{\beta}  + 1  
\Big\} \times
1_{E(x_0, \, 2r)}(y) 
\frac{1}{|x_0 - y|^{\beta} 
\langle y \rangle^{\gamma}},     \eqno(\rm A.9)
$$
in absolute value, for all $x$ with
$|x - x_0| \le r$.
Since $\beta + \gamma >n$, 
by assumption of the lemma, we see that
the function in (A.9) belongs to 
$L^1({\mathbb R}^n)$. Therefore, the
Lebesgue dominated convergence theorem is
applicable to the right hand side of
(A.8) and shows that
$$
 \lim_{x \to x_0}
\big(\Phi_e(x) - \Phi_e(x_0) \big) = 0.
$$
Combining this with (A.6) and (A.7),
we deduce that
$$
 \limsup_{x \to x_0}
\big| \Phi(x) - \Phi(x_0) \big| 
\le 2 \varepsilon.     
$$
Since $\varepsilon$ was arbitrary,
this completes the proof of the step 1.
\hfill$\square$

\vspace{20pt}

To establish the desired inequalities, 
we make another decomposition of 
$\Phi (x)$:

\begin{eqnarray*}
\Phi (x)=\Phi_1(x)+\Phi_2(x)+\Phi_3(x),
\end{eqnarray*}
where
\begin{eqnarray*}
\Phi_1(x) &=& 
 \int_{|y| \le |x|/2}
\frac{1}{|x - y|^{\beta} 
\langle y \rangle^{\gamma}} \, dy,   \\
\noalign{\vskip 6pt}
\Phi_2(x) &=& 
 \int_{|x|/2 <|y| \le 2|x|}
\frac{1}{|x - y|^{\beta} 
\langle y \rangle^{\gamma}} \, dy,   \\
\noalign{\vskip 6pt}
\Phi_3(x) &=& 
 \int_{2|x| < |y|}
\frac{1}{|x - y|^{\beta} 
\langle y \rangle^{\gamma}} \, dy.  
\end{eqnarray*}
Since $\Phi (x)$ is bounded on
each compact subset of ${\mathbb R}^n$
by continuity of  $\Phi (x)$,
it is sufficient to get estimates of
$\Phi_i$'s for $|x|\ge 1$.

\vspace{20pt}

{\it Step 2. \ For $|x|\ge 1$, we have
\begin{eqnarray*}
|\Phi_1(x)|
\le  C_{\beta\gamma n}
\begin{cases}  
 |x|^{-(\beta + \gamma -n )}  
      &\text{  if $\; 0 < \gamma < n$,}   \\
{}    &\\
 |x|^{-\beta} 
   \log(1 + \langle x \rangle)
    &\text{  if $\; \gamma = n$, }  \\
{}    &\\
|x|^{-\beta} 
  &\text{  if $ \; \gamma >n$.  }
\end{cases}
\end{eqnarray*}
}

\vspace{5pt}

{\it Proof.}
Note that $|x - y|\ge |x|/2$ if
$|y| \le |x|/2$. 
This fact implies that
$$
\Phi_1(x) \le 2^{\beta} \, |x|^{-\beta}
\int_{|y| \le |x|/2}
\frac{1}{\, \langle y \rangle^{\gamma} \,} 
\, dy.
$$
If $0 < \gamma < n$, then
we get, using spherical polar coordinates,
\begin{eqnarray*}
\int_{|y| \le |x|/2}
\frac{1}{\, \langle y \rangle^{\gamma} \,} 
\, dy
\le \omega_n \!
 \int_0^{|x|/2} r^{-\gamma + n -1} \, dr
= 
\frac{\, \omega_n \, 2^{-\gamma + n}\,}{n -\gamma}
|x|^{-\gamma + n},
\end{eqnarray*}
where $\omega_n$ denotes the area 
of the unit sphere in ${\mathbb R}^n$.
Similarly, if $\gamma = n$,  we then have
\begin{eqnarray*}
\int_{|y| \le |x|/2}
\frac{1}{\, \langle y \rangle^{\gamma} \,} 
\, dy
&\le&
 \omega_n \!
 \int_0^{|x|/2} 
{\Big( 
\frac{ 1 + r }{\sqrt{2}} 
\Big)}^{\!\!\! -\gamma} \, r^{n -1} \, dr   \\
\noalign{\vskip 4pt}
&\le&
\omega_n \, 2^{\gamma/2} 
\int_0^{|x|/2} (1+r)^{-1} \, dr  \\
\noalign{\vskip 4pt}
&\le& 
\omega_n \, 2^{\gamma/2}
\log \! \big( 1 + \frac{\,|x|\,}{2} \big),  
\end{eqnarray*}
where we have used the 
inequality $\langle y \rangle \ge 
(1+ |y|)/ \sqrt{2}$.
If $\gamma > n$,  we evidently have
\begin{eqnarray*}
\int_{|y| \le |x|/2}
\frac{1}{\, \langle y \rangle^{\gamma} \,} 
\, dy
\le
\int_{{\mathbb R}^n}
\frac{1}{\, \langle y \rangle^{\gamma} \,} 
\, dy < + \infty.
\end{eqnarray*}
Summing up, we conclude that
the desired inequalities for $\Phi_1(x)$ hold.
\hfill$\square$

\vspace{20pt}

{\it Step 3. \ \ For $|x|\ge 1$, we have
\begin{eqnarray*}
|\Phi_2(x)|
\le  C_{\beta\gamma n} \,
 |x|^{-(\beta + \gamma -n )}.  
\end{eqnarray*}
}

\vspace{5pt}

{\it Proof.} Let $B^*(x)$ and $E^*(x)$ be 
the sets defined by
\begin{eqnarray*}
B^*(x) &:=& \{  \, y \in {\mathbb R}^n \; | \,\;
  | x - y | < \frac{\, |x| \,}{2} \; \}, \\
E^*(x) &:=& \{  \, y \in {\mathbb R}^n \; | \,\;
  \frac{\, |x| \,}{2} < |y| \le 2|x|,
 \, | x - y | \ge \frac{\, |x| \,}{2} \;
\} .
\end{eqnarray*}
Then we have
\begin{eqnarray*}
\Phi_2(x)= \int_{B^*(x)} \frac{1}{|x - y|^{\beta} 
\langle y \rangle^{\gamma}} \, dy
 + 
\int_{E^*(x)}  \frac{1}{|x - y|^{\beta} 
\langle y \rangle^{\gamma}} \, dy.
\end{eqnarray*}
Since $B^*(x)$ is a subset of the annulus
$\{ \, y \; | \,\;
 |x|/2 < |y| \le 2|x| \;
\}$,
it follows that
\begin{eqnarray*}
\langle y \rangle \ge 
\frac{1}{\, \sqrt{2}\,} (1 + |y| ) 
\ge  \frac{1}{\, 2\sqrt{2}\,} |x|
\qquad ( \, \forall y \in B^*(x) \, ),
\end{eqnarray*}
which gives
\begin{eqnarray*}
\int_{B^*(x)} \frac{1}{|x - y|^{\beta} 
\langle y \rangle^{\gamma}} \, dy
&\le& 2^{3\gamma /2} \, |x|^{-\gamma}
\int_{|x-y| < |x|/2} 
\frac{1}{|x - y|^{\beta} } \, dy    \\
\noalign{\vskip 4pt}
&=& 
\frac{\,2^{3\gamma /2 + \beta -n} \,\omega_n \,}
{n - \beta} \, |x|^{-\gamma -\beta + n}.
\end{eqnarray*}
If $y \in E^*(x)$, then
$|x-y|\ge |y|/4$, therefore
\begin{eqnarray*}
\int_{E^*(x)} \frac{1}{|x - y|^{\beta} 
\langle y \rangle^{\gamma}} \, dy   
&\le& 4^{\beta}
\int_{|x|/2 < |y| \le 2|x|} 
\frac{1}{|y|^{\beta + \gamma} } \, dy    \\
\noalign{\vskip 4pt}
&=& 
\frac{\,4^{\beta} \,\omega_n \,}
{\beta + \gamma -n} \, 
(2^{\beta +\gamma - n}
- 2^{-\beta -\gamma + n} ) \,
|x|^{ -\beta -\gamma + n}.
\end{eqnarray*}
Summing up, we obtain the desired
inequality for $\Phi_2(x)$.
\hfill$\square$

\vspace{20pt}

{\it Step 4. \  For $|x|\ge 1$, $\Phi_3(x)$ 
satisfies the same inequality as 
$\Phi_2(x)${\rm:}
\begin{eqnarray*}
|\Phi_3(x)|
\le  C_{\beta\gamma n} \,
 |x|^{-(\beta + \gamma -n )}. 
\end{eqnarray*}
}

\vspace{5pt}

{\it Proof.} 
If $2|x| < |y|$, then it follows that
\begin{eqnarray*}
|x -y| \ge |y| - |x| \ge 
  \frac{\, |y| \,}{2}.
\end{eqnarray*}
Hence we have
\begin{eqnarray*}
\Phi_3(x)   \le   2^{\beta}
\int_{2|x| < |y|} 
\frac{1}{\,|y|^{\beta +\gamma}\,} \, dy    
= 
\frac{\, 2^{-\gamma + n}\,\omega_n \,}
{\beta + \gamma -n}
\,|x|^{-\beta -\gamma + n}.
\end{eqnarray*}
This completes the proof.
\hfill$\square$

\vspace{15pt}

It is evident that Lemma A.1 follows
from the steps 1--4.


\end{document}